\documentclass[11pt]{amsart}
\usepackage{amscd,amsmath,amsthm,amssymb}
\usepackage{mathtools}
\usepackage{mathrsfs}
\usepackage[all,cmtip]{xy}
\usepackage{xcolor}
\usepackage{enumitem} 
\usepackage{array}
\usepackage{caption}
\usepackage{mathabx}
\usepackage{csquotes}
\usepackage{url}
\usepackage{xparse}
\usepackage{circuitikz}
\ctikzset{bipoles/resistor/height=0.15}
\ctikzset{bipoles/resistor/width=0.4}

\definecolor{cadmiumgreen}{rgb}{0.0, 0.42, 0.24}
\usepackage[
colorlinks, citecolor=cadmiumgreen,
pdfauthor={Omid Amini, Noema Nicolussi}, 
pdftitle={Moduli of hybrid curves I},
pdfstartview ={FitV},
]{hyperref}

\usepackage[
alphabetic,
msc-links,
nobysame,
lite,
]{amsrefs} 

\usepackage{tikz, float} 
\usetikzlibrary {positioning}
\usetikzlibrary{matrix, arrows}
\usetikzlibrary{calc,decorations.markings}
\usetikzlibrary{shapes,snakes}
\usepackage{a4wide}

\usepackage[T1]{fontenc}
\usepackage[utf8]{inputenc}
\usepackage[english]{babel}
\usepackage[leqno]{amsmath}
\usepackage{tabularx}
\usepackage{url}
\usepackage{hyperref}
\usepackage{tikz-cd}
\usepackage{csquotes}
\usepackage{mathtools}
\usepackage{manfnt}
\usepackage{mathrsfs}

\DeclareRobustCommand{\SkipTocEntry}[5]{}

\newtheorem{thm}{Theorem}[section]

\newtheorem{lem}[thm]{Lemma}

\newtheorem{prop}[thm]{Proposition}

\newtheorem{question}[thm]{Question}

\theoremstyle{definition}
\numberwithin{equation}{section}

\newcommand{\R}{\mathbb{R}}

\DeclareMathOperator{\rk}{rank} 

\renewcommand{\H}{\mathcal H}

\newcommand{\abs}[1]{\lvert #1\rvert}

\numberwithin{equation}{section}

\tikzstyle{Cwhite}=[scale = .8,circle, fill = white, minimum size=3mm] 
\tikzstyle{Cgray}=[scale = .4,circle, fill = gray, minimum size=3mm] 
\tikzstyle{Cblack2}=[scale = .4,circle, fill = black, minimum size=5mm] 
\tikzstyle{Cblack}=[scale = .7,circle, fill = black, minimum size=3mm]
\tikzstyle{C0}=[scale = .9,circle, fill = black!0, inner sep = 0pt, minimum size=3mm]
\tikzstyle{C1}=[scale = .7,circle, fill = black!0, inner sep = 0pt, minimum size=3mm]
\tikzstyle{Cred}=[scale = .4,circle, fill = red, minimum size=3mm]

\newcommand{\Spec}{{\rm Spec \,}}

\newcommand{\Acal}{{\mathcal A}}

\newcommand{\sE}{{\mathcal E}}

\newcommand{\sH}{{\mathcal H}}

\newcommand{\Ocal}{{\mathcal O}}

\newcommand{\C}{{\mathbb C}}

\renewcommand{\H}{{\mathbb H}}

\newcommand{\N}{{\mathbb N}}

\newcommand{\Z}{{\mathbb Z}}

\newcommand{\gr}{\mathrm{gr}}

\newcommand{\p}{\mathbf{p}}

\DeclareMathOperator{\Ext}{Ext}

\DeclareMathOperator{\id}{Id}
\DeclareMathOperator{\usExt}{\underline{\sE\!\mathit{xt}}}
\DeclareMathOperator{\usHom}{\underline{\sH\!\mathit{om}}}
\renewcommand{\Im}{\operatorname{Im}}

\makeatletter
\newsavebox\myboxA
\newsavebox\myboxB
\newlength\mylenA

\newcommand*\overbar[2][0.75]{%
	\sbox{\myboxA}{$\m@th#2$}%
	\setbox\myboxB\null
	\ht\myboxB=\ht\myboxA%
	\dp\myboxB=\dp\myboxA%
	\wd\myboxB=#1\wd\myboxA
	\sbox\myboxB{$\m@th\overline{\copy\myboxB}$}
	\setlength\mylenA{\the\wd\myboxA}
	\addtolength\mylenA{-\the\wd\myboxB}%
	\ifdim\wd\myboxB<\wd\myboxA%
	\rlap{\hskip 1\mylenA\usebox\myboxB}{\usebox\myboxA}%
	\else
	\hskip -0.5\mylenA\rlap{\usebox\myboxA}{\hskip 0.5\mylenA\usebox\myboxB}%
	\fi}
\makeatother

\newcommand{\comp}[1]{\overbar[.5]{#1}} 

\newcommand{\mghyb}[1]{\mg_{{\hspace{-.04cm}#1}}^\hyb}   
 
\newcommand{\mgbarhybr}[2]{\mgbar_{{\hspace{-.09cm}#1}}^{^{\scaleto{\Hyb(#2)}{4.8pt}}}}

\newcommand{\mg}{\mathscr M} 
\newcommand{\mggbar}[1]{\comp{\mathscr M}_{{\hspace{-.04cm}#1}}} 
\newcommand{\mgbar}{\comp{\mathscr M}} 

\usepackage{scalerel}

\newcommand{\Hyb}{{\mathrm{hyb}}}         
\newcommand{\hyb}{{^{\scaleto{\Hyb}{4.4pt}}}} 
\newcommand{\hybr}[1]{{^{\scaleto{\Hyb(#1)}{4.8pt}}}}  

\newcommand{\Trop}{\mathrm{trop}} 
\newcommand{\trop}{{^{\scaleto{\Trop}{4pt}}}}   

\newcommand{\rsf}{\mathcal S} 

\newcommand{\Zh}{{_{\scaleto{\mbox{Zh}}{4.5pt}}}} 

\newcommand{\Ar}{{_{\scaleto{\mbox{Ar}}{4.5pt}}}} 

\newcommand{\cT}{\mathcal{T}}

\newcommand{\base}{B} 

\newcommand{\bp}{\flat} 

\newcommand{\genusfunction}{\mathfrak g} 
\newcommand{\graphgenus}{{h}} 
\newcommand{\mc}{{\mathcal{M}\mathscr{C}}} 
\newcommand{\genus}{h} 

\newcommand{\subface}{\leq}

\newcommand{\basisa}{\mathscr A}
\newcommand{\basisb}{\mathscr B}

\newcommand{\Log}{\mathrm{Log}}

\newcommand{\mgr}{\mathscr G} 

\renewcommand{\lg}{\mathscr G} 
\newcommand{\filter}{\mathscr F} 
\newcommand{\filt}{\mathrm{F}}  
\newcommand{\dfilter}{\mathscr E} 
\newcommand{\Can}{\mathrm{can}}
\newcommand{\can}{{{\scaleto{\Can}{2.5pt}}}}

\newcommand{\hcurve}{\mathscr C}  
\newcommand{\hcurvef}{\mathscr C^\hyb} 

\newcommand{\grm}{\mathrm{gr}} 
\newcommand{\proj}{{\scaleto{\mathscr K}{5.6pt}}} 
\newcommand{\cont}[1]{{{\proj_{\hspace{-.04cm}{\scaleto{#1}{5.4pt}}}}}}
\newcommand{\continv}[1]{{\proj_{\hspace{-.04cm}{\scaleto{#1}{4pt}}}}^{\scaleto{\hspace{-.1cm}-1}{4.9pt}}}

\newcommand{\pr}{\mathfrak p} 

\newcommand{\forget}{\mathfrak q}

\newcommand{\unicurve}{\mathscr C}

\newcommand{\thy}{\mathbf{t}} 
\newcommand{\Pty}{\mathbf{P}} 

\newcommand{\s}{\mathbf{s}} 

\newcommand{\marking}{{\mathfrak m}}
\newcommand{\countmarking}{{\mathfrak n}}
\newcommand{\contract}[2]{#1\big/#2}

\newcommand{\aut}{\mathrm{Aut}}

\newcommand{\Ical}{\mathcal I}

\newcommand{\rest}[1]{\raisebox{-1pt}{$\vert$}_{#1}}

\newcommand{\mA}{ {\mathcal{A}} } 
\newcommand{\mB}{ {\mathcal{B}} } 

\theoremstyle{definition}
\newtheorem{defii}[thm]{Definition}
\newenvironment{defi}
{\pushQED{\qed}\defii}
{\popQED\enddefii}
\newtheorem{remm}[thm]{Remark}
\newenvironment{remark}
{\pushQED{\qed}\remm}
{\popQED\endremm}
\newtheorem{exx}[thm]{Example}
\newenvironment{example}
{\pushQED{\qed}\exx}
{\popQED\endexx}

\numberwithin{equation}{section}

\newcommand{\inn}{\mathring}          

\newcommand{\Bg}{B_{\scaleto{G}{4pt}}}

\newcommand{\transpose}{{\scaleto{\mathrm{T}}{4.5pt}}}

\newcommand{\Pifs}{\Pi}
\newcommand{\st}{\bigm|} 

\renewcommand{\setminus}{\smallsetminus}

\makeatletter
\let\@oldinfty\infty
\newcommand{\@sminfty}{{\scaleto{\@oldinfty}{2.8pt}}} 
\renewcommand{\infty}{{\mathchoice%
		{\displaystyle{\@oldinfty}}%
		{\textstyle{\@oldinfty}}%
		{\scriptstyle{\@sminfty}}%
		{\scriptscriptstyle{\@sminfty}}}
}
\makeatother

\newcommand{\ndim}{{\scaleto{\mathrm{N}}{6pt}}}
\makeatletter
\let\@oldndim\ndim
\newcommand{\@smndim}{{\scaleto{\@oldndim}{4.4pt}}} 
\renewcommand{\ndim}{{\mathchoice%
		{\displaystyle{\@oldndim}}%
		{\textstyle{\@oldndim}}%
		{\scriptstyle{\@smndim}}%
		{\scriptscriptstyle{\@smndim}}}
}
\makeatother

\NewDocumentCommand{\ssub}{O{0pt} O{.8} m t! e{_^}}{
	#3%
	\IfValueT{#5}{
		\IfBooleanTF{#4}{\sb{\hspace{#1}\scaleobj{#2}{#5}}}{\sb{#5}}
	}
	\IfValueT{#6}{\sp{#6}}
}

\ExplSyntaxOn
\NewDocumentCommand{\tossub}{o o m}{
	\expandafter\let\csname old\cs_to_str:N #3\endcsname#3
	\renewcommand#3%
	{\ssub[#1][#2]{\csname old\cs_to_str:N #3\endcsname}}
}
\ExplSyntaxOff

\newcommand{\ssR}{\ssub{\mathbb R}!}
\newcommand{\ssZ}{\ssub{\mathbb Z}!}
\tossub\sigma
\tossub\mg
\tossub\hcurve
\tossub\omega

\newcommand{\imi}{\mathrm{i}}

\begin{document}
\title[Moduli of hybrid curves I: Variations of canonical measures]{Espace de modules des courbes hybrides I : \\
Variations des mesures canoniques \\ \vspace{.4cm} Moduli of hybrid curves I: \\ Variations of canonical measures}

\author{Omid Amini}
\address{CNRS - Centre de math\'ematiques Laurent Schwartz, \'Ecole Polytechnique}
\email{\href{omid.amini@polytechnique.edu}{omid.amini@polytechnique.edu}}

\author{Noema Nicolussi}
\address{Institute of Analysis and Number Theory, Faculty of Mathematics, Physics and Geodesy, Graz University of Technology}
\email{\href{nicolussi@math.tugraz.at}{nicolussi@math.tugraz.at}}
\thanks{\emph{To Maryam Mirzakhani,} \vspace{.1cm} {with admiration for her sense of beauty in mathematics.}}

\begin{abstract}
Le présent article est le premier d'une série de travaux consacrés à l'étude de la géométrie asymptotique des surfaces de Riemann et de leurs espaces de modules. 

Nous introduisons l'espace de modules des courbes hybrides comme une nouvelle compactification de l'espace de modules des courbes, raffinant celle construite par Deligne et Mumford. Il s'agit de l'espace de modules pour les objets géométriques multi-échelles qui mélangent la géométrie complexe et la géométrie tropicale et non archimédienne de  rang supérieur, reflétant à la fois des caractéristiques discrètes et continues.

Nous définissons des mesures canoniques sur les courbes hybrides qui généralisent les mesures d'Arakelov-Bergman sur les surfaces de Riemann et les mesures de Zhang sur les graphes métriques.

Nous montrons ensuite que la famille universelle de courbes hybrides munies de leurs mesures canoniques est une famille continue d'espaces mesurables au-dessus de cet espace de modules hybride. Ce résultat fournit un lien précis entre la mesure de Zhang non archimédienne et les variations des mesures d'Arakelov--Bergman dans les familles de surfaces de Riemann, répondant ainsi à une question ouverte depuis les travaux pionniers de Zhang sur l'accouplement admissible dans les années quatre-vingt-dix.

\vspace{.2cm}

The present paper is the first in a series devoted to the study of asymptotic geometry of Riemann surfaces and their moduli spaces.

We introduce the moduli space of hybrid curves as a new compactification of the moduli space of curves, refining the one obtained by Deligne and Mumford. This is the moduli space for multiscale geometric objects which mix complex and higher rank tropical and non-Archimedean geometries, reflecting both discrete and continuous features.

We define canonical measures on hybrid curves which combine and generalize Arakelov--Bergman measures on Riemann surfaces and Zhang measures on metric graphs.

 We then show that the universal family of canonically measured hybrid curves over this moduli space varies continuously. This provides a precise link between the non-Archimedean Zhang measure and  variations of Arakelov--Bergman measures in families of Riemann surfaces, answering a question which has been open since the pioneering work of Zhang on admissible pairing in the nineties.

 \end{abstract}

\maketitle

\setcounter{tocdepth}{1}

\tableofcontents

\section{Introduction}
 In this paper we study \emph{canonical measures} on Riemann surfaces  and their \emph{tropical and hybrid limits}.  
 
  By canonical measure on a Riemann surface we mean the Arakelov--Bergman measure $\mu_\Ar$ defined in terms of holomorphic one-forms. For a compact Riemann surface $S$ of positive genus $g$, this is the positive density measure of total mass $g$ on $S$ given by
  \[\mu_{\Ar} \coloneqq  \frac{\imi}{2}\sum_{j=1}^g \omega_j \wedge \bar{\omega}_j,\]
  where $\imi=\sqrt{-1}$ and $\omega_1, \dots,\omega_g$ form an orthonormal basis for the space of holomorphic one-forms on $S$, with respect to the  Hermitian inner product
  \[ \langle \eta_1, \eta_2\rangle \coloneqq \frac \imi{2} \int_S \eta_1 \wedge \bar{\eta}_2\]
  for pairs of holomorphic one-forms $\eta_1, \eta_2$ on $S$. 

Let $\mg_g$ be the moduli space of curves of genus g and let $\mggbar{g}$ be its Deligne–Mumford compactification, consisting of stable curves of genus $g$. The question we address in this paper is the following:

\begin{question}\label{question:main}
Consider a sequence of smooth compact Riemann surfaces $S_j$ of genus $g$ such that the corresponding points $s_j$ of $\mg_g$ converge to a point in $\mggbar{g}$. Let $\mu^\can_j$ be the canonical measure of $S_j$. What is the asymptotic behavior of the measures $\mu^\can_j$?
\end{question}
Since the measures $\mu^\can_j$ live on distinct surfaces, the first problem to handle consists in giving a precise mathematical meaning to the question. Moreover, once this has been taken care of, the answer to the question appears to be sensitive to the \emph{speed} and \emph{direction} of the convergence of the sequence of points $(s_j)$. Formalizing these points leads to the definition of new geometric objects called \emph{hybrid curves} and their \emph{moduli spaces}.

Although discovered in context with canonical measures, hybrid curves and their moduli spaces provide a framework in which several other questions can be answered as well, see Section~\ref{ss:FurtherWork} for a brief outlook on further developments.

\subsection{Hybrid curves}
Let $s_\infty$ be the limit of $(s_j)$ and let $S_\infty$ be the stable Riemann surface corresponding to $s_\infty$. Denote by $G=(V,E)$ the dual graph of $S_\infty$: its vertices are in bijection with the irreducible components of $S_\infty$ and its edges are in bijection with the nodes (singular points) of $S_\infty$. If an edge $e$ corresponds to a node $p$, then its endpoints are given by the two components  of $S_\infty$ (possibly identical, if $e$ is a loop) containing $p$. The genera of components of $S_\infty$ define a \emph{genus function} $\genusfunction \colon V\to \N$. The pair $(G, \genusfunction)$ is called the \emph{stable dual graph} of $S_\infty$.

With the help of the stable dual graph, it will be possible to capture additional information on the speed and direction of convergence of the sequence.

We can view the speed of convergence as a way to distinguish the \emph{relative order of appearance of the singular points} in the limit, when the sequence $S_1, S_2, \dots$ approaches the stable Riemann surface $S_\infty$. This leads to an ordered partition $\pi = (\pi_1, \dots, \pi_r)$ of the edge set $E$, where $\pi_1$ are the edges corresponding to the \emph{fastest appearing nodes}, $\pi_2$ are the ones which are \emph{fastest among the remaining nodes}, and so on. We will call the ordered partition $\pi$ a \emph{layering} of the graph, the elements of $\pi$ the \emph{layers}, and the graph $G$ endowed with the layering a \emph{layered graph}.

Capturing the direction corresponds to choosing edge length functions $\ell_1 \colon \pi_1 \to \R_+, \dots,$ $\ell_r \colon \pi_r \to \R_+$ on the layers, each of them well-defined up to multiplication by a positive scalar. Here and everywhere else, $\R_+$ denotes the set of strictly positive real numbers.

\noindent We have now arrived at the data underlying the definition of a hybrid curve: 
\begin{itemize}[leftmargin = 2em]
\item a stable curve $S$ with stable dual graph $G=(V,E,\genusfunction)$ and
\item an ordered partition $\pi=(\pi_1, \dots, \pi_r)$ of $E$, for some $r\in \mathbb N$, and
\item an edge length function $\ell \colon E \to \R_+$. 

\end{itemize}

Given the triple $(S, \ell, \pi)$, we define the associated \emph{layered metrized complex} as the metrized complex $\mc$ from~\cite{AB15}, obtained as the \emph{metric realization} of $(S, \ell)$, and enriched with the data of the layering $\pi$.

We recall that $\mc$ is obtained by taking first the normalization $\widetilde S$ of $S$, which is by definition the disjoint union of the (normalization of the) irreducible components of $S$, then taking an interval $\Ical_e$ of length $\ell_e$ for each edge $e \in E$ associated to a node $p_e$ in $S$, and finally gluing the two extremities of $\Ical_e$ to the two points in the normalization $\tilde S$ corresponding to the node $p_e$ in $S$ (see Figure~\ref{fig:HybridCurveIntro}). 

\begin{figure}[!t]
\centering
    \scalebox{.28}{\input{kite2.pspdftex}}
\caption{An example of a hybrid curve of rank three. The graph of the underlying stable Riemann surface has five vertices and seven edges. Its edges are partitioned into three sets $\pi_1,\pi_2,$ and $\pi_3$.}
\label{fig:HybridCurveIntro}
\end{figure}

We then define the following \emph{conformal equivalence relation} (\emph{at infinity}) on layered metrized complexes. Two layered metrized complexes $\mc$ and $\mc'$ are conformally equivalent if there exists an isomorphism of the semistable curves underlying $\mc$ and $\mc'$ so that under this isomorphism, they have the same underlying graph $G=(V, E)$ (this will be automatic), the same ordered partition $\pi =(\pi_1, \dots, \pi_r)$ on $E$, and for each layer $\pi_j$, $j=1, \dots, r$ of $\pi$, there is a positive number $\lambda_j >0$ such that
\[\ell\rest{\pi_j} = \lambda_j \, \ell'\rest{\pi_j}
\]
for the respective edge length functions $\ell, \ell' \colon  E \to \R_{+}$.

A \emph{hybrid curve} $\hcurve^\hyb$ is a conformal equivalence class of layered metrized complexes.  In other words, the edge length function $\ell$ is well defined  only up to multiplying its restrictions $\ell\rest{\pi_j}$ to layers by positive real numbers. Imposing $\sum_{e\in \pi_j} \ell_e=1$ for each $j=1, \dots, r$ leads to a unique choice of $\ell$. The integer $r$ is called the \emph{rank} of $\hcurve^\hyb$. 

Note that each metrized complex $\mc$ gives rise to a hybrid curve $\hcurve^\hyb$, namely the conformal equivalence class of $\mc$ endowed with the \emph{trivial layering} (the trivial ordered partition consisting of one element $\pi = (E)$). In particular, hybrid curves provide an enrichment of the category of metrized curve complexes (with normalized edges lengths).

\subsection{Moduli space $\mghyb{g}$ of hybrid curves of genus $g$} 
In order to answer Question~\ref{question:main}, we will construct the \emph{moduli space $\mg_g^\hyb$ of hybrid curves of genus $g$} that provides a new compactification of $\mg_g$ naturally lying over the Deligne--Mumford compactification $\mggbar{g}$. One may view $\mg_g^\hyb$ as a \emph{higher rank tropical refinement} of $\mggbar{g}$, that allows to study problems involving both Archimedean and tropical non-Archimedean aspects of Riemann surfaces and their families. The precise meaning of this is subject to our work~\cite{AN-hybrid-green} and its forthcoming sequels (see Section~\ref{ss:FurtherWork}).

The points in $\mg_g^\hyb$ correspond to \emph{hybrid limits of Riemann surfaces of genus $g$}. The topology reflects aspects which are  reminiscent of both analytic and Zariski topologies in $\mggbar{g}$.  Roughly speaking, our construction replaces each point of $\mggbar{g}$ representing a non-smooth stable Riemann surface $S$ by infinitely many hybrid points representing the hybrid curves with $S$ as their underlying stable Riemann surface.

We outline the construction here and refer to Section~\ref{sec:hybrid_moduli_space} for the details.

\subsubsection{Hybrid spaces of higher rank} The definition of the hybrid moduli space $\mg_g^\hyb$ is based on a general construction of hybrid spaces of higher rank. To a complex manifold $B$ and a simple normal crossing divisor $D$, we associate a hybrid space $B^{\hyb}$ by enriching the points $x$ on the divisor $D$ with additional simplicial coordinates. In the following, we briefly summarize this procedure (see Section~\ref{sec:hybrid_spaces} for more information).

Let $B$ be an $\ndim$-dimensional complex manifold, for a positive integer $\ndim$, and consider a simple normal crossing divisor  $D = \bigcup_{e \in E} D_e$ in $B$. That is, we require that $(D_e)_{e \in E}$ is a finite family of smooth, connected and closed submanifolds of codimension one in $B$ such that for any subset $F \subseteq E$, the stratum 
\[
 D_F \coloneqq  \bigcap_{e \in F} D_e
\]
is either empty or a smooth submanifold of codimension $|F|$ with only finitely many connected components. We define the \emph{inner stratum} of $F$ as the subset
\[
{\inn D_F} \coloneqq  D_F \setminus \bigcup_{e \notin F} D_e. 
\]
Note that for $F = \varnothing$, we recover the \emph{open part} ${\inn D_\varnothing} = B \setminus D \eqqcolon B^\ast$. Moreover, the inner strata provide a partition of $B$,
\[
B = \bigsqcup_{\substack{ F \subseteq E }} {\inn D_F}.
\]

Let $\Pifs(F)$ be the set of \emph{ordered partitions} of a subset $F \subseteq E$. Each element $\pi \in \Pifs(F)$ is an ordered sequence $\pi = (\pi_1, \pi_2, \dots, \pi_r)$, for some $r\in \N$, consisting of non-empty, pairwise disjoint subsets $\pi_i \subseteq F$  with $\bigsqcup_{i=1}^r \pi_i = F$. 

For an ordered partition $\pi = (\pi_i)_{i=1}^r$ in $\Pifs(F)$, $F \subseteq E$, the \emph{hybrid stratum} $D_\pi^\hyb$ is defined as
\[
{D}^{\hyb}_\pi \coloneqq   \inn D_F \times \inn\sigma_{\pi_1} \times \dots \times \inn\sigma_{\pi_r}, \]
where for a finite set $A$, we denote by $\sigma!_A$ the standard simplex in $\R^A$, with vertices the standard basis of $\R^A$, and by $\ssub{\inn\sigma}!_A$ the relative interior of $\sigma!_A$. For $\pi = \pi_\varnothing$, that is, the empty ordered partition of  the empty set $F= \varnothing \subseteq E$, we set ${D}^{\hyb}_{\pi_\varnothing} \coloneqq  \inn D_{\varnothing} = B^\ast$. The \emph{hybrid space} $B^{\hyb}$ is the disjoint union
\[
B^{\hyb} \coloneqq \bigsqcup_{\substack{  F \subseteq E }} \bigsqcup_{\substack{ \pi \in \Pifs(F)}} {D}^{\hyb}_\pi = B^\ast \sqcup \bigsqcup_{\substack{ \varnothing \subsetneq F \subseteq E }} \bigsqcup_{\substack{ \pi \in \Pifs(F)}} {D}^{\hyb}_\pi
\]
equipped with the \emph{hybrid topology} (see Section~\ref{sec:hybrid_spaces}).

\subsubsection{Definition of the hybrid moduli space} In the definition of the moduli space $\mg_g^\hyb$, we apply the above to the versal deformation spaces of stable Riemann surfaces $S$. In this context $B = \Delta^{\ndim}$ is a polydisc and the divisor is a union of coordinate hyperplanes $D =  \bigcup_{e \in E}  \{z_e = 0\}$, where $z_e$, $e \in E$, are coordinates parametrizing the appearance of the nodes of $S$. The polydisk $B$ provides an étale chart around the point $s$ that represents $S$ in the Deligne-Mumford stack $\mggbar{g}$. We call the resulting hybrid space $B^{\hyb}$ a \emph{hybrid \'etale chart}. The topology on $B^{\hyb}$ formalizes the aforementioned layered appearance of the nodes. The hybrid moduli space $\mg_g^\hyb$ is then obtained by glueing the hybrid \'etale charts.

   The hybrid moduli space $\mg_g^\hyb$ naturally comes with the \emph{universal family of hybrid curves of genus $g$}, defined more precisely on hybrid \'etale charts $B^\hyb$,
 \[\hcurvef_g \to \mg_g^\hyb,  \qquad \qquad \rsf^\hyb \to B^{\hyb}.
 \]
The fiber $\rsf_{\thy}^\hyb$ over any point $\thy \in B^{\hyb}$ corresponds to the hybrid curve represented by $\thy$. The universal family in turn inherits a natural hybrid topology which makes the projection maps $\hcurvef_g \to \mg_g^\hyb$ and $\rsf^\hyb \to B^{\hyb}$ continuous. This leads to commutative diagrams of continuous maps
 \[
 \begin{tikzcd}
\hcurvef_g \dar \rar &  \mghyb{g}   \dar \\  
\hcurve_g \rar& \mggbar{g}
 \end{tikzcd}, \qquad  \qquad 
 \begin{tikzcd}
\rsf^\hyb  \dar \rar &  B^{\hyb}   \dar \\  
\rsf \rar& B
 \end{tikzcd}
 \] 
with $\hcurve_g $ refering to the universal family of stable Riemann surfaces of genus $g$ over $\mggbar{g}$, defined over an étale chart $B \to \mggbar{g}$ by the map $\rsf \to B$ with fiber $\rsf_t$ at $t\in B$ given by the corresponding stable Riemann surface to $t$.  The constructions are thus compatible with the classical picture.
We refer to Section~\ref{sec:hybrid_moduli} for more details.

\subsection{Hybrid canonical measures} In order to identify all the possible limits in Question~\ref{question:main}, we will define a \emph{canonical measure} on a hybrid curve. This measure combines and generalizes both Archimedean canonical measures on Riemann surfaces and non-Archimedean measures on metric graphs. The latter arise naturally in non-Archimedean Arakelov geometry, when working with analytification of curves over non-Archimedean fields and their skeleta.

The canonical measure on a metric graph is the one introduced by Zhang in his work on admissible pairing~\cite{Zhang}. Zhang's measure $\mu_\Zh$ is a weighted combination of Lebesgue measures on edges, with weights carrying interesting information on spanning trees (see Section~\ref{sec:measures} for several equivalent definitions). Moreover, when viewed on the analytified curve, it carries essential arithmetic geometric information~\cite{Zhang, Zhang2, Cinkir, Ami-W, deJong}.

Consider a hybrid curve $\hcurve^\hyb$, defined as the metric realization of a stable curve $S$ with dual graph $G =(V,E)$, endowed with a layering $\pi=(\pi_1, \dots, \pi_r)$, and a layerwise normalized edge length function $\ell$. We define the \emph{canonical measure} $\mu^\can$ on $\hcurve^\hyb$ as the sum 
\[\mu^\can \coloneqq  \mu_\Ar + \mu^1_\Zh + \dots + \mu^r_\Zh\]
where

\begin{enumerate}[leftmargin = 2em]
\item the measure $\mu_\Ar$, \emph{the Archimedean part of the canonical measure}, restricts to the Arakelov--Bergman measure on each positive-genus component of $S$, viewed inside $\hcurve^\hyb$, and vanishes elsewhere;

\item the measure $\mu^j_\Zh$,  \emph{the $j^{\mathrm{th}}$ graded non-Archimedean part of the canonical measure}, is supported on the intervals $\Ical_e$, for edges $e\in \pi_j$, in $\hcurve^\hyb$. On each such interval $\Ical_e$, it has the same restriction as  the Zhang measure of the metric graph obtained by removing all edges which appear in lower layers $\pi_1, \dots, \pi_{j-1}$ and then contracting all edges which appear in upper layers $\pi_{j+1}, \dots, \pi_r$. Since these graphs appear as result of deletion and subsequent contraction, we call them \emph{graded minors} of $\hcurve^\hyb$, appealing to Robertson--Seymour theory of graph minors. (Any minor of $G$ appears in this way as a graded minor associated to some ordered partition of $E$.) 
\end{enumerate}

Examples of canonical measures are given in Figure~\ref{fig:canonical_measures}. Note that when $\hcurve^\hyb$ corresponds to a smooth curve,  then $\mu^\can =\mu_\Ar$. If $\hcurve^\hyb$ is a metrized curve complex $\mc$ with the trivial ordered partition $\pi=(E)$, then $\mu^\can = \mu_\Ar + \mu_\Zh$ is the Arakelov--Zhang measure on $\mc$. The latter underlies the equidistribution of Weierstrass points on non-Archimedean curves~\cite{Ami-W}. 

When the rank of hybrid curves is larger than $2$, this definition becomes very sensitive to the order of the layers, that is, a change in order (while keeping the same edge lengths) can drastically change the canonical measure.  
\subsection{Continuity of the canonical measures over the hybrid moduli space} The main theorem of this paper reads as follows.
 
\begin{thm}\label{thm:mainglobal-intro} For each hybrid étale chart $B^{\hyb}$, the family of canonically measured hybrid curves $(\rsf^\hyb,\mu^\can)$ forms a continuous family of measured spaces over $B^{\hyb}$. In particular, the universal family of canonically measured hybrid curves $(\hcurvef_g, \mu^\can)$ forms a continuous family of measured spaces over the hybrid moduli space $\mg_g^\hyb$. \end{thm}

Continuity in the above statement is defined in the distributional sense. That is, for any hybrid étale chart $B^\hyb \to \mg_g^\hyb$ with the corresponding family $\rsf^\hyb \to B^\hyb$, and any continuous function $f \colon \rsf^\hyb \to \R$, the induced function $F \colon B^\hyb\to \R$ defined by \emph{integration along fibers}
\begin{align*}
&F(\thy) \coloneqq  \int_{\rsf^\hyb_{\thy}} f_{|_{\rsf^\hyb_{\thy}}} \, d\mu^\can_\thy, \qquad  \thy \in B^\hyb,
\end{align*}
is continuous on $B^\hyb$. Here, $\mu_\thy^\can$ denotes the canonical measure on the hybrid curve $\rsf^\hyb_\thy$, the fiber of $\rsf$ over the point $\thy\in B^\hyb$.

The second statement in Theorem~\ref{thm:mainglobal-intro} asserts that the universal family $\hcurvef_g \to \mg_g^\hyb$ verifies the analogous distributional continuity. Note that the fiber $\hcurvef_\thy$ over $\thy \in  \mg_g^\hyb$ is the topological quotient of the hybrid curve represented by $\thy$ by the action of its automorphism group. The canonical measure $\mu_\thy^\can$ on the fiber $\hcurvef_\thy$ is defined as the pushforward measure of the canonical measure on the corresponding hybrid curve by the quotient map. The claimed continuity follows directly from the first statement of the theorem, since $\hcurvef_g \to \mg_g^\hyb$ is locally a topological quotient of the hybrid family $\rsf^\hyb \to B^{\hyb}$ over a hybrid étale chart $B^{\hyb}$ (see Section~\ref{sec:hybrid_moduli}).

 Theorem~\ref{thm:mainglobal-intro} gives a complete answer to Question~\ref{question:main}. If the sequence $(s_j)$ converges in $\mghyb{g}$ to a point $\thy$ corresponding to the hybrid curve $\rsf^\hyb_\thy$, then the canonical measures $\mu_j^\can$ of $S_j$ converge  to the canonical measure of $\rsf^\hyb_\thy$. Since the hybrid moduli space $\mghyb{g}$ is a compact Hausdorff space, the converse is also true. In particular, Theorem~\ref{thm:mainglobal-intro} proves that the canonical measure cannot be extended continuously to the singular curves in the universal family over $\mggbar{g}$.

In Sections~\ref{sec:hybrid_spaces} and~\ref{sec:hybrid_moduli}, we will construct a tower of hybrid spaces  interpolating between $\mggbar{g}$ and $\mghyb{g}$. Instead of all hybrid curves of a fixed genus $g$, we restrict to those having {\em rank at most $k$}. This means that their layerings $\pi = (\pi_1, \dots, \pi_r)$ have at most $k$ parts, i.e., $r\leq k$. The $k$-th hybrid space $\mgbarhybr{g}{k}$ compactifies $\mg_g$ by adding these hybrid curves as its boundary part. However,  to ensure compactness, the notions are slightly relaxed to allow some vanishing edge lengths (see  Sections~\ref{sec:hybrid_spaces} and~\ref{sec:hybrid_moduli}).  Altogether, this leads to the following tower of hybrid compactifications
\begin{align}\label{eq:tower-moduli}
\mggbar{g} \longleftarrow \mgbarhybr{g}{1} \longleftarrow \mgbarhybr{g}{2} \longleftarrow \dots \longleftarrow \mgbarhybr{g}{\ndim}= \mg_{g}^\hyb
\end{align}
for $\ndim = 3g-3$.  We refer to $\mgbarhybr{g}{k}$ as \emph{the rank $k$ hybrid compactification of $\mg_g$} and call it the \emph{compactified moduli space of hybrid curves of rank {at most $k$}} (because of the possible presence of length zero edges). Then, $\mg_g^{\hybr{k}}$, the \emph{moduli space of hybrid curves of rank {at most $k$}}, is an open subset. The first hybrid space $\mgbarhybr{g}{1}$ can be interpreted as a compactification of $\mg_g$ by means of metrized complexes.

 None of the intermediate spaces in the tower verify the continuity property in Theorem~\ref{thm:mainglobal-intro}.

Theorem~\ref{thm:mainglobal-intro} implies similar continuity results for other families of Riemann surfaces. Consider a family of stable Riemann surfaces $\rsf \to X$, over a complex manifold  $X$, whose discriminant locus (i.e., the locus of points in $X$ whose fiber in the family is not smooth)  forms a simple normal crossing divisor  $D$ in $X$. Assume that $\rsf \to X$ comes from a toroidal map $f \colon X \to \mggbar{g}$, meaning that we have an isomorphism $\rsf \simeq f^*(\hcurve_g)$ locally on étale charts.  From the toroidal assumption we get a map $f^\hyb \colon X^\hyb \to \mghyb{g}$, defined locally on hybrid étale charts. 
The hybrid space $X^\hyb$ associated to $X$ and the divisor $D$ lies above $X$ and comes with the family of hybrid curves $\rsf^\hyb = f^{\hyb^{*}}(\hcurvef_g)$, which is obtained by removing the non-smooth fibers over $D$ and adding appropriate hybrid curves as fibers over $X^{\hyb} \setminus (X \setminus D)$. Initially, this is defined on hybrid étale charts, but since the original family is defined over $X$, the family of hybrid curves descends over $X^\hyb$ (see the discussion in Section~\ref{sec:universal_hybrid_curve}).

\begin{thm} \label{thm:families}The family of measured spaces $(\rsf^\hyb_\thy, \mu^\can_\thy)_{\thy \in X^\hyb}$ is continuous. \end{thm}

One specific case of Theorem~\ref{thm:families} concerns families varying over a one-dimensional base, i.e., when $\dim(X)=1$. In this case, the divisor $D$ consists of a finite set of points $t_0, \dots, t_n$ in $X$ and our hybrid construction simply replaces the stable curves $\rsf_{t}$, $t\in D$, with the corresponding metrized complex $\mc_{t}$. 
Let $\mgr_{t}$ be the underlying stable metric graph of $\mc_{t}$ and consider the projection map $\mc_t \to \mgr_t$. The pushout of the canonical measure on $\mc_t$ coincides with the Zhang measure on the stable metric graph $\mgr_t$. In this specific situation, we get the following theorem independently proved by Sanal Shivaprasad~\cite{Shiva} using  methods from non-Archimedean analysis.

\begin{thm}[Shivaprasad~\cite{Shiva}]\label{thm:main_onedimension} Let $X$ be a complex curve, and consider a one parameter family $\rsf \to X$ of stable complex curves of genus $g$ over $X$ with smooth generic fiber. For a point  $t_0 \in X$ with singular fiber $\rsf_{t_0}$, let $\mc_{t_0}$ be the corresponding metrized complex. Then, as $t$ tends to $t_0$, the measured spaces $(\rsf_t, \mu^\can_t)$ converge weakly to $(\mc_{t_0}, \mu^\can_{t_0})$. 
\end{thm}

In particular, in a one parameter family as above, the Zhang measure on the stable metric graph $\mgr_{t_0}$  may be viewed as the limit of the canonical measures $\mu_t^\can$ when $t$ tends to $t_0$. A special consequence of Theorem~\ref{thm:main_onedimension}, \emph{the total mass of the limit measure on a node of $\rsf_{t_0}$},  was obtained by R.\ de Jong~\cite{deJong}. Wentworth~\cite{Wentworth} treated the case where the stable Riemann surface has a single node, i.e., the dual graph has only one edge, and the base space $X$ has dimension $\dim(X) = 1$. More recently, the case where $\rsf_{t_0}$ is a curve of compact type (that is, the dual graph is a tree), or a rational curve with nodes (that is, the dual graph is a rose) was treated by Ng--Yeung~\cite{NY20}. 

As a comparison to Theorem~\ref{thm:mainglobal-intro}, we note that other natural families of measures associated to Riemann surfaces, namely hyperbolic and Narasimhan-Simha measures, behave completely differently. In both cases, the measures extend continuously to the boundary of the Deligne-Mumford compactification, see~\cite{Abi77} (and recent work~\cite{PS}) for the former, and~\cite{Shi20} for the latter. In the hyperbolic case,  the limit measure on a stable Riemann surface is the sum of hyperbolic measures on its punctured irreducible components, punctures being given by the nodes.

\subsection{Applications and outlook} \label{ss:FurtherWork}  The present paper forms the basis for a series of works where we study the asymptotic geometry of Riemann surfaces  and their moduli spaces.

In the sequel~\cite{AN-hybrid-green}, we introduce a notion of Laplace operator and formulate a Poisson equation on hybrid curves. We then obtain the asymptotics of solutions to Poisson equations on degenerating families of Riemann surfaces in terms of solutions to the hybrid Poisson equation.  Using Theorem~\ref{thm:mainglobal-intro}, we describe the behavior of the Arakelov--Bergman Green function on Riemann surfaces close to the boundary of their moduli space. This solves essentially completely a problem arising in Arakelov geometry and mathematical physics~\cite{HGB19, deJong, Faltings,Jorgenson90,Wentworth}. Moreover, \cite{AN-hybrid-green} is related to the recent work of Yang Li~\cite{YangLi} on SYZ conjecture and Monge-Amp\`ere equation. In this work, Yang Li reduces the SYZ conjecture in one parameter families to the existence of solutions (still conjectural) to some skeletonized Monge-Amp\`ere equation. He then obtains the dominating terms of the solution to the nearby complex Monge-Amp\`ere equation. 
On Riemann surfaces, the Monge-Amp\`ere equation corresponds to the Poisson equation. In~\cite{AN-hybrid-green}, we obtain a full asymptotic expansion for solutions, including higher order terms, for multiparameter families of degenerating Riemann surfaces.

In the third part of our series, we establish for hybrid curves analogs of fundamental theorems on smooth compact Riemann surfaces (e.g., Riemann--Roch, Abel--Jacobi, and Poincaré-Lelong  theorems). In particular, we associate Jacobian tori to hybrid curves and interpret them as multi-scale Gromov--Hausdorff limits of Jacobians of degenerating Riemann surfaces. For this purpose, we apply the framework of higher rank inner products developed in~\cite{AN-higherrank-voronoi}

The canonical measure and Arakelov--Bergman Green function are used to define invariants of Riemann surfaces, e.g., Faltings, Kawazumi--Zhang and Wilms invariants~\cite{Faltings84, Zhang2, Wilms, dJS21, Wilms21} and modular graph functions in string theory~\cite{HGB19}. In another forthcoming work, we apply the results of this paper and \cite{AN-hybrid-green} to study the asymptotics of modular graph functions on degenerating Riemann surfaces.

Hybrid spaces introduced in this paper refine previous constructions of hybrid spaces \cite{MS1, Berkovich, Acampo75, AB15, BJ17} and go beyond 
by producing higher rank compactifications which are themselves hybrid on their boundaries. The constructions presented here are related to higher rank non-Archimedean geometry. In this regard, a geometric study of higher rank valuations appears in~\cite{AI}, and the underlying higher rank polyhedral geometry is further developed by Hernan Iriarte in his PhD thesis~\cite{Iri21}.

\subsection{Notes  and references to related work} 

Properties of hyperbolic surfaces and their behavior close to the boundary of the moduli spaces are largely studied in the literature, see e.g.,~\cite{Wolp83, Wolp87, Wolp90, Wolp92, Wolp10, Ji93, GZ97} for a sample of results on the type of questions studied in this paper. The geometric behavior of Riemann surfaces close to the boundary of their moduli spaces is studied primarily in connection with Arakelov geometry and mathematical physics in~\cite{Ara74, Faltings84, Smit88, Jorgenson90, Wentworth, JW92, JK06, JK09, Zhang2, Faltings, deJong,  Went08,  HGB19, Wilms, dJS21}.

This paper fits  into the current research aiming to explain large scale limits of classical (e.g., algebraic and complex) geometry. In the past decade, this has given rise in particular to the development of tropical geometry and has spread into diverse fields of mathematics, ranging from algebraic and complex geometry to mathematical physics, combinatorics, computer science and optimization theory. Our results in this paper and its sequels show the interest of developing a higher rank version of tropical and non-Archimedean geometries.

The question of understanding non-Archimedean limits of complex manifolds, in particular, has been an active area of research. Pioneering work on the subject goes back already to the work of Bergman~\cite{Bergman}, Bieri--Groves \cite{BG}, Gelfand--Kapranov--Zelvinsky \cite{GKZ},  Passare and Rullg{\aa}rd \cite{PR04, Rug01}, Einsiedler--Kapranov--Lindand \cite{EKL}, Mikhalkin \cite{Mik04}, Jonsson \cite{Jonsson}, and others. Recently, it has emerged in the work of Kontsevich and Soibelman~\cite{KS}, in relation with mirror symmetry, in a conjectural formulation of the non-Archiemdean Calabi metric as the \emph{tropical limit} of maximal degenerations of Calabi--Yau varieties.  This has resulted in a series of works linking non-Archimedean and tropical geometry to metric limits of complex geometry, see e.g.,~\cite{KS, GW, GTZ, GTZ2, BFJ15, Sus18, Oda-Oshi, BJ18, YangLi}.

The first appearance of hybrid spaces can be traced back in the works of A'Campo \cite{Acampo75}  on monodromy zeta functions, Morgan--Shalen~\cite{MS1}, in the study of degenerations of hyperbolic structures on surfaces to real trees, and more recently in the work of Berkovich~\cite{Berkovich}, in providing a non-Archimedean interpretation of the weight zero part of the limit mixed Hodge structures, in degenerating families of complex algebraic varieties. Metrized complexes were used in~\cite{AB15} as a tool to study degenerations of linear series,  and found applications in the study of curves and their moduli spaces~\cite{ABBR, BJ16, FJP}. Recently, hybrid spaces have found applications in describing limits of Calabi-Yau measures in one parameter families of Calabi-Yau varieties~\cite{BJ17, Shiva2}, in the study of degenerations of K\"ahler-Einstein measures~\cite{PS}, and in dynamics in the works~\cite{Fav18, DKY, DF, Poi22,Poi22-2}. Parker has introduced the formalism of exploded manifolds in symplectic geometry~\cite{Par12} which is reminiscent of the works \cite{Acampo75, AB15} in that setting.

In connection with this paper, an equidistribution theorem for limits of Weierstrass points on families of complex curves was proved in~\cite{Ami-W}. The Archimedean version of this result is a theorem of Neeman and Mumford~\cite{Neeman}. A deeper link between the Zhang and Arakelov measures was naturally suggested by these results, the recent works of Faltings~\cite{Faltings} and R.\ de Jong~\cite{deJong} on non-Archimedean limits of Faltings $\delta$-invariants, that of~\cite{ABBF} on non-Archimedean limits of height pairing, and that of Baik--Shokrieh--Wu~\cite{SW19, BSW20} on Kazhdan theorems for metric graphs and Riemann surfaces.

In another direction related to the results of this paper, tropicalization of moduli spaces of curves has been studied by Abramovich--Caporaso--Payne~\cite{ACP}, and by Chan--Galatius--Payne in context with the cohomology of the moduli space of curves~\cite{CGP1}. 
Tropical moduli spaces have been used in~\cite{Oda, Oda2, Oda-Oshi}, in connection with the constructions of~\cite{BJ17}, to construct compactifications of moduli spaces of curves, Abelian varieties, and $K3$ surfaces. The hybrid moduli space $\mg_g^{^{\Hyb(1)}}$ which appears through the constructions given in Sections~\ref{sec:relations} and~\ref{sec:towermoduli} provides in a sense a refinement of these compactifications, replacing the dual complex with the dual metrized complex, making a more precise link between~\cite{Acampo75, AB15} and \cite{BJ17}. 
Moduli spaces of graphs have also been studied in the pioneering work of Culler--Vogtmann~\cite{CV86}.

Finally, we mention that degenerations of holomorphic differential forms and integrals to non-Archimedean differential forms and integrals have been studied in~\cite{CT10, Lagerberg, CLD, DHL, BGJK}, and connections to asymptotic Hodge theory appear in \cite{AP, AP-CS, GS06,  IKMZ,  Rud20}.

\subsection{Organization of the paper} This paper is structured as follows. Section~\ref{sec:prel} is of preliminary character, where we collect necessary definitions. In Section~\ref{sec:hybrid_spaces}, we define the hybrid space $B^\hyb$ associated to a complex manifold $B$ with SNC divisor $D$. The moduli space of hybrid curves $\mghyb{g}$ and its universal curve $\unicurve^\hyb_g$ are introduced in Section~\ref{sec:hybrid_moduli}. Section~\ref{sec:measures} contains the definitions of canonical measures on different geometric objects. Section~\ref{sec:tropical} concerns a purely tropical version of Theorem~\ref{thm:mainglobal-intro}. In Section~\ref{sec:monodromy} and Section~\ref{sec:period}, we study the period matrix of Riemann surfaces undergoing degeneration to a hybrid curve. Sections~\ref{sec:generic} and \ref{sec:proof} contain the proof of Theorem~\ref{thm:mainglobal-intro}. In Appendix~\ref{sec:AppendixTopology}, we define the topology on the versal hybrid family of a stable curve.

\subsection{Basic notation} \label{sec:basic-notations}
A curve in this paper means a complex projective algebraic curve. The analytification of a complex curve is a  compact Riemann surface.  We deliberately use the terminology curve and Riemann surface to switch back and forth between the algebraic and complex analytic setting, however, we essentially work in the latter setting. 

$\R_+$ denotes the set of strictly positive real numbers. We set $\imi \coloneq \sqrt{-1}$.

For a non-negative integer $n$, the symbol $[n]$ means the set $\{1,\dots, n\}$. "Almost all" in this paper means for all but a finite number of exceptions.

For a subset $F \subseteq E$, we denote by $F^c$ the complement of $F$ in $E$.

The letter $G$ is used for graphs.  We use the symbols $\mgr$,  $\hcurve^\trop$, $\mc$ and $\hcurve^\hyb$ for metric graphs, tropical curves, metrized complexes and hybrid curves, respectively. Moreover, when working with the family of curves over a versal deformation space, we use the notation $\rsf^\hyb$ for the corresponding family of hybrid curves over the hybrid deformation space $B^\hyb$.

In context with asymptotics, we will also make use of the standard Landau symbols. Let $X$ be a topological space, $x$ a point in $X$ and $g \colon U \setminus \{x\} \to (0, \infty)$ a function on a punctured neighborhood of $x$. Recall that the expressions $O(g)$ and $o(g)$ are used to abbreviate functions $f$ such that the quotient $|f(y)| / g(y)$ remains bounded and goes to zero, respectively, as $y$ goes to $x$ in $X$ (here, $f$ is defined on a punctured neighborhood of $x$ as well).

For an $n\times m$ matrix $A$, and a pair $i\in[n]$ and $j\in[m]$, the $(i,j)$-coordinate of $A$ is denoted by $A(i,j)$.

\subsection*{Acknowledgments}%
 We thank the referees for their careful reading of our manuscript, and for their helpful comments which enabled us to improve the presentation. We would like to thank Mirko Mauri for pointing out the reference~\cite{Acampo75}.

O.A. is part of the project ANR-18-CE40-0009.  N.N. acknowledges financial support by the Austrian Science Fund (FWF) under Grants No. P 28807, W 1245 and J 4497.


\section{Preliminaries}\label{sec:prel}
In this preliminary section we introduce the combinatorial objects appearing in the upcoming sections. These are graphs, ordered partitions, layered graphs and their graded minors, spanning trees, and augmented graphs.  We then provide a metric version of these notions, introducing metric graphs and tropical curves. We should clarify that our use of the latter notion  is a {multi-scale higher rank} refinement of the current usage in the literature. Tropical curves introduced here encode in addition the choice of a \emph{layering}, a concept introduced below, and {are uniformized} to have normalized edge lengths on each layer. This naturally leads to the definition of hybrid curves, as {uniformized} metrized complexes endowed with a layering on their underlying graphs.

\subsection{Graphs}  Part of the terminology and notation concerning graphs in this paper is standard either in graph theory or algebraic geometry. The aim of this section is to recall some of these notions. For more information, we refer to classical text books on graph theory~\cite{Bol, BM, Diestel} and algebraic geometry~\cite{ACGH}. All our graphs are finite. We allow parallel edges and loops.

 The vertex set and the edge set of a graph $G$ are denoted by $V(G)$ and $E(G)$, respectively. If $G$ is clear from the context, we simply use $V$ and $E$ and write $G=(V, E)$ to indicate the vertex and edge sets. 
 
 Two vertices $u$ and $v$ in $G$ are called \emph{adjacent} if there is an edge with extremities $u$ and $v$. In this case we write $u \sim v$. An edge $e$ and a vertex $v$ are \emph{incident} and we write $e\sim v$ if $v$ is an extremity of $e$.

 The \emph{degree} of a vertex $v$ is defined by 
 \[\deg(v) \coloneqq  \Bigl|\bigl\{e\in E \,\bigl|\, e \sim v\bigr\}\Bigr| + \Bigl|\bigl\{ e\in E \,\bigl|\, e \sim v \textrm{ and $e$ is a loop} \bigr\}\Bigr|.\]  
In other words, $\deg(v)$ counts the number of half-edges incident to $v$, so every loop contributes twice to the degree of the vertex.

 A \emph{subgraph} of $G=(V,E)$ is a graph $H=(W, F)$ with $W \subseteq V$ and $F \subseteq E$. We say that a subgraph $H=(W, F)$ of $G$ is \emph{spanning} if $W=V$, i.e., if it contains all vertices of $G$. For a collection of spanning subgraphs $H_1=(V, E_1), \dots, H_r = (V, E_r)$ of $G$, we denote by $H_1 \cup \dots \cup H_r$ the spanning subgraph of $G$ with edge set $E_1 \cup \dots \cup E_r$.

For a subset of edges $F \subseteq E$, we denote by $G[F]$ the spanning subgraph of $G$ with edge set $F$, i.e., $G[F] =(V, F)$.

By a \emph{spanning tree} of a graph $G=(V,E)$ we mean a spanning subgraph $T$ of $G$ without any cycle and with the same number of connected components as $G$. Equivalently,  $T$ has the maximum possible number of edges without containing any cycle. We denote by $\cT(G)$ the set of all spanning trees of $G$. Spanning trees are the central combinatorial structures underlying the definition of the canonical measure (see Section~\ref{ss:MeasuresGraphs}).

There are two simple operations that can be defined on graphs. These are called \emph{deletion} and \emph{contraction} and lead to the definition of minors of a graph. Despite their simplicity, these operations profoundly reflect the \emph{global properties of graphs}, when seen as a \emph{whole}, for instance in the context of Robertson--Seymour's theory of graph minors. Special minors associated to ordered partitions of edge sets, that we will call \emph{graded minors}, are introduced below and play a central role in this paper. 

Given a graph $G=(V, E)$ and an edge $e\in E$, we denote by $G-e$ the spanning subgraph $H$ of $G$ with edge set $E\setminus\{e\}$. In this case we say $H$ is obtained from $G$ by \emph{deleting} $e$. 
Moreover, we denote by $\contract{G}{e}$ the graph obtained by \emph{contracting} the edge $e$: it has as vertex set the set $U$ obtained by identifying the two extremities of  $e$, and as edge set, the set $E \setminus \{e\}$ viewed in $U$ via the projection map $V \to U$.
More generally, for a subset $F \subseteq E$, we denote by $G-F$ and $\contract{G}{F}$ the graphs obtained by deletion and contraction of $F$ in $G$, respectively: $G-F$ is the spanning subgraph of $G$ with edge set $E\setminus F$, and $\contract{G}{F}$ is the graph obtained by contracting one by one all edges in $F$. (One shows that the order does not matter.)  Later on, we will introduce a notion of contractions for \emph{augmented graphs}, which are graphs endowed with a genus function on vertices. In this case, the contraction remembers the sum of the \emph{genera of connected components of the contracted part}, thus keeping the genus of the graph constant before and after contraction.

A \emph{minor} of a graph $G$  is a graph $H$ which can be obtained by contracting a subset of edges in a subgraph $H'=(W, F)$ of $G$.

For a graph $G$, we denote by $H_1(G)$ the first homology group of $G$ with $\Z$ coefficients. The rank of $H_1(G)$, which coincides with the first Betti number of $G$ viewed as a topological space, is called \emph{genus} of the graph. In this paper we use the letter $\graphgenus$ when referring to the genus of graphs, and reserve the letters $g$ and $\genusfunction$ for the genus of augmented graphs, algebraic curves, and their analytic and hybrid variants.

\noindent The genus $\graphgenus$ of a graph $G=(V, E)$ is equal to $\abs{E} - \abs{V}+ c(G)$ where $c(G)$ denotes the number of connected components of $G$.

\subsection{Ordered partitions} \label{subsec:OrdPart} Let $E$ be a finite set. An \emph{ordered partition} of a (non-empty) subset $F \subseteq E$ is an ordered sequence $\pi = (\pi_1, \pi_2, \dots, \pi_r)$ of non-empty, pairwise disjoint subsets $\pi_i \subseteq F$, $i\in[r]$, such that
\begin{equation}
	F = \bigsqcup_{i=1}^r \pi_i.
\end{equation}
The integer $r$ is called the \emph{rank} of $\pi$. 
If $F = \varnothing$, then we make a slight abuse of notation and allow $\pi_\varnothing = \varnothing$ as its only ordered partition. The set of ordered partitions of a subset $F \subseteq E$ is denoted by $\Pi(F)$ (in particular, $\Pi(\varnothing) = \{\pi_\varnothing\}$). 

The ordered partitions of a subset $F$ are in bijective correspondence with its \emph{filtrations}:  in this paper, a \emph{filtration} $\filter$ of \emph{depth} $r$ on a (non-empty) subset $F \subseteq E$ is a strictly increasing sequence of non-empty subsets
\[\filter: \quad  F_1 \subsetneq F_2  \subsetneq F_3 \subsetneq \dots \subsetneq F_r=F\] 
for $r \in \mathbb N$. Clearly, each filtration $\filter$ on a subset $F \subseteq E$ defines an ordered partition $\pi = (\pi_i)_{i=1}^r$ of $F$ by setting
\[
	\pi_i = F_i \setminus F_{i-1}, \qquad i=1,\dots,r,
\]
with $F_0\coloneqq \varnothing$ by convention, and this gives a bijection. Sometimes we use the notation $\filter^\pi$  to emphasize the associated ordered partition $\pi$. Again, for $F = \varnothing$ we allow $\filter^\varnothing \coloneqq  \varnothing$ as the only filtration. 

Given an ordered partition $\pi=(\pi_1, \dots, \pi_r)$, we set $E_\pi \coloneqq  \pi_1 \cup \dots \cup \pi_r$, the subset of $E$ underlying the ordered partition.

We now define a natural partial order  on the set of ordered partitions (equivalently, filtrations) of subsets of a given set $E$. A filtration $\filter$ is a \emph{refinement} of the filtration $\filter'$, if $\filter' \subseteq \filter$ as a set. That is, denoting by $F_1 \subsetneq \dots \subsetneq F_r=F$ the elements of $\filter$, and by $F_1' \subsetneq \dots \subsetneq F_{s}' =F'$ the elements of $\filter'$, for each $i\in[s]$, there exists $j\in[r]$ such that $F_i'=F_j$.  We say that an ordered partition $\pi$ is a \emph{refinement} of an ordered partition $\pi'$, and write $\pi' \preceq \pi$, if this holds true for the corresponding filtrations  $\filter$ and  $\filter'$. It is easy to see that $\pi' \preceq \pi$ exactly when $\pi$ is obtained from $\pi'$ by performing the following two steps:
\begin{itemize}[leftmargin=2em]
\item[(i)] Replacing each set $\pi_i'$ in the ordered sequence $\pi' = (\pi_i')_{i=1}^r$ by an ordered partition ${\varrho}^i = ({\varrho}_k^i)_{k=1}^{s_i}$ of $\pi_i'$.
\item[(ii)] Possibly adding an ordered partition $\varrho=(\varrho_k)_{k=1}^s$ of a subset $F^\infty \subseteq E\setminus E_{\pi'}$ at the end of the ordered sequence $\pi'$. 
\end{itemize}

(The use of the symbol $F^\infty$ will become clear when we link these constructions to geometry.)

We stress that in the definition of a refinement, the ordered partitions and filtrations are not necessarily defined over the same subset of $E$ (see Example \ref{ex:op}). For instance, $\pi_\varnothing \preceq \pi$ for all ordered partitions $\pi$. However, if $\pi' \preceq \pi$ (or equivalently, $\filter' \subseteq \filter$), then this implies $E_{\pi'} \subseteq E_\pi$. It is also clear that "$\preceq$" defines a partial order on $ \Pi \coloneqq \bigcup_{F\subseteq E} \Pi(F)$, the set of all ordered partitions of subsets $F \subseteq E$.
\begin{example} \label{ex:op}
Assume that $E = \{e,f \}$ for two elements $e \neq f$. Then, the ordered partitions of subsets $F \subseteq E$ are given by
\begin{align*}
	\Pi ( \varnothing ) &= \big \{ \pi_\varnothing \big  \}, \qquad \Pi ( \{e\}) = \Big \{ \big (\{e \} \big )\Big  \}, \qquad
	\Pi( \{f\}) =\Big  \{\big  (\{f \} \big ) \Big \}, \\
	\Pi(E) &= \Big \{ \big (\{e\}, \{ f\} \big ),\big  (\{f\}, \{ e\} \big ),\big  (\{e, f\}\big  )  \Big  \}.
\end{align*}
The partial order "$\preceq$" on $\Pi = \bigcup_{F\subseteq E} \Pi(F)$ consists of the following relations
\begin{align*}
	\big ( \{e\} \big ) & \preceq \big  (\{e\}, \{ f\} \big ),  \qquad
	\big ( \{f\}\big )  \preceq \big  (\{f\}, \{ e\} \big ), \\
	\big ( \{e,f\}\big ) &\preceq \big  (\{e\}, \{ f\} \big ), \qquad
	\big ( \{e,f\}\big ) \preceq \big (\{f\}, \{ e\} \big ),
\end{align*}
and the trivial relations that $\pi_\varnothing \preceq \pi$ and $\pi \preceq \pi$ for all $\pi \in \Pi$.
\end{example}

\subsection{Layered graphs and their graded minors}\label{ss:GradedMinors} A \emph{layered graph}  is the data of a graph $G=(V, E)$ and an ordered partition $\pi=(\pi_1, \dots, \pi_r)$ on the set of edges $E$. We call $\pi$ a \emph{layering} and the elements $\pi_j$ of $\pi$ the \emph{layers}.
The integer $r$ is called the \emph{rank} of the layered graph $G$, which thus equals the rank of the ordered partition $\pi$.  By an abuse of the notation, we use the same letter $G$ to denote the layered graph $(V, E, \pi)$. The genus of a layered graph is defined as that of its underlying graph. 

Let now $G=(V, E, \pi)$ be a layered graph with the ordered partition $\pi = (\pi_1, \dots, \pi_r)$ of $E$. We define the \emph{graded minors} of $G$ as follows.

Consider the filtration 
\[\filter^\pi_\bullet \colon F_0 =\varnothing \subsetneq F_1 \subsetneq F_2 \subsetneq \dots \subsetneq F_r = E\] 
associated to $\pi$. For each $j$, we let 
\[E^j_\pi\coloneqq E\setminus F_{j-1} = \pi_j \cup \pi_{j+1} \cup \dots \cup \pi_r \]
and obtain a decreasing filtration on $E$,
\[
\dfilter_\pi^\bullet \colon E_\pi^1 = E \supsetneq E_\pi^2 \supsetneq \dots E_\pi^r \supsetneq E_\pi^{r+1} = \varnothing.
\]
This leads to a decreasing sequence of spanning subgraphs of $G$
\[
G \eqqcolon G^1_\pi \supset G^2_\pi \supset G_\pi^3 \supset \dots \supset G_\pi^{r} \supset G_\pi^{r+1}= (V, \varnothing)
\]
where for each $j$,  $G_\pi^j\coloneqq (V, E_\pi^j)$ has edge set $E_\pi^j$.

For each integer $j=1, \dots, r$, the $j^{\mathrm{th}}$ \emph{graded minor} of $G$, denoted by $\grm_\pi^j(G)$, is obtained by contracting all the edges of $E_{\pi}^{j+1}$ 
in $G_\pi^j$, i.e., 
\begin{equation} \label{eq:DefGradedMinors}
\grm_\pi^j(G) \coloneqq  \contract{G^{j}_\pi}{E^{j+1}_\pi} = \contract{G[F_{j-1}^c]}{F_j^c} = \contract{G[\pi_j\cup \dots \cup\pi_r]}{\pi_{j+1} \cup \dots \cup \pi_r}.
\end{equation}

\begin{figure}[!t]
\centering
    \includegraphics[scale =0.3]{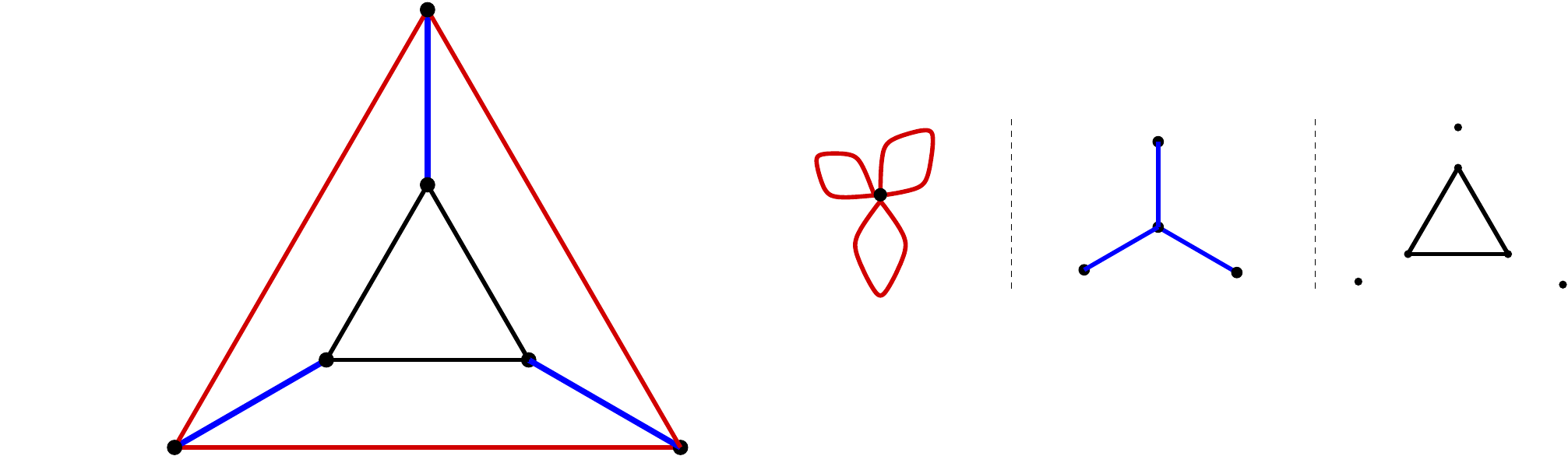}
\caption{A layered graph with ordered partition $\pi=({\color{red}\pi_1}, {\color{blue}\pi_2}, \pi_3)$, with three layers, and its graded minors $\grm_\pi^1(G)$, $\grm_\pi^2(G)$ and $\grm_\pi^3(G)$, from left to right.}
\label{fig:layered_graph}
\end{figure}

It follows that $\grm_\pi^j(G)$ has edge set equal to $\pi_j$. More precisely, we should view the elements of $\pi_j$ in the contracted graph, but we simplify and identify the two sets in the writing. An example is depicted in Figure~\ref{fig:layered_graph}.
We denote by $V_\pi^j$ the vertex set of $\grm_\pi^j(G)$. If $\pi$ or $G$ is understood from the context, we will simply write $\grm^j(G)$, $\grm^j_\pi$, $V^j$, \emph{etc.} We denote by $\cont{j} \colon G_\pi^j \to \gr^j_\pi(G)$ the contraction map.

\subsection{Genus formula for layered graphs} Let $G=(V, E, \pi)$ be a layered graph. We denote by $c_\pi^j$ the number of connected components and by $h^j_\pi$ the genus of $\grm^j_\pi(G)$. These are related by $h^j_\pi = |\pi_j| - |V^j_\pi|+c^j_\pi$.

 \begin{prop}[Genus formula]\label{prop:genusformula} \label{prop:genusformula} Let $G$ be a layered graph of genus $\graphgenus$ with layering $\pi=(\pi_1, \dots, \pi_r)$. Then, we have the equality
 \[\graphgenus = \sum_{j=1}^r h_\pi^j.\]
 \end{prop}
 In preparation for the proof, consider $\gr^1_\pi(G) = (V_\pi^1, \pi_1) $. For each vertex $u\in V_\pi^1$, let $V_u \subseteq V$ be the set of vertices in $V$ that are contracted to $u$, i.e., $V_u = \continv{1}(u)$. Denote by $E_u$ the set of edges in $E_\pi^2$ with both end-points in $V_u$. Define $G_u\coloneqq (V_u, E_u)$.
 
 \begin{prop} The subgraph $G_u = (V_u, E_u)$ of $G$ is connected.
 \end{prop}
 \begin{proof} The contraction of the edges of $G_u$ results in a graph with a single vertex $u$. This can only happen if the graph is connected. \end{proof}
 The ordered partition $\pi$ of $E$ induces an ordered partition $\pi_u$ of $E_u$: we take first the ordered partition
 \[\pi_1\cap E_u =\varnothing, \quad \pi_2\cap E_u, \quad  \pi_3\cap E_u, \quad \dots, \quad \pi_r \cap E_u,\]
 and then remove all the empty sets from the sequence. The pair $(G_u, \pi_u)$ is a layered graph. Note that $\pi_u$ has rank at most $r-1$ since the first  set $\pi_1 \cap E_u$ is empty.

\begin{proof}[Proof of Proposition~\ref{prop:genusformula}]  This can be obtained by induction on the rank of the ordered partition $\pi$. If the partition $\pi$ consists of a single set $\pi_1=E$, the proposition is obvious. Otherwise, let $G_\pi^2=(V,  E_\pi^2)$ be the spanning subgraph of $G$ with edge set $E_\pi^2$. The set $E_\pi^2$ is a disjoint union of the sets $E_u$ for $u\in V^1_\pi$. It follows that $G^2_\pi$ is the disjoint union of graphs $G_u$ for $u\in V_\pi^1$. By the previous proposition, the number of connected components of $G_\pi^2$ is $|V^1_\pi|$, and its genus is thus given by 
\begin{equation}\label{eq:genusformula}
\genus(G^2_\pi) = |E_\pi^2| - |V| + |V^1_\pi| = \sum_{u\in V_\pi^1} \genus(G_u).
\end{equation}
Using the induction hypothesis for the layered graphs $G_u$, we conclude.
\end{proof}


\subsection{Spanning trees of layered graphs} \label{ss:LayeredSpanningTrees}
Let $G=(V, E)$ be a graph endowed with an ordered partition $\pi = (\pi_1, \dots, \pi_r)$ of $E$. Consider the graded minors $\gr^1_\pi(G), \dots, \gr^r_\pi(G) $ introduced in \eqref{eq:DefGradedMinors}.

For each $j=1, \dots, r$, let $T_j$ be a spanning tree of $\gr^j_\pi(G)$ with edge set $A_j \subset \pi_j$, so that $T_j = (V^j_\pi, A_j)$. Let $A = A_1 \cup \dots \cup A_r$. The following proposition can be regarded as a refinement of the genus formula.

\begin{prop} \label{prop:spanning_tree_layered} The subgraph $T = (V, A)$ of $G$ is a spanning tree.
\end{prop} 

\begin{defi}[Spanning trees of layered graphs] \label{def:LayeredSpanningTree} A \emph{spanning tree of a layered graph} $G=(E, V, \pi)$ is a spanning tree $T= (V, A)$ of $G$ which is obtained as above, i.e., $A = A_1 \cup \dots \cup A_r$ for $A_j$ the edge set of a spanning tree of $\gr^j_\pi(G)$. We denote by $\cT_\pi(G)$ the set of all such spanning trees. 
\end{defi}
\begin{remark} We have $\cT_\pi(G) \subseteq \cT(G)$ but the two sets are different in general.  An example is given in Figure~\ref{fig:layered_spanning_trees}.
\end{remark}

\begin{figure}[!t]
\centering
    \scalebox{.28}{\input{layered_spanning_trees.pspdftex}}
\caption{A layered graph $(G, \pi)$ with two layers and its four spanning trees $T_1, T_2, T_3,$ and $T_4$. The underlying non-layered graph $G$ has five spanning trees, those of $(G, \pi)$ and the spanning tree $T$ depicted in the figure.}
\label{fig:layered_spanning_trees}
\end{figure}

\begin{proof}[Proof of Proposition~\ref{prop:spanning_tree_layered}] The sets $A_1, \dots, A_r$ are disjoint and $\abs{\pi_j \setminus A_j} = h^j_\pi$. By Proposition~\ref{prop:genusformula}, we see that $\abs{E \setminus A} = h$. So it suffices to prove that $T$ contains no cycles. For the sake of a contradiction, suppose $T$ has a cycle $C\subset G$, and let $j$ be the smallest index with $E(C) \cap \pi_j \neq\varnothing$. But then, the set of edges $E(C)\cap \pi_j$ form a cycle in $\gr^j_\pi(G)$, which is a contradiction. 
\end{proof}


\subsection{Marked graphs} Let $n$  be a non-negative integer. A \emph{graph with $n$} (\emph{labeled}) \emph{marked points} is a graph $G=(V, E)$ endowed with the data of labels $1, \dots, n$ placed on its vertices. More formally, this is encoded in the marking function $\marking \colon [n] \to V$.

To a graph with $n$ marked points, we associate the \emph{counting function} $\countmarking \colon V \to \mathbb N \cup \{0\}$ 
which counts the numbers of labels at each vertex. That is, $\countmarking(v)$ is the number of elements $j \in [n]$ with $\marking(j)=v$.


\subsection{Augmented and stable graphs} We now consider graphs endowed with a \emph{genus function} on the vertex set that we call \emph{augmented graphs}, using the terminology of \cite{ABBR}. Note that the same objects are sometimes called weighted graphs in the literature.

An \emph{augmented graph} is a graph $G=(V, E)$ endowed with the data of an integer-valued \emph{genus function} $\genusfunction \colon V \to \mathbb N \cup\{0\}$ on the vertex set. The integer $\genusfunction(v)$ is called the genus of the vertex $v \in V$. The genus of an augmented graph $G=(V, E, \genusfunction)$ is defined as
\[g\coloneqq  h + \sum_{v\in V} \genusfunction(v),\]
where $h=\abs{E} - \abs{V} + c(G)$ is the genus of the underlying graph $(V, E)$. 

An \emph{augmented graph with $n$ marked points} is an augmented graph $(V, E, \genusfunction)$ with a marking $\marking \colon [n] \to V$.

We mostly consider augmented (marked) graphs which verify the following stability condition.   A \emph{stable graph with $n$ marked points} is a quadruple $G=(V, E, \genusfunction, \marking)$ with marking and genus functions  $\marking\colon [n] \to V$ and $\genusfunction \colon V\to \mathbb N \cup\{0\}$, respectively,  subject to the following \emph{stability condition}: 
 
 \begin{itemize}[leftmargin = 2em]
\item for every vertex $v$ of genus zero, we have $\deg(v) + \countmarking(v) \geq 3$, and
\item for every vertex $v$ of genus one, we have $\deg(v) + \countmarking(v)\geq 1$.
\end{itemize}

 A \emph{stable graph}  is a triple $(V, E, \genusfunction)$ (without any marking) that verifies these conditions for the counting function $\countmarking \equiv0$.
 
By an abuse of the notation, we use $G$ both for the stable (marked) graph  and its underlying graph $(V, E)$, forgetting the genus function (and marking).

The notion of deletion and contraction can be extended to stable marked graphs. (The construction works more generally without the stability condition. However, in our applications, we restrict to the stable case.) Let $G = (V, E, \genusfunction, \marking)$ be a stable marked graph of genus $g$ with $n$ marked points. For an edge $e \in E$, we define the stable marked graph $\contract{G}{e} = :G' = (V', E', \genusfunction', \marking')$ as follows. First, the underlying graph $G'$ of $\contract{G}{e}$ is the graph obtained by contracting $e$ in $G$. Let $v_e$ be the vertex of $G'$ which corresponds to the identification of the two extremities $u$ and $v$ of $e$. The genus function $\genusfunction'$ of $\contract{G}{e}$ is then defined by 
\[\genusfunction'(w) = \begin{cases} \genusfunction(w) & \textrm{if $w \neq v_e$} \\
\genusfunction(u)+\genusfunction(v) & \textrm{if $w = v_e$ and $e$ is not a loop, i.e., $u\neq v$}\\
\genusfunction(u)+1 & \textrm{if $w=v_e$ and $e$ is a loop, i.e., $u=v$}.\end{cases}\]

In the presence of a marking function $\marking \colon  [n] \to V$, we define the marking function $\marking'$ of $G'$ as $\marking' \coloneqq  \proj \circ \marking: [n] \to V'$, where $\proj \colon G \to G'$ denotes the contraction map.

 The following proposition is straightforward. 
 
 \begin{prop} Notations as above, the stable graph $\contract{G}{e}$ has genus $g$.
 \end{prop}

 Finally, we define layered stable (marked) graphs as stable (marked)  graphs $G=(V, E, \genusfunction, \marking)$ endowed with an ordered partition $\pi$ on the edge set $E$.
 
 \subsection{Partial order on stable (marked) graphs of given genus} \label{sec:PartialOrderMarkedGraphs}
Let $n$ be a non-negative integer. We define a partial order $\subface$ on the set of stable marked graphs of genus $g$ with $n$ marked points as follows. For two stable marked graphs $G$ and $H$, we say $G \subface H$ if $G$ can be obtained  by a sequence of edge-contractions from $H$, in the sense of the previous paragraph.

\subsection{Metric graphs} \label{ss:MetricGraphs} We briefly recall the definition of a \emph{metric graph} and how it arises as a \emph{metric realization} of a graph with associated \emph{edge lengths}.

Let $G=(V,E)$ be a finite graph and $\ell\colon E \rightarrow \mathbb R_{+}$ an edge length function that assigns a positive real number $\ell_e$ to every edge $e \in E$. To such a pair $(G, \ell)$, we associate a metric space $\mgr$ as follows: by assigning each edge a direction and calling one of its vertices the initial vertex $v_i$ and the other one the terminal vertex $v_{t}$, every edge $e\in E$ can be identified with a copy of the interval $\Ical_e= [0,\ell_e]$ (the left and right endpoint correspond to $v_i$ and $v_t$). The topological space $\mgr$ is obtained by further identifying the ends of edges corresponding to the same vertex $v$ (in the sense of a topological quotient). The topology on $\mgr$ is metrizable by the so-called {\em path metric}: the distance between two points $x,y \in\mgr$ is defined as the arc length of the shortest path connecting them.

A \emph{metric graph} is a compact metric space arising from the above construction for some pair $(G, \ell)$ of a graph $G$ and edge length function $\ell$. The metric graph $\mgr$ is then called the \emph{metric realization} of the pair $(G, \ell)$. Conversely, a pair $(G, \ell)$ such that its metric realization is isometric to a fixed metric graph $\mgr$ is called a \emph{finite graph model} of $\mgr$. Clearly, any metric graph has infinitely many finite graph models. However, as we explain now, for stable metric graphs it is possible to find a \emph{minimum model}.

  An \emph{augmented metric graph} is a metric graph with a genus function $\genusfunction \colon \mgr \to \N \cup\{0\}$ which is zero outside a finite set of points in $\mgr$.
A \emph{metric graph with $n$ marked points} is a metric graph $\mgr$ endowed with a marking function $\marking \colon [n] \to \mgr$. We define the counting function $\countmarking \colon \mgr \to [n]$ similar to the case of graphs.

A \emph{finite graph model of an augmented metric graph with $n$ marked points} is a finite graph model $(G=(V,E), \ell)$ of the underlying metric graph such that $V$ contains all points of positive genus and all marked points.

An augmented metric graph $\mgr$ with $n$ marked points is called \emph{stable} if every point $x$ of $\mgr$ of degree one verifies  $\deg(x) +\genusfunction (x)  + \countmarking(x) \geq 3$, and any connected component of $\mgr$ entirely consisting of points of degree two (i.e., homeomorphic to circle) has either a point of positive genus or a marked point. The following is well-known. 
\begin{prop} Any stable metric graph has a unique minimum model, that is a finite graph model, whose vertex set is contained in the vertex set of any other finite graph model.
\end{prop}
\begin{proof}
The vertex set of any finite graph model contains all points $x$ with $\deg(x) \neq 2$, all points $x$ with $\genusfunction(x)>0$, and all points $x$ with $\countmarking(x) >0$. However, under the stability condition, these points themselves form the vertex set of a finite graph model.
\end{proof}

\subsection{Layered metric graphs and tropical curves}  A \emph{layered augmented metric graph} (\emph{with marking}) is a pair $(\mgr, \pi)$ consisting of an augmented metric graph $\mgr$ (with marking) together with a finite graph model $(G, \ell)$ and an ordered partition $\pi$ on the edge set $E(G)$ of $G$.  In this situation, $(G, \pi)$ is called the \emph{combinatorial type} of the augmented metric graph $(\mgr, \pi)$.

Two layered augmented metric graphs (with marking) $\mgr$ and $\mgr'$ are called \emph{conformally equivalent at infinity} if they have the same combinatorial type $(G, \pi)$, the same genus function on the vertex set $V(G)$ and for each layer $\pi_j$ of $\pi$, there is a positive number $\lambda_j >0$ such that
\[\ell\rest{\pi_j} = \lambda_j \, \ell'\rest{\pi_j}
\]
for the edge length functions $\ell, \ell' \colon  E(G) \to \R_{+}$. In other words, the layered augmented metric graphs $(\mgr', \pi')$ equivalent to some fixed $(\mgr, \pi)$ are obtained by multiplying the lengths for every layer $\pi_j$ by constants $\lambda_j >0$. 

In the following we will refer to this as  \emph{the conformal equivalence relation} (at infinity being assumed at all time).

A \emph{tropical curve} (\emph{with marking}) is a conformal equivalence class of layered metric graphs (with markings). Equivalently, this is a pair $\hcurve^\trop = (\mgr, \pi)$ consisting  of an augmented metric graph $\mgr$ (with marking) with an ordered partition $\pi=(\pi_1, \dots, \pi_r)$ on the edge set $E$ of a finite graph model $(G, \ell)$ of $\mgr$ such that in addition, the normalization conditions
\[\sum_{e\in \pi_j} \ell_e =1, \quad j=1, \dots, r\]
hold true.
A layered metric graph  (with marking) is called \emph{stable} if its underlying metric graph (with marking) is stable. A tropical curve (with marking) is called \emph{stable} if any, and so each, of its layered metric graph (with marking) representatives are stable. For a stable tropical curve (with marking), we \emph{always} assume that the layering is defined on the minimum model of the underlying metric graph.

\begin{remark} \label{rem:NewDefinitionTropicalCurves}
 In~\cite{AN-hybrid-green} we will extend the definition of tropical curves. 
The tropical curves appearing in the current paper are those having \emph{full sedentarity}, that is, all their edges  live at infinity (the word \emph{infinity} in the conformal equivalence relation refers to this). A tropical curve in general can have a \emph{finite part} meaning edges which do not live at infinity.  Considering them all at the same time will allow to define a compactified moduli space of higher rank tropical curves in~\cite{AN-hybrid-green}. The notion of rank in the definition of tropical curves refers to the number of layers at infinity.
\end{remark}

\subsection{Layered metrized complexes and hybrid curves} \label{sec:hybrid_curves} We first recall the definition of metrized complexes, see~\cite{AB15} for more details.  

A {\em metrized curve complex $\mc$} consists of the
following data:
\begin{itemize}[leftmargin = 2em]
\item A finite graph $G=(V, E)$.
\item A metric graph $\mgr$ with a model $(G, \ell)$ for a length function $\ell\colon E \to \R_{+}$.  
\item For each vertex $v\in V$, a smooth projective complex curve $C_v$. 
\item For each vertex $v \in V$, a bijection $e \mapsto p^e_v$ between the edges incident to $v$ (with loop edges counted twice) and a subset 
$\mathcal A_v = \{ p^e_v \}_{e \ni v}$ of $C_v(\C)$.
\end{itemize}
By an abuse of the notation, we use the same letter $\mc$ for the metrized complex and its \emph{geometric realization} defined as follows. For each edge $e \in E$, let $\Ical_e = [0, \ell_e]$ be an interval of length $\ell_e$. For each vertex $v \in V$, by an abuse of the notation, denote equally by $C_v$ the analytification, which is a compact Riemann surface. For each vertex $v$ of an edge $e$, we identify the corresponding extremity of $\Ical_e$ with the marked point $p^e_v$. This identifies $\mc$ as the disjoint union of the Riemann surfaces $C_v$, $v \in V$, and intervals $\Ical_e$, $e \in E$, glued by these identifications. $\mc$ is naturally endowed with the quotient topology.

To any metrized complex, one associates the underlying graph $G = (V,E)$ and the genus function $\genusfunction \colon V \to \mathbb N\cup\{0\}$ which maps $v \in V$ to the genus $\genusfunction(v)$  of $C_v$. The metric realization of this augmented graph and the edge length function $\ell$ is called the underlying augmented metric graph of $\mc$.

A \emph{metrized complex with $n$ marked points} is a metrized complex $\mc$ with a marking function $\marking \colon [n] \to \mc$ such that for each $j\in [n]$, $\marking(j)$ lies outside the union of the intervals $\Ical_e$, $e \in E$.%

A \emph{layered metrized complex} is a pair $(\mc, \pi)$ consisting of a metrized complex  $\mc$ and  a layering $\pi$ on the edge set of $\mc$. 
Any layered metrized complex naturally gives rise to a layered augmented metric graph.

We define \emph{the conformal equivalence relation} (at infinity) on layered metrized complexes (with marking) in the same way as for layered metric graphs.

Finally, a \emph{hybrid curve} (\emph{with marking}) $\hcurve^\hyb$  is a conformal equivalence class of layered metrized complexes (with marking). Equivalently, we can define a hybrid curve (with marking) $\hcurve^\hyb$ as a layered metrized complex $(\mc, \pi=(\pi_1, \dots, \pi_r))$ whose edge length function verifies the normalization condition $\sum_{e\in \pi_j} \ell_e =1$ for all $ j=1, \dots, r.$

Any hybrid curve $\hcurve^\hyb$ naturally gives rise to a tropical curve $\hcurve^\trop$: we take the underlying augmented metric graph $\mgr$ of $\mc$ with normalized edge lengths, together with the layering $\pi$ of the edge set $E$ of the graph $G$, which is a finite graph model of $\mgr$. 

In this paper, unless otherwise stated, we only consider stable metrized complexes and hybrid curves (with marking), which are those whose underlying augmented (marked) graph is stable.

\section{Hybrid spaces of higher rank} \label{sec:hybrid_spaces}
The aim of this section is to present a construction of hybrid spaces of higher rank. Given the data of a complex manifold $B$ and a simple normal crossing divisor  $D$, we associate to it a hybrid space $B^{\hyb}$ by equipping the points $t \in D$ with additional simplicial coordinates. Specializing to the case of a polydisc $B = \Delta^{\ndim}$ with the SNC divisor given by the coordinate axes, we obtain a suitable base space for the hybrid versal deformation in Section \ref{sec:hybrid_base}.

\subsection{The setting} \label{ss:setup}
In the following, let $B$ be an $\ndim$-dimensional complex manifold, for an integer $\ndim$, and let $D = \bigcup_{e \in E} D_e$ be a simple normal crossing divisor. More precisely, we require that $(D_e)_{e \in E}$ is a finite family of smooth, connected and closed submanifolds of codimension one in $B$ such that for any subset $F \subseteq E$, the intersection 
\begin{equation} \label{eq:divstrat} D_F \coloneqq  \bigcap_{e \in F} D_e
\end{equation}
is either empty or a smooth submanifold of codimension $|F|$ (with only finitely many connected components). We refer to the intersection $D_F$ as the \emph{stratum} associated to the subset $F \subseteq E$. The \emph{inner stratum} is the subset
\begin{equation} \label{eq:interiorcomp}
{\inn D_F} \coloneqq  D_F \setminus \bigcup_{e \notin F} D_e = \bigl\{t \in B \,\st\, E_t = F \bigr\}
\end{equation}
where for $t \in B$, we set
\begin{equation}
E_t \coloneqq  \bigl\{e \in E\,\st \, t \in D_e \bigr\}.
\end{equation}
Notice that for $F = \varnothing$, we simply recover the \emph{open part}
\[
 {\inn D_\varnothing} = B \setminus D \eqqcolon B^\ast.
\]
The inner strata \eqref{eq:interiorcomp} form a partition of $B$, that is,
\begin{equation} \label{eq:partitionB}
B = \bigsqcup_{\substack{ F \subseteq E }} {\inn D_F}.
\end{equation}
An \emph{adapted coordinate neighborhood}, also called a \emph{system of local parameters}, for $D$ around a point $t \in B$ is a pair $(U,z)$ where $U$ is an open neighborhood of $t$ in $B$ and $z=(z_i)_{i=1}^{\ndim}$ are local coordinates on $U$ with

\begin{itemize}[leftmargin = 2em]
\item [(i)] $|z_i|< 1$ on $U$ for all $i = 1,\dots,\ndim$, and

\item [(ii)] $D_e \cap U = \varnothing$ for all $e \notin E_t$, and

\item [(iii)] $D_e \cap U =\bigl\{s \in U\,\st \; z_e(s) = 0 \bigr\}$ for all $e\in E_t$. More precisely, for each $e \in E_t$, there is an associated coordinate $z_e \coloneqq  z_{i(e)}$ on $(U,z)$ such that this holds true.
\end{itemize}

We will also need the following elementary fact:
\begin{prop} \label{prop:adapcord}
Assume that $(U, z)$ and $(U',z')$ are adapted coordinate neighborhoods for points $t, t' \in B$, respectively. Then, for all $e \in E_t \cap E_{t'}$,
\[
	z_e ' = g_e z_e \qquad \text{on } U \cap U'
\]
for a non-vanishing holomorphic function $g_e \colon U \cap U' \to \C$.
\end{prop}

\subsection{Definition of the hybrid space} \label{ss:bhyb}
Before formally defining the hybrid space $B^{\hyb}$, we briefly motivate the construction. The idea is to enrich the points $t \in D$ with additional \emph{simplicial coordinates} in order to distinguish different ways of approaching them in the complex manifold $B$.

Namely, fix some $t \in B$ and consider a small neighborhood of $t$ in $B$. Upon choosing local coordinates $(z_i)_i$, we can assume that the divisors $D_e$, $e \in E_t$, are given by $D_e = \{z_e = 0 \}$. This trivially implies that
\begin{equation} \label{eq:expalpha}
	\lim_{s \to t}  \log(|z_e(s)|)  = - \infty, \qquad e \in E_t.
\end{equation}
However, the asymptotic behavior of the normalized logarithmic coordinates,
\begin{equation}
\Log_e(s)\coloneqq \frac{\log(|z_e(s)|)}{\sum_{\hat e\in E_t} \log(|z_{\hat e}(s)|)},\qquad e\in E_t,
\end{equation}
carries non-trivial information. This hints that sequences $(s_n)_n$ in $B^\ast$ converging to $t$ in $B$ with different limit values
\begin{equation}
	y_e \coloneqq  \lim_{n \to \infty} \Log_e (s_n) \in [0,1]
\end{equation}
should be further distinguished and hence converge to different points in $B^{\hyb}$. Since $\sum_e \Log_e(s)=1$, the limit $y\coloneqq (y_e)_e$ lies in the simplex $\sigma!_{E_t}=\bigl\{(y_e)_e \in \ssR_{\ge0}^{E_t}\,\st \,\sum_{e} y_e=1\bigr\}$. On the other hand, there might be a subset $E_t^2 \subsetneq E_t$ with
\begin{equation}
	y_e = \lim_{n \to \infty} \Log_e (s_n) = 0, \qquad e \in E_t^2.
\end{equation}
These logarithmic coordinates are dominated by the others and their limit $y_e = 0$ carries only little information. To analyze them further, we renormalize and consider the limits
\begin{equation} \label{eq:expomega}
y_e^2\coloneqq \lim_{n \to \infty} \frac{\log(|z_e(s_n)|)}{\sum_{\hat e \in E_t^2} \log(|z_{\hat e}(s_n)|)}, \qquad e \in E_t^2.
\end{equation}
Repeating the same steps, we end up with subsets $E_t \eqqcolon E_t^1 \supsetneq E_t^2\supsetneq  \cdots \supsetneq E_t^r\supsetneq E_t^{r+1} =\varnothing$  and corresponding limit values $ y\eqqcolon y^1 \in \sigma!_{E_t^1}$,  $\cdots$, $y^r \in \sigma!_{E_t^r}$.
The idea behind the topology of $B^{\hyb}$ is as follows: sequences $(s_n)_n \subseteq B^\ast$ approaching the same point $t$ in $B$ have different limit points in $B^{\hyb}$, according to the asymptotic behavior given by the values $(y^1,\cdots, y^r)$.

We now proceed with the formal definition of the space $B^{\hyb}$. Recall from Section \ref{subsec:OrdPart} that $\Pi(F)$ denotes the set of ordered partitions of a subset $F \subseteq E$ (in particular, $\Pi(\varnothing) = \{\pi_\varnothing\}$). We denote by $\Pi$ the union of $\Pi(F)$ for $F \subseteq E$.

For any subset $F \subseteq E$, the standard simplex $\sigma!_F$ in $\R^F$ is denoted by
\[  \label{eq:simplex}
\sigma!_F =  \Bigl\{ (x_e)_e \in  \ssR_{\ge 0}^F \, \,\bigl|\,\, \sum_{e \in F } x_e = 1\Bigr\}
 \]
and its relative interior $\ssub{\inn \sigma}!_F$ is given by
\[
\ssub{\inn \sigma}!_F =  \Bigl\{ (x_e)_e \in  \sigma!_F \, \,\bigl|\,\, x_e > 0 \text{ for all } e \in F\Bigr\}.
 \]
For an ordered partition $\pi =(\pi_1, \dots, \pi_r) \in \Pi(F)$ of a subset $F \subseteq E$, we define its \emph{hybrid stratum}  as 
\begin{equation} \label{eq:hybridstratum}
{D}^{\hyb}_\pi \coloneqq  \inn D_\pi \times \inn  \sigma_{\pi}   \coloneqq   \inn D_{\pi}  \times \inn\sigma_{\pi_1} \times \inn \sigma_{\pi_2} \dots  \times \inn \sigma_{\pi_r},
\end{equation}
where $\inn D_{\pi} \coloneqq  \inn D_{F}$. For $\pi = \pi_\varnothing$ (i.e., the empty ordered partition of $F= \varnothing$), we set
\[
	{D}^{\hyb}_{\pi_\varnothing} \coloneqq  \inn D_{\varnothing} = B^\ast.
\]
As a set, we introduce the \emph{hybrid space} $B^{\hyb}$ as the disjoint union
\begin{equation} \label{eq:Bhyb}
B^{\hyb} \coloneqq  \bigsqcup_{\substack{  F \subseteq E }} \bigsqcup_{\substack{ \pi \in \Pi(F)}} {D}^{\hyb}_\pi = B^\ast \sqcup \bigsqcup_{\substack{ \varnothing \subsetneq F \subseteq E }} \bigsqcup_{\substack{ \pi \in \Pi(F)}} {D}^{\hyb}_\pi.
\end{equation}
In the sequel, elements of $B^{\hyb}$ will be written as pairs $\thy = (t,x)$, meaning that $t \in \inn D_\pi$ and $x \in \inn \sigma_\pi$ for some ordered partition $\pi \in \Pi$.

In the last section of this paper, we will also consider a specific type of subspaces of $B^\hyb$. Namely, fix two ordered partitions $\pi$ and $\pi'$ in $\Pi$ satisfying $\pi' \preceq \pi$ (w.r.t. to the partial order "$\preceq$" on $\Pi$ introduced in Section~\ref{subsec:OrdPart}). Then, we denote by $B^{\hyb}_{[\pi, \pi']}$ the subspace 
\begin{equation} \label{eq:BInterval}
B^{\hyb}_{[\pi, \pi']} \coloneqq  \bigsqcup_{\substack{\varrho \in \Pi \\ \pi'\preceq \varrho \preceq \pi } } {D}^{\hyb}_\varrho \subseteq B^\hyb 
\end{equation}
corresponding to the closed interval $[\pi',  \pi] = \{ \varrho \in \Pi \colon \pi'\preceq \varrho \preceq \pi\}$ in $\Pi$.


\subsection{The hybrid topology} \label{ss:hybtop}

It remains to equip $B^{\hyb}$ with a suitable topology. Motivated by \eqref{eq:expalpha} -- \eqref{eq:expomega}, this will be done by using normalized logarithmic coordinates. Let $(W, z)$ be an adapted coordinate neighborhood for some point $t \in B$ (see Section \ref{ss:setup}). For a subset $F \subseteq E_t$, consider the following $\Log$-map 
\begin{equation}
\begin{array}{cccc}
\Log_F \colon & {W}\setminus \bigcup_{e \in F} D_e &\longrightarrow & \ssub{\inn \sigma}!_F \\
& s &\mapsto & \Big (\frac{ \log|z_e(s)|}{\sum_{\hat e\in F}\log |z_{\hat e}(s)|}\Big)_{e \in F}.
\end{array}
\end{equation}
We will define a topology on $B^{\hyb}$ by specifying for each point $\thy \in B^{\hyb}$ a system of neighborhoods $\mathcal{U}(\thy)$. It suffices to define a neighborhood base $\bigl(\,U^\varepsilon(\thy)\,\bigr)_{\varepsilon >0}$ around each $\thy \in B^{\hyb}$. The total system of neighborhoods $\mathcal{U} (\thy)$ will then consist of all sets $V$ with $V \supseteq U^\varepsilon(\thy)$ for some $\varepsilon > 0$. 

Hence, assume that $\thy \in B^{\hyb}$ and $\varepsilon > 0$ are fixed. Suppose that  $\thy = (t,x) \in {D}^{\hyb}_\pi$ for the ordered partition $\pi$ and fix some adapted coordinate neighborhood $(W,z)$ of $t$. Bearing in mind the decomposition \eqref{eq:Bhyb}, we will define $U^\varepsilon(\thy)$ by a similar decomposition. More precisely, let 
\[
	U^\varepsilon (\thy) \coloneqq  \bigcup_{\pi' \preceq \pi} U_{\pi'}^\varepsilon,
\]
where the subset $U_{\pi'}^\varepsilon \subseteq {D}^{\hyb}_{\pi'}$ of the hybrid stratum ${D}^{\hyb}_{\pi'}$ is defined as follows: \\
Assume that $\pi' = (\pi_i')_{i=1}^r$ and that $\pi' \preceq \pi$. Recall that this implies that $E_{\pi'} \subseteq E_\pi$ and that $\pi$ has the following form (see Section \ref{subsec:OrdPart})
\begin{equation} \label{eq:preceqdecomp}
\pi = \Big( \varrho^1, \varrho^2, \dots, \varrho^r, \varrho \Big ) = \Big( (\varrho^i)_{i=1}^r, \varrho \Big ),
\end{equation}
where $\varrho^i =(\varrho_k^i)_{k=1}^{r_i}$ is an ordered partition of $\pi_i'$ for $i=1,\dots,r$ and the last part $\varrho =(\varrho_k)_{k=1}^{r_\infty}$ is an ordered partition of $E_{\pi} \setminus E_{\pi'}$. Note that we allow $E_{\pi} = E_{\pi'}$.

A point $\mathbf{s} = (s,y)$ in the stratum ${D}^{\hyb}_{\pi'}$ belongs to $U_{\pi'}^\varepsilon$ if the following set of conditions is satisfied:

\begin{itemize}[leftmargin = 2em]
\item [(i)] For each complex coordinate $z_i$ on $W$, that is, for $i=1,\dots,\ndim$, 
\[
	 |z_i(s) - z_i(t)| < \epsilon.
\]

\item [(ii)] Assume that $E_{\pi'} \subsetneq E_\pi$, that is, $\varrho =(\varrho_k)_{k=1}^{r_\infty}$ is non-trivial. Then, 

\[
\max_{k =1, \dots, r_\infty-1} \left \{
\frac{\sum_{e\in\varrho_{k+1}}\log|z_e(s)|}{\sum_{e\in\varrho_k}\log|z_e(s)|}\right\} < \varepsilon \quad \text{and} \quad \max_{k =1, \dots, r_\infty}  \Big\{
\|\Log_{\varrho_k}(s) -x_{\varrho_k} \|_\infty \Big \} < \varepsilon,
\]
where $x_{\varrho_k}$ is the point $x_{\varrho_k} = (x_e)_{e \in \varrho_k}$ in the simplex $\inn \sigma_{\varrho_k}$.

\item[(iii)] For each $i \in \{1,\dots,r\}$, consider the ordered partition $\varrho^i =(\varrho_k^i)_{k=1}^{r_i}$ of $\pi_i'$. Then,
\[
	\max_{k =1, \dots, r_i-1} \left\{ \frac{\sum_{e \in \varrho_{k+1}^i} y_e}{ \sum_{e \in \varrho_{k}^i } y_e} \right \}  < \varepsilon \qquad  \text{ and }  \qquad \max_{k =1, \dots, r_i} \max_{e \in \varrho_k^i} \Big|\frac{y_e}{\sum_{\hat e \in \varrho_k^i} y_{\hat e}} -x_e \Big | < \varepsilon.
\]
\end{itemize}

Notice that by Proposition~\ref{prop:adapcord}, the neighborhood system $\mathcal{U}(\thy)$ is independent of the choice of the adapted coordinate neighborhood $(W,z)$ (although the $U^\varepsilon(\thy)$'s are not). 
\begin{thm} 
There is a unique topology on $B^{\hyb}$ such that $\mathcal{U}(\thy)$ coincides with the system of neighborhoods for any  $\thy \in B^{\hyb}$. 
\end{thm}
\begin{proof} It suffices to verify that the set systems $\mathcal{U}(\thy)$, $\thy \in B^{\hyb}$, satisfy the axioms of a neighborhood system and this is straightforward.
\end{proof}
Although the definition of the topology on $B^{\hyb}$ is slightly involved, there is a simple description of the convergence of sequences.

\begin{prop} \label{prop:convseq}
Let $\thy \in B^{\hyb}$ and assume $\thy=(t,x) \in D_\pi^{\hyb}$ for the ordered partition $\pi \in \Pi$. Suppose that $(\thy_n)_n = (t_n, x_n)_n$ is a sequence in $B^{\hyb}$. Then, the following statements hold true:
\begin{itemize}[leftmargin = 2em]
\item [(a)] If $\thy_n$ converges to $\thy$ in $B^{\hyb}$, then, almost all $\thy_n$ belong to hybrid strata $D_{\pi'}^{\hyb}$ of ordered partitions $\pi'$ with $\pi' \preceq \pi$.

\item [(b)] Assume that $(\thy_n)_n \subseteq D_{\pi'}^{\hyb}$ for some fixed ordered partition $\pi' \preceq \pi$. In particular, $\pi$ is of the form \eqref{eq:preceqdecomp}. Then, $\thy_n$ converges to $\thy$ in $B^{\hyb}$ if and only if the following conditions hold:
\end{itemize}
\begin{itemize}
\item [(i)] $t_n$ converges to $t$ in $B$.
\item [(ii)] Assume that $E_{\pi'} \subsetneq E_\pi$, that is, $\varrho = (\varrho_k)_k$ is non-trivial (see \eqref{eq:preceqdecomp}). Then,
\begin{equation}  \label{eq:conv1} \lim_{n\to\infty}\Log_{\varrho_k}(t_n)=x_{\varrho_k}
\end{equation}
for all $k$, where $x_{\varrho_k} = (x_e)_{e \in \varrho_k} \in \inn \sigma_{\varrho_k}$.
 Moreover, for coordinates $e \in \varrho_k$ and $e' \in \varrho_{k'}$ with $k < k'$,
\begin{equation} \label{eq:conv2}
\lim_{n\to\infty}\frac{\log|z_{e'}(t_n)|}{\log|z_e(t_n)|}=0.
\end{equation}
In the above statements, $(W,z)$ is a fixed (or equivalently, any) adapted coordinate neighborhood of $t$.

\item [(iii)] For every $i = 1,\dots,r$, let $\varrho^i =(\varrho^i_k)_k$ be the ordered partition of $\pi_i'$ in \eqref{eq:preceqdecomp}. Then, for all $k$ and all $e \in \varrho_k^i$, we have
\begin{equation} \label{eq:conv3}
\lim_{n\to\infty} \frac{x_{n,e}}{\sum_{\hat e \in \varrho_k^i} x_{n, \hat e}} =x_e.
\end{equation}
Moreover, if $e \in \varrho_k^i$ and $e' \in \varrho_{k'}^i$ with $k < k'$, then,
\begin{equation} \label{eq:conv4}
\lim_{n\to\infty}\frac{x_{n,e'}}{x_{n,e}}=0.
\end{equation}

\end{itemize}
 \end{prop}
\begin{proof}
The claims follow from the definition of the topology on $B^{\hyb}$.
\end{proof}

 \begin{remark}
 Proposition \ref{prop:convseq} confirms that the topology on $B^{\hyb}$ formalizes the ideas of \eqref{eq:expalpha}--\eqref{eq:expomega}. Indeed, \eqref{eq:conv1} and \eqref{eq:conv2} correspond exactly to the ordering of the singular coordinates $z_e$, $e\in E_t$, around a point $t$ in groups (according to their logarithmic growth) and then renormalizing. The other conditions \eqref{eq:conv3} and \eqref{eq:conv4} represent the same idea for the simplicial coordinates.
 \end{remark}
 The following proposition summarizes basic properties of the topology. 
 
 \begin{prop} \label{prop:TopologicalPropertiesHybridSpaces} $B^{\hyb}$ has the following topological properties.
 \begin{itemize}[leftmargin = 2em]
  \item [(i)] $B^{\hyb}$ is a second countable, locally compact Hausdorff space and in particular, $B^{\hyb}$ is metrizable. The convergent sequences are described in Proposition \ref{prop:convseq}. Moreover, if $B$ is compact, then $B^\hyb$ is compact.   
  
   \item [(ii)] The natural projection map $\operatorname{pr}\colon B^{\hyb} \to B$, $\thy = (t,x) \mapsto t$ is continuous. Moreover, $B$ is homeomorphic to the quotient space $B^{\hyb}/\sim$, where $\thy \sim \mathbf{s}$ for $\thy, \mathbf{s} \in B^{\hyb}$ if and only if $\operatorname{pr}(\thy) = \operatorname{pr}(\mathbf{s})$.

 \item [(iii)] The subset $B^\ast \subseteq B^{\hyb}$ is open in $B^{\hyb}$. Moreover, the topology on $B^\ast$ induced by $B^{\hyb}$ coincides with the one from $B$.

 \item [(iv)] In general, the hybrid stratum ${D}^{\hyb}_{\pi}$ of an ordered partition $\pi\in \Pi$ is neither open nor closed in $B^{\hyb}$. On each ${D}^{\hyb}_{\pi}$, the topology induced by $B^{\hyb}$ coincides with the product topology induced from~\eqref{eq:hybridstratum}.
 \end{itemize}
 \end{prop}
\begin{proof}
We prove in detail that $B^{\hyb}$ is locally compact. All other claims are elementary to verify. Fix a point $\mathbf{s} = (s,y)$ in $B^\hyb$. We have to prove that $\mathbf{s}$ has a compact neighborhood. Let $A \subset B$ be a compact neighborhood of $s = \operatorname{pr}(\mathbf{s}) \in B$. Since the projection map $\operatorname{pr}\colon B^{\hyb} \to B$ is continuous (see $(ii)$), the preimage $A^{\hyb} \coloneqq \operatorname{pr}^{-1}(A)$ is a neighborhood of $\mathbf{s}$ in $B^\hyb$. In the following, we prove that $A^\hyb \subset B^\hyb$ is compact.

Since $B^\hyb$ is second countable (see $(i)$), it suffices to show that every sequence $(\thy_n)_n$ in $A^{\hyb}$ has a convergent subsequence. Restricting to a subsequence, we may assume that all $\thy_n$, $n \in\N$, lie in the hybrid stratum ${D}^{\hyb}_{\pi'}$ of some ordered partition $\pi' = (\pi_1',\dots, \pi_r')$ of a subset $F' = E_{\pi'}$ of $E$. We write $\thy_n = (t_n,x_n)$ with $t_n =\operatorname{pr}(\thy_n) \in B$ and $x_n = (x_{n,e})_{e \in F'}  \in\inn  \sigma_{\pi'}$, with $\inn\sigma_{\pi'} = \inn\sigma_{\pi'_1} \times \inn \sigma_{\pi'_2} \dots  \times \inn \sigma_{\pi'_r}$. By compactness of $A \subset B$, we may assume that the points $t_n = \operatorname{pr}(\thy_n)$ converge in $B$ to a point $t \in A$ as $n \to \infty$. Set $F \coloneqq E_t$.

In the following, we construct an ordered partition $\pi$ of $F$ and a point $x \in \inn \sigma_{\pi}$ such that, up to passing to a subsequence, $\thy_n$ converge to the hybrid point $\thy \coloneqq  (t,x) \in {D}^{\hyb}_{\pi}$.

We begin with the construction of the ordered partition $\pi$. Using our convention from Section~\ref{subsec:OrdPart}, we will define $\pi$ of the form
\begin{equation} \label{eq:FormOrderedPartition}
\pi = \big( \varrho^1, \varrho^2, \dots, \varrho^r, \varrho \big ),
\end{equation}
where $\varrho^i = (\varrho^i_k)_k$ is an ordered partition of $\pi_i'$, $i=1, \dots,r$, and $\varrho = (\varrho_k)_k$ is an ordered partition of $F^\infty \coloneqq  F \setminus F'$, so that $\pi$ refines $\pi'$.

Fix an adapted coordinate neighborhood $(U,z)$ around $t$ with coordinates $z_e$, $e \in F$, corresponding to the components $D_e$, $e \in F$. Consider the quotients
\[
q_{e,f}(n) = \frac{|\log|z_e(t_n)||}{|\log|z_f(t_n)||}, \qquad e,f \in F^\infty,
\]
and for every $\pi_i'$, $i =1, \dots, r$, the quotients
\[
q_{e,f}(n) = \frac{x_{n,e}}{x_{n,f}}, \qquad e,f \in \pi_i'.
\]
Taking a subsequence of $(\thy_n)_n$, we assume that all limits $q_{e,f} \coloneqq  \lim_{n \to \infty} q_{e,f}(n)$ exist in $[0,+ \infty]$.

We now define the desired ordered partition $\varrho$ of $F^\infty$ as follows. If $F^\infty = \varnothing$, we set $\varrho = \varnothing$. Otherwise, we define the subset $\varrho_1 \subset F^\infty$ by
\[
\varrho_1 \coloneqq \left\{ e \in F^\infty \, \st \, q_{e,f} \neq 0 \text{ for all } f \in F^\infty\right\}.
\]
Note that $\varrho_1 $ is non-empty. Indeed, arguing by contradiction, suppose that for every $e \in F^\infty$, there exists $e' \in F^\infty$ with $q_{e,e'} =0$. Then, we get a contradiction by
\begin{align*}
1 = \frac{\sum_{e \in F^\infty} |\log|z_e(t_n)||}{\sum_{e \in F^\infty} |\log|z_e(t_n)||} = \frac{ \sum_{e \in F^\infty}  q_{ e, e'}(n) |\log|z_{e'}(t_n)||}{\sum_{e \in F^\infty} |\log|z_e(t_n)||} \le  \sum_{e \in F^\infty} q_{ e, e'}(n) \to 0 
\end{align*}
for $n \to \infty$. We then proceed by induction and repeat the same construction for the set $F^\infty \setminus \varrho_1$, and so on. More precisely, setting
	\[
	\varrho_{k+1} \coloneqq \Big \{ e \in F^\infty \setminus \bigcup_{l=1}^{k} \varrho_l \, \st \, q_{e,f} \neq 0 \text{ for all $f \in F^\infty \setminus \bigcup_{l=1}^{k} \varrho_l$} \Big \}
	\]
we inductively obtain an ordered partition $\varrho = (\varrho_k)_{k=1}^{r_\infty}$ of $F^\infty$ of some length $r_\infty \in \{1, \dots |F^\infty|\}$. Moreover, using the relations $q_{e,f}(n) = 1/q_{f,e}(n)$ and $q_{e,f}(n) = q_{e,g}(n)q_{g,f}(n)$, one can deduce that for all sets $\varrho_k$, $k = 1, \dots, r_\infty$, we have
\begin{equation} \label{eq:PropertiesPartition}
	q_{ e,f} \in (0, +\infty) \text{ for $e,f \in \varrho_k$} \quad \text{and} \quad q_{ e,f} = + \infty \text{ for $e \in \varrho_k$ and $f \in \varrho_l$ with $k < l$}.
\end{equation}

The desired ordered partition $\varrho^i$ of the set $\pi_i'$, $i=1, \dots, r$, is constructed analogously. We set
\begin{align*}
\varrho^i_1 &\coloneqq \bigl \{ e \in \pi_i' \, \st \, q_{e,f} \neq 0 \text{ for all } f \in \pi_i' \bigr \}, \quad \textrm{ and inductively define}\\
\varrho^i_{k+1} &= \Big \{ e \in \pi_{i}'\setminus \bigcup_{l=1}^{k} \varrho^i_l \, \st \, q_{e,f} \neq 0 \text{ for all } f \in \pi_i' \setminus \bigcup_{l=1}^{k} \varrho^i_l \Big \},
\end{align*}
thereby obtaining an ordered partition $\varrho^i = (\varrho^i_k)_{k=1}^{r_i}$ of $\pi_i'$ of some length $r_i \in \{1, \dots |\pi_i'|\}$. Moreover, similar to the case of $\varrho$ treated above, we have for every set $\varrho^i_k$, $k = 1, \dots, r_i$, that
\begin{equation} \label{eq:PropertiesPartition2}
	q_{e,f} \in (0, +\infty) \text{ for $e,f \in \varrho^i_k$} \quad \text{and} \quad q_{e,f} = + \infty \text{ for $e \in \varrho^i_k$ and $f \in \varrho^i_l$ with $k < l$}.
\end{equation}
Altogether, we get the desired ordered partition $\pi$ of $F$ by \eqref{eq:FormOrderedPartition}.

It remains to construct the limiting point $\thy \in D_\pi^\hyb$ of the sequence $(\thy_n)_n$. We first define an element $x =  (x_e)_{e \in F}$ of $\inn \sigma_\pi$ as follows. If $e \in \varrho_k$ for some $k$, then,
\[
	x_e \coloneqq  \lim_{n \to \infty} \frac{|\log|z_e(t_n)||}{\sum_{f \in \varrho_k} |\log|z_f(t_n)||} = (\sum_{f \in \varrho_k} q_{ f,e})^{-1}.
\]
If $e \in \varrho^i_k$ for some $i=1, \dots, r$ and $k=1, \dots, r_i$, then
\[
x_e \coloneqq  \lim_{n \to \infty} \frac{x_{n,e}}{\sum_{f \in \varrho^i_k} x_{n,f}} = (\sum_{f \in \varrho^i_k} q_{f,e})^{-1}.
\]
Finally, consider the point $\thy \coloneqq  (t, x) \in  {D}^{\hyb}_{\pi}$. Combining properties \eqref{eq:PropertiesPartition} and \eqref{eq:PropertiesPartition2} with Proposition~\ref{prop:convseq}, we conclude that $\lim_{n \to \infty} \thy_n = \thy$ in $B^\hyb$. The proof is complete.
\end{proof}

 \subsection{Family of hybrid compactifications}\label{sec:relations}
 Our constructions above refine and combine previous constructions of hybrid spaces. It leads moreover to a family of (relative) hybrid compactifications. Since this is thoroughly discussed in our work~\cite{AN-hybrid-green}, we only briefly explain it here. 
For an integer $r \in \N$, we define the following hybrid space $B^{^{\Hyb(r)}}$. Let $\Pi_r$ be the set of all ordered partitions of subsets of $E$ which have rank bounded by $r$. For each subset $F \subseteq E$ and ordered partition $\pi = (\pi_1, \dots, \pi_k)$ of $F$ with $k\leq r$, define the stratum
\[D_F^{^{\Hyb(r)}} \coloneqq \begin{cases}
\inn D_F \times \inn \sigma_{\pi_1} \times \dots \times \inn \sigma_{\pi_{r-1}} \times \sigma_{\pi_r} & \quad \textrm{if } k=r  \\
\inn D_F \times \inn \sigma_{\pi_1} \times \dots \times \inn \sigma_{\pi_{k-1}} \times \inn \sigma_{\pi_k} &  \quad \textrm{if } k <r.
\end{cases}\] 
Note that if $k=r$, then we allow coordinates in the last simplex $\sigma_{\pi_r}$ to take value zero. For $\pi=\pi_\varnothing$, we set $D_{\pi_\varnothing}^{^{\Hyb(r)}} = B^*$.

We introduce $B^{^{\Hyb(r)}}$ as the disjoint union
\begin{equation} \label{eq:Bhyb(r)}
B^{^{\Hyb(r)}} \coloneqq  \bigsqcup_{\substack{  F \subseteq E }} \bigsqcup_{\substack{ \pi \in \Pi_r(F)}} {D}^{^{\Hyb(r)}}_\pi = B^\ast \sqcup \bigsqcup_{\substack{ \varnothing \subsetneq F \subseteq E }} \bigsqcup_{\substack{ \pi \in \Pi_r(F)}} {D}^{^{\Hyb(r)}}_\pi
\end{equation} 
and define the hybrid topology on $B^{^{\Hyb(r)}}$ similarly as on $B^{\hyb}$. Analogous to the proof of Proposition~\ref{prop:TopologicalPropertiesHybridSpaces}, one can show that $B^{^{\Hyb(r)}}$ is a locally compact Hausdorff space. Moreover, if $B$ is compact, then $B^{^{\Hyb(r)}}$ is compact as well. 

Note that by definition, $B^\hyb = B^{^{\Hyb(|E|)}}$. Also, we remark that the hybrid boundary $D^{^{\Hyb(1)}}\coloneqq  B^{^{\Hyb(1)}} \setminus B^*$ of $ B^{^{\Hyb(1)}}$ coincides with the \emph{metrized complex} $\mc(D)$ associated to $D$ defined as follows. We first take the disjoint union of sets $D_F \times \sigma!_F$ for all non-empty subsets $F \subseteq E$. For each pair of non-empty subsets $F_1 \subseteq F_2$, we then identify the two sets $D_{F_1} \times \sigma!_{F_1}$ and $D_{F_2} \times \sigma!_{F_2}$ along the common subset  $D_{F_2} \times \sigma!_{F_1}$. The resulting space endowed with the quotient topology is the metrized complex $\mc(D)$. 

We further define $B^{^{\Hyb(0)}}$ as the hybrid space from~\cite{BJ17}, that is, $B^{^{\Hyb(0)}} \coloneqq B^* \sqcup \Sigma(D)$ where $\Sigma(D)$, the dual complex of $D$, is the simplicial complex that for every non-empty subset $F\subset E$ and each connected component $Z \subseteq D_F$, contains one copy $\sigma!_Z$ of the standard simplex $\sigma!_F$.

For each $j>i$, we get a forgetful map $\forget_{j>i}: B^{^{\Hyb(j)}} \to B^{^{\Hyb(i)}}$. For $i\neq 0$, this is defined as follows. The map $\forget_{j>i}$ restricts to identity on any stratum $D_\pi^{^{\Hyb(j)}}$ with $\pi \in \Pi_i$. For an element $\pi=(\pi_1, \dots, \pi_k) \in \Pi_k$ with $k\geq i+1$,  let $\pi' \coloneqq  (\pi_1, \dots, \pi_{i-1}, \pi_i \cup \dots \cup \pi_k)$. Then, the restriction of $\forget_{j>i}$ to $D_\pi^{^{\Hyb(j)}}$ has image in the stratum $D_{\pi'}^{^{\Hyb(i)}}$,  
\[\forget_{j>i} \colon D_\pi^{^{\Hyb(j)}} \to D_{\pi'}^{^{\Hyb(i)}}. \]
It is given by sending a point $\thy = (t,x=(x^1, \dots, x^{k}))$ of $D_\pi^{^{\Hyb(j)}}$ to the point $(t, y)$ with $y=(x^1, \dots, x^{i-1}, y^i=(x^i, 0, \dots, 0))$. For $i =0$, let $\forget_{j>0} \coloneqq  \forget_{1>0} \circ \forget_{j>1}$ where $\forget_{1>0}$ is defined by contracting the proper complex strata. 
\\
 This leads to a tower of hybrid spaces and continuous maps
\begin{align}\label{eq:tower}
  B^{^{\Hyb(0)}} \longleftarrow B^{^{\Hyb(1)}}  \longleftarrow \dots \longleftarrow B^{^{\Hyb(\abs E - 1)}} \longleftarrow  B^{^{\Hyb(\abs E) }} = B^\hyb
\end{align}
interpolating between $B^\hyb$ and $B^{^{\Hyb(0)}}$.


 \section{Hybrid deformation spaces and moduli of hybrid curves}\label{sec:hybrid_moduli}

The aim of this section is to construct the moduli spaces of (marked) hybrid curves $\mg_g^\hyb$ and $\mg_{g,n}^\hyb$ as well as the universal hybrid curves $\unicurve^\hyb_g$ and $\unicurve^\hyb_{g,n}$ lying above them. Although we consider coarse moduli spaces here, in proving our results, we reduce to the hybrid replacements of versal deformation spaces, which provide \'etale charts for the fine moduli space. Universal hybrid curves are defined on these charts. This approach is enough for our main purpose in this paper, that concerns the study of variations of canonical measures, as well as in our forthcoming work that concern asymptotic of analytic and complex geometry on Riemann surfaces close to the boundary of their moduli spaces.  We hope to elaborate on hybrid Deligne--Mumford stacks in our future work, in the context of providing a valuation theoretic approach to hybrid geometry.

\smallskip

 We proceed by applying the construction of Section~\ref{sec:hybrid_spaces} to the deformation space of a stable curve (that provides local \'etale charts in  $\mgbar_{g,n}$ around the corresponding point). The resulting local hybrid charts are used in order to define the hybrid topology. In order to define $\mg_g^\hyb$, set theoretically, we define hybrid strata associated to stable graphs and then take their disjoint union. Finally, we define the universal hybrid curve and describe its topology.

\subsection{Deformations of stable curves} \label{sec:deformations}

This section recalls standard results about deformations and degenerations of complex stable curves, with or without markings. These will be used in the construction of the hybrid moduli spaces as well as in the study of the monodromy and period maps (see Section~\ref{sec:monodromy}). We follow closely our main reference~\cite{ABBF} and refer to~\cite{ACGH, DM69, Hof84} for more details.

Let  $S_0$ be a stable curve with dual graph $G=(V, E, \genusfunction)$. For $v\in V$, we denote by $C_v$ the (normalization of the) corresponding irreducible component of $S_0$. Let $\genusfunction: V \to \ssZ_{\geq0}$ be  the \emph{genus function} which associates to each vertex $v$ of $G$ the genus $\genusfunction(v)$ of $C_v$. We denote by $g$ the arithmetic genus of $S_0$, which is equal to 
\[g = \graphgenus+ \sum_{v\in V} \genusfunction(v).\]
Here $\graphgenus$ denotes the first Betti number of $G$, i.e., $\graphgenus = \rk H_1(G)$. Each edge $e =uv$ corresponds to a point of intersection of $C_u$ and $C_v$ in $S_0$. These points are denoted by $p^e_u$ and $p^e_v$ in $C_u$ and $C_v$, respectively.  In this way, for each vertex $v$ in $G$, we get a \emph{marked Riemann surface of genus $\genusfunction(v)$ with $\deg(v)$ distinct marked points $p^e_v$} in one-to-one correspondence with half-edges $e$ of the graph incident to $v$. The collection of these marked points on $C_v$ is denoted by $\mathcal A_v$.

\subsubsection{Formal deformations}\label{sec:formal_deformations}  Let $S_0$ be a stable curve of arithmetic genus $g$. Standard results in deformation theory provide a smooth formal scheme $\widehat{B} = \text{Spf }\C[[t_1,\dotsc,t_{\ndim}]]$,  and a versal formal family of curves $\pr \colon \widehat{\rsf} \to \widehat{B}$ with fiber $\rsf_0$ over the point $0\in \widehat{B}$ isomorphic to $S_0$.

The total space $\widehat{\rsf}$ is formally smooth over $\C$ and the tangent space $T$ to $\widehat{B}$ at $0$ can be identified with the Ext group $\Ext^1(\Omega^1_{S_0}, \Ocal_{S_0})$ in the category of sheaves over $S_0$.  
Locally for the \'etale topology at a node of $S_0$,  we have
\begin{equation}\label{eq:node_local}
S_0 \simeq \Spec R, \quad \textrm{for } R=\C[x,y]/(xy).
\end{equation}
This implies  we have locally an isomorphism $\Omega^1_{S_0} \simeq Rdx\oplus Rdy/(xdy+ydx)$. Moreover, the element $xdy \in \Omega^1_{S_0}$ is killed by both $x$ and $y$. This means that $\Omega^1_{S_0}$ has a non-trivial torsion subsheaf supported at the singular points of $S_0$.  

By local-to-global Ext-spectral sequence, 
we get the short exact sequence 
\begin{equation}\label{eq:5term}
0 \to H^1(S_0, \usHom(\Omega^1_{S_0},\Ocal_{S_0})) \to T \to \Gamma(S_0, \usExt^1(\Omega^1_{S_0},\Ocal_{S_0})) \to 0.
\end{equation}

By the above local presentation of $\Omega^1_{S_0}$, the $\usExt$ sheaf on the right hand side of this equation is locally given by the extension 
\begin{align*}
0 \to R &\longrightarrow Rdx\oplus Rdy \to \Omega^1_R \to 0 \\
1 &\longmapsto xdy+ydx. 
\end{align*} 
In other words,  $\usExt^1(\Omega^1_{S_0},\Ocal_{S_0})$ is
the skyscraper sheaf which contains one copy of $\C$ supported at each
singular point.  This implies that 
\begin{equation}\label{eq:skyscraper}
\Gamma\bigl(S_0, \usExt^1(\Omega^1_{S_0},\Ocal_{S_0})\bigr) \simeq \C^E.
\end{equation}

By~\cite{DM69},  the global sections $\Gamma\bigl(S_0,
\usExt^1(\Omega^1_{S_0},\Ocal_{S_0})\bigr)$ correspond to smoothings of
the singular points of $S_0$. Deformations in which the dual graph remains the same (and so only the marked curves $(C_v,\mathcal A_v)$ deform) are represented by the subspace $H^1(S_0,  \usHom(\Omega^1_{S_0},\Ocal_{S_0}))\subset T$. 

 For any pair of non-negative integers $a$ and $n$, the moduli space $\mg_{a,n}$ of curves of genus $a$ with $n$ marked points has dimension $3a-3+n$. This leads to the following equality of  dimensions 
\begin{equation*}
\dim H^1(S_0,  \usHom(\Omega^1_{S_0},\Ocal_{S_0})) = \sum_{v\in
  V} (3\genusfunction(v) -3 + \deg(v)), \qquad \textrm{ and }
\end{equation*}
\begin{align*}
\ndim=\dim T &= \sum_{v\in V} (3\genusfunction(v) + \deg(v)-3) +|E|= 3g-3.
\end{align*}

The above analysis shows that deformations of $S_0$ which preserve the singularity at $p^e$ are given by a divisor $\widehat{D}_e \subset \widehat{B}$, and correspond to the tangent directions in $T$ that have zero as the $e$-th coordinate via \eqref{eq:5term} and \eqref{eq:skyscraper}. Denoting by $z_e \in \Ocal_{\widehat{B}}$ the local parameters around 0 for the divisors $\widehat{D}_e$, $e\in E$, the projection map $T \to \C^E \xrightarrow{pr_e} \C$ is seen to be defined by $dz_e$, and the surjective map $T \to \C^E$ in \eqref{eq:5term} is given precisely by the differentials of these local parameters. It thus follows that we have a collection of principal divisors $\widehat{D}_e\subset \widehat{B}$ indexed by the edges of $G$ which meet transversally. Moreover, the intersection of these divisors in $\widehat B$ is precisely the locus of those deformations of $S_0$ which keep the dual graph fixed.    More generally, we get the following correspondence. For any subset $F \subseteq E$, the intersection of the divisors $\widehat{D}_e$ for $e\in F$ is the locus of those deformations whose dual graph contains all the edges $F \subseteq E$. This is naturally stratified into those deformations whose dual graph is obtained by contracting  a subset of edges $A \subseteq E \setminus F$. This dual graph is denoted by $\contract{G}{A}$. In addition to local parameters $z_e$ for edges $e\in E$, we have $\ndim-|E|$ more local parameters which correspond to the deformations which preserve the dual graph.

Similarly, for $(S_0, q_1, \ldots, q_n)$  a complex stable curve of arithmetic genus $g$ with $n$ marked points, there exists a formal disc $\widehat{B}$ of dimension $\ndim=3g-3+n$ and a versal formal deformation $\pr \colon \widehat{\rsf} \to \widehat{B}$ such that the tangent space to $0 \in \widehat{B}$ is identified with $\Ext^1(\Omega^1_{S_0}, \Ocal_{S_0}(-\sum_{i=1}^n q_i))$. The fiber at $0$ of the family is isomorphic to $S_0$, and  in addition, the family comes with sections $\sigma_i \colon \widehat{B} \to \widehat{\rsf}$ such that $\sigma_i(0)=q_i$. 

The local study from above is carried out in a similar way. It leads to local parameters $z_e$, indexed by the edges of the dual graph $G = (V,E)$ of $S_0$, and formal divisors $\widehat{D}_e\subset \widehat{B}$, $e \in E$,  which meet transversally. In addition to the local parameters $z_e$, $e\in E$, we have $\ndim- |E| = 3g-3+n-|E|$ more local parameters that correspond to the deformations of the stable marked curve which preserve the dual graph. That is, they correspond to deformations of the marked curves $(C_v, \mathcal A_v)$. Here, $\mathcal A_v$ consists of the nodes $p^e_v$, $e \sim v$, and those marked points $q_1, \dots, q_n$ that lie on $C_v$.

As in the unmarked case, we get the following correspondence. For any subset $F \subseteq E$, the intersection of the divisors $\widehat{D}_e$ for $e\in F$ is the locus of those deformations whose dual graph contains all edges $F \subseteq E$. This is naturally stratified into those deformations whose dual stable graph with $n$ markings is obtained by contracting  a subset of edges $A \subseteq E \setminus F$, that is, it is given by  $\contract{G}{A}$ (see Section~\ref{sec:prel}). Recall that the marking function in $\contract{G}{A}$ is defined as the composition of the map
$\marking \colon [n] \to V(G)$ with the projection map $\proj \colon V(G) \to V(\contract{G}{A})$. We denote this marked graph still with  $\contract{G}{A}$, as it will be clear from the context whether markings are involved or not.

\subsubsection{Analytification} From the local theory above using formal schemes, we get analytic deformations $\rsf \to B$ of stable Riemann surfaces (with markings if $n\neq 0$) over a polydisc $B$ of dimension $\ndim=3g-3 +n$. 
This means that
\[B = \underbrace{\Delta \times \Delta \times \dots \times \Delta}_{\ndim \textrm{ times }}
\]
for $\Delta$ a small disk around $0$ in $\C$.
Moreover, the formal divisors above give rise to analytic divisors $D_e \subset B$ which are defined by equations $\{z_e=0\}$ for analytic local parameters $z_e$ associated to edges $e\in E$. Once again, $D_e$ is the locus of all points $t\in B$ such that in the family $\pr \colon \rsf \to B$, the fiber $\pr^{-1}(t)$ is a Riemann surface which has a singular point corresponding to $e$. We denote by $p^e \colon D_e \to \rsf$ the associated section, mapping $t \in D_e$ to the singular point $p^e(t)$ on the fiber $\rsf_t$.

In the following, $B$ will denote the base of the versal deformation of the stable curve $S_0$, or the stable marked curve $(S_0, q_1, \dots, q_n)$. We let $B^* \coloneqq B \setminus \bigcup_{e \in E} D_e$ be the locus of points whose fibers in the family are smooth Riemann surfaces and set $\rsf^* \coloneqq  \pr^{-1}(B^*)$. 

For more details related to the discussion in this section, we refer to~\cite[Chapter XI]{ACGH} and the recent work of Hubbard-Koch~\cite{HK14} (that plays a crucial role in the further investigation of the geometry of hybrid curves and their moduli spaces in the sequel~\cite{AN-hybrid-green}).

\subsection{Hybrid deformation space}  \label{sec:hybrid_base}
We assume w.l.o.g. that  $E = \{1,\dots,|E|\} \subseteq \{1,\dots,\ndim\}$, Let  $B = \Delta^{\ndim}$ be a polydisc in $\C^{\ndim}$ and $D = \bigcup_{e \in E} D_e$, where the divisors $D_e$, $e \in E$, are given by 
\[
	D_e = \left\{ z \in  \Delta^{\ndim} \st z_e = 0 \right\}.
\]
The corresponding hybrid space $B^{\hyb}$ will serve as the base space for the hybrid versal deformation. Notice that in this setting, adapted coordinate neighborhoods $(U,z)$ can be defined in a particularly simple way: for each $t \in B$, we can choose $U = B \,\setminus \bigcup_{e: t_e \neq 0} D_e$ and $z=(z_i)_i$ as the standard coordinates on $U \subseteq \Delta^{\ndim}$.

\subsection{Universal hybrid curve over hybrid deformation space} \label{sec:versal_hybrid_curve}
In this section, we introduce the \emph{universal hybrid family} $\rsf^\hyb $ over the hybrid deformation space $B^\hyb$ associated to a stable curve  $S_0$.  The fiber $\rsf_{\thy}^\hyb$ over a point ${\thy} = (t, x) \in B^{\hyb}$ is the hybrid curve represented by $\thy$. This has the underlying stable Riemann surface $\rsf_t$ (the fiber in the original versal family $\pr\colon \rsf \to B$) and has edge length function given by the coordinates $x_e$, $e \in E_t$, of $x$. Here, $E_t$ is naturally identified with the set of edges of the dual graph of $\rsf_t$. Moreover, the ordered partition $\pi =(\pi_k)_k$ of $E_t$ is such that $\thy \in D_\pi^\hyb$. In particular, since $x \in \inn \sigma_{\pi}$, the normalization condition $\sum_{e \in \pi_k} x_e = 1$ holds. In short, we get the hybrid family over $B^\hyb$ by replacing non-smooth fibers of $\pr\colon \rsf \to B$ by hybrid curves.

The rigorous definition of the family of hybrid curves $\pr^\hyb\colon \rsf^\hyb \to B^\hyb$ is slightly involved. However, since it plays a key role in our works, we present the construction in some detail. In order to simplify the exposition, we omit any mention of the markings.

Let $S_0$ be a stable curve  with dual graph $G=(V, E, \genusfunction)$ and $\pr\colon \rsf \to B$ its versal deformation family (see Section \ref{sec:deformations}). The family $\rsf^\hyb$ will be defined over the hybrid deformation space $B^{\hyb} = \bigsqcup_\pi D_\pi^\hyb$ from Section \ref{sec:hybrid_base}. We will introduce $\rsf^\hyb$ by a decomposition
\begin{equation} \label{eq:S_hyb_decomposition}
\rsf^\hyb \coloneqq \rsf^\ast \sqcup \bigsqcup_{\substack{ \varnothing \subsetneq F \subseteq E }} \,\, \bigsqcup_{\pi \in \Pi(F)} \rsf_\pi^\hyb \end{equation}
where each $\pr^\hyb_\pi\colon \rsf_\pi^\hyb \to D_\pi^\hyb$ is a family of hybrid curves with ordered partition $\pi$ defined over $D_\pi^\hyb$. The subfamily $\rsf^\ast \coloneqq  \rsf_{\pi_\varnothing}^\hyb$ will coincide with the original family of Riemann surfaces $\rsf^\ast = \rsf\rest{B^\ast}$. Moreover, we will define a hybrid topology on $\rsf^\hyb$.

\subsubsection{Normalized families} Recall from Section \ref{sec:deformations} that the fiber $\rsf_t$ of each $t \in B$ is a stable Riemann surface with exactly $|E_t|$ singular points $p^e(t)$, $e\in E_t$. Its stable dual graph is the graph $G_{t} \coloneqq G / (E \setminus E_t)$ obtained by contracting all edges in $E \setminus E_t$. Moreover, $\rsf_t$ has precisely $|V(G_{t})|$ irreducible components $C_{v,t}$, $v \in V(G_{t})$.

The sections $p^e \colon D_e \to \rsf$, $e\in E$, are continuous. The complement of their images
\begin{equation} \label{eq:P^0}
	P^\circ \coloneqq  \rsf \, \setminus \bigcup_{e \in E} p^e(D_e) \subseteq \rsf,
\end{equation}
consists of all points in $\rsf$ which are also smooth points of their fiber.

Let $p = p^e(s)$ be a point on the section $p^e \colon D_e \to \rsf$ for some $e \in E$. Then, there exists a small neighborhood $U$ of $p$ in $\rsf$ and local coordinates (note that $\rsf$ is smooth as an analytic space)
\begin{equation} \label{eq:standard_coordinates}
	z=\Big ( (z_i)_{i \neq e}, z_u^e, z_v^e \Big)
\end{equation}
with the following properties:

\begin{itemize}[leftmargin = 2em]
\item [(i)] All coordinates have absolute value smaller than one, that is, $|z_u^e|, |z_v^e|$ and $|z_i|< 1$ on $U$ for all $i=1,\dots \ndim$, $i \neq e$, and

\item [(ii)] the coordinates are compatible with the coordinates on the base $B$ by the projection map $\pr \colon \rsf \to B$, in the sense that for all $q \in U$,
\[
	\pr(q)_i = z_i(q), \qquad i \neq e,
\]
and $\pr(q)_e = z_u^e(q) z_v^e (q)$,  and

\item [(iii)] for each base point $t \in \pr(U)$ with $t \in D_e$, we have
\[
C_{u,t} \cap U = \left\{ q\in U\,\, \st \,\,  \pr(q) = t \text{ and } z_v^e(q) = 0 \right\}
\] and similarly for $C_{v,t}$ and $z_u^e$. Here, $C_{u,t}$ and $C_{v,t}$ are the components of the fiber $\rsf_t$ corresponding to the vertices $u, v \in V(G_t)$ with $e = uv$.
\end{itemize}

In the following, we call a pair $(U,z)$ consisting of a neighborhood $U$ and local coordinates $z$ on $U$ with the above properties a \emph{standard coordinate neighborhood} of $p \in p^e(D_e)$ on $\rsf$.

Fix an ordered partition $\pi$ of a non-empty subset $F \subseteq E$. As the first step in the definition of the aforementioned family $\rsf_\pi^\hyb \to D_\pi^\hyb$ (see \eqref{eq:S_hyb_decomposition}), we construct the \emph{Riemann surface components} for the hybrid curves.

By the preceding discussion, for $t \in \inn D_F$, the fiber $\rsf_t$ is a stable Riemann surface with stable dual graph $G_F \coloneqq G / (E \setminus F) $. Moreover, $\rsf_t$ has $|V(G_F)|$ components $C_{v,t}$, $v \in V(G_F)$, and $|F|$ singular points $p^e(t)$, $e \in F$. In what follows, we identify $F$ with the edge set of the graph $G_F$.

We consider the restricted family $\rsf_F\coloneqq  \rsf\rest{\inn D_F} \to \inn D_F$ and construct the corresponding \emph{normalized family} 
\[
	\widetilde{\rsf}_F \to \inn D_F.
\]
That is, in the fibers $\rsf_{F}\rest{t} = \rsf_t$, $t \in \inn D_F$, we replace each node $p^e(t)$ of an edge $e= uv \in F$ by two different points $p_u^e(t)$ and  $p_v^e(t)$ lying on the components $C_{u,t}$ and $C_{v,t}$, thereby disconnecting the components. This amounts to replacing points $q$ with $z_u^e(q) = z_v^e(q) =0 $ by two new points, as in the normalization of a single fiber.  Each fiber $\widetilde{\rsf}_{F}\rest{t} $ of $\widetilde{\rsf}_F$ is the disjoint union of the smooth Riemann surfaces $C_{v,t}$ for $v \in V(G_F)$. For each edge $e = uv \in F$, there are two continuous sections $p^e_u, p^e_v \colon \inn D_F \to \widetilde{\rsf}_F$ such that $p^e_u(t)$ and $p^e_v(t)$ lie in the components $C_{u,t}$ and $ C_{v,t}$.  Moreover, $\rsf_F$ is obtained from $\widetilde{\rsf}_F$ by glueing the sections $p^e_u$ and $p^e_v$, for all $e= uv \in F$, over any point $t\in \inn D_F$.  The construction of $\widetilde{\rsf}_F$ clearly leaves the smooth points $p \in P^\circ$ invariant (see \eqref{eq:P^0}). More formally,
\begin{equation} \label{eq:smooth_points_normalized}
	 \left\{ q \in P^\circ \, \st \, \pr(q)  \in \inn D_F \right\}    \subseteq \widetilde{\rsf}_F.
\end{equation}

Since $D_\pi^\hyb = \inn D_F \times \inn \sigma_\pi$, we can naturally extend $\widetilde{\rsf}_F \to \inn D_F$ to
\begin{equation}  \label{eq:hybrid_versal_1}
	\widetilde{\rsf}_F \times \inn \sigma_\pi \to D_\pi^\hyb.
\end{equation}
For each edge $e = uv$ in $F$, the sections $p^e_u, p^e_v \colon \inn D_F \to \widetilde{\rsf}_F$ extend to sections
\begin{equation} \label{eq:hybrid_sections_1}
\tilde{p}^e_u, \tilde{p}^e_v \colon D_\pi^\hyb   \to \widetilde{\rsf}_F \times \inn \sigma_\pi.
\end{equation}

\subsubsection{Metric graph part of hybrid curves} \label{ss:MetricGraphsUniversalCurve}
Next, we will construct the \emph{versal family of intervals} for the aforementioned family $\rsf_\pi^\hyb \to D_\pi^\hyb$. Recall that we had fixed an ordered partition $\pi$ of a non-empty subset $F \subseteq E$. For each edge $e =uv$ in $F$, define
\begin{equation} \label{eq:Ical_pi^e}
	\Ical_\pi^e= \Big \{ (x, \lambda_u, \lambda_v) \in \inn \sigma_{\pi} \times \ssR_{\ge 0}^2\,\, \bigl| \,\, \lambda_u + \lambda_v = x_e \Big \}  \subseteq \inn \sigma_\pi \times \ssR_{\ge 0}^2.
\end{equation}
The projection to the first factor $\operatorname{pr}_1 \colon \Ical_\pi^e \to \inn \sigma_\pi$ is continuous and the fiber over $x=(x_f)_{f \in F} \in  \inn \sigma_\pi$ can be identified with the  interval $[0, x_e]$ of length $x_e >0$. Moreover, $\Ical_\pi^e \to \inn \sigma_\pi$ is a fiber bundle, topologically isomorphic to the trivial bundle $(\inn \sigma_\pi \times [0,1]) \to  \inn \sigma_\pi$. We obtain two continuous sections
\begin{align*}
s_u^e \colon & \inn \sigma_\pi  \longrightarrow  \Ical_\pi^e,   \qquad (x_f)_f\mapsto  (x, 0, x_e), \quad \textrm{and}  \\
s_v^e \colon & \inn \sigma_\pi  \longrightarrow  \Ical_\pi^e,     \qquad (x_f)_f \mapsto  (x, x_e , 0).
\end{align*}

Since $D_\pi^{\hyb} = \inn D_F \times \inn \sigma_\pi$, we can naturally extend $\Ical_\pi^e \to \inn \sigma_\pi$ to
\begin{equation} \label{eq:hybrid_versal_2}
	\inn D_F \times \Ical_\pi^e \to D_\pi^{\hyb}.
\end{equation}
For each edge $e = uv$ in $F$, the sections $s_u^e, s_v^e \colon \inn \sigma_\pi  \to \Ical_\pi^e$ extend to sections
\begin{equation} \label{eq:hybrid_sections_2}
\tilde{s}_u^e, \tilde{s}_v^e \colon D_\pi^\hyb   \to \inn D_F \times \Ical_\pi^e.
\end{equation}

\subsubsection{Construction of $\rsf^\hyb$}
In the following, we combine the previous constructions to define the families $\rsf_\pi^\hyb \to D_\pi^\hyb$, $\pi \in \Pi$, and $\rsf^\hyb \to B^\hyb$.

Fix an ordered partition $\pi$ of some subset $F \subseteq E$. From the preceding considerations, we get the \emph{normalized hybrid family} 
\begin{equation} \label{eq:HybridNormalizedFamily}
	\widetilde{\rsf}_\pi^{\, \hyb} \coloneqq  \Big( \widetilde{\rsf}_F \times \inn \sigma_\pi  \Big) \sqcup  \Big(\bigsqcup_{e \in F} \inn D_F \times \Ical_\pi^e \Big).
\end{equation}

By \eqref{eq:hybrid_versal_1}--\eqref{eq:hybrid_sections_2}, we see that $\widetilde{\rsf}_\pi^{\, \hyb}$ comes with a map
\[
	\widetilde\pr^\hyb_\pi\colon {\widetilde \rsf}_\pi^{\, \hyb} \to D_\pi^{\hyb}
\]
and continuous sections $\tilde{p}_u^e, \tilde{p}_v^e, \tilde{s}_u^e, \tilde{s}_v^e \colon D_\pi^{\hyb} \to {\widetilde \rsf}_\pi^{\, \hyb}$ for all edges $e = uv$ in $F$.

We define $\rsf_\pi^{\hyb}$ by glueing these sections pairwise, that is, as the quotient  
\begin{equation} \label{eq:definition_S_pi^hyb}
	{\rsf}_\pi^{\hyb} \coloneqq  {\widetilde \rsf}_\pi^{\, \hyb} / \sim,
\end{equation}
where over each base point $\thy \in D_\pi^{\hyb}$, we identify $\tilde{p}_u^e(\thy) $ with $\tilde{s}_u^e(\thy)$ and $\tilde{p}_v^e(\thy)$ with $\tilde{s}_v^e(\thy)$, for all $e = uv$ in $F$. The resulting points in ${\rsf}_\pi^{\hyb}$ are denoted by $\p_u^e(\thy)$ and $\p_v^e(\thy)$.

Endowing ${\rsf}_\pi^{\hyb}$ with the quotient topology, we get a continuous map 
\begin{equation} \label{eq:hybrid_projection_pi}
	\pr^\hyb_\pi\colon {\rsf}_\pi^{\hyb} \to D_\pi^{\hyb}.
\end{equation}
By construction, the fiber ${\rsf}_\pi^{\hyb}\rest{\thy}$ at some $\thy = (t,x) \in D_\pi^{\hyb}$ is the metrized complex $\mc$ defined by the stable Riemann surface $\rsf_t$ and the edge lengths of $x$ on its dual graph $G_F$. To turn ${\rsf}_\pi^{\hyb} \rest\thy$ into a hybrid curve, we add the data of the ordered partition $\pi$. Altogether, ${\rsf}_\pi^{\hyb} \to D_\pi^{\hyb} $ is a family of hybrid curves with underlying graph $G_F$ and layering $\pi$. The continuous sections $\p_u^e, \p_v^e \colon D_\pi^{\hyb} \to {\rsf}_\pi^{\hyb}$ give the attachment points of the intervals representing the edge $e = uv$ in the hybrid curves. For the empty ordered partition $\pi_\varnothing$, we simply recover the original family
\[
\rsf^\ast \coloneqq  \rsf_{\pi_\varnothing}^\hyb = \rsf\rest{B^\ast}
\] over the open part $D_{\pi_\varnothing} = B^\ast$.

Notice that $\rsf_\pi^\hyb$ consists of three different types of points:
\begin{equation} \label{eq:point_types}
 \rsf_\pi^\hyb = \Pty_\pi^0 \sqcup \Pty_\pi^1 \sqcup \Pty_\pi^2
\end{equation}
for the following three sets:
\begin{itemize}[leftmargin = 2em]
\item $\Pty_\pi^0$ is the set of all points in $\rsf_\pi^{\hyb}$, which are smooth points of some Riemann surface component of their fiber. Formally, we can write $\Pty_\pi^0$ as (see \eqref{eq:smooth_points_normalized})
\begin{equation} \label{eq:P1}
	\Pty_\pi^0 \coloneqq  \bigsqcup_{\thy=(t,x) \in D_\pi^{\hyb}} (\rsf_t \cap P^\circ)\times \{x\} = ( \rsf_F \cap P^\circ ) \times \inn \sigma_\pi \subseteq \rsf_\pi^{\hyb}.
\end{equation}
\item $\Pty_\pi^1$ consists of all attachment points of intervals in the fibers,
\begin{equation} \label{eq:P2}
	\Pty_\pi^1 \coloneqq  \bigsqcup_{\thy \in D_\pi^{\hyb}} \bigsqcup_{e =uv \in F} \Big \{ \p_u^e(\thy), \p_v^e(\thy) \Big \}.
\end{equation}
\item $\Pty_\pi^2$ is the set of all points in $\rsf_\pi^{\hyb}$, which lie in the strict interior of some interval of their fiber. That is, formally (see \eqref{eq:hybrid_versal_2}),
\begin{equation} \label{eq:P3}
	\Pty_\pi^2 \coloneqq   \Big(\bigsqcup_{e \in F} \inn D_F \times \Ical_\pi^e \Big) \setminus \Pty_\pi^1 \subseteq \rsf_\pi^\hyb.
\end{equation}
\end{itemize}

Finally, as a set, we introduce the family $\rsf^{\hyb}$ by \eqref{eq:S_hyb_decomposition}, that is, we set
\begin{equation} 
\rsf^\hyb \coloneqq  \rsf^\ast \sqcup \bigsqcup_{\substack{ \varnothing \subsetneq F \subseteq E }} \,\, \bigsqcup_{\pi \in \Pi(F)} \rsf_\pi^\hyb.
\end{equation}
Since $B^{\hyb} = \bigsqcup_{\pi} D_\pi^{\hyb}$, the maps $\pr^\hyb_\pi$ in \eqref{eq:hybrid_projection_pi} glue to a \emph{hybrid projection map}
\begin{equation} \label{eq:hybrid_projection}
	\pr^\hyb \colon {\rsf}^{\hyb} \to B^{\hyb}.
\end{equation}
Note that on $\rsf^\ast \subseteq \rsf^\hyb$, the hybrid projection $\pr^\hyb$ coincides with the projection $\pr\colon \rsf^\ast \to B^\ast$ for the versal family $\pr \colon \rsf \to B$. If there is no risk of confusion, we will shorten notation and use the same letter $\pr$ for both maps.

\subsubsection{The topology on $\rsf^\hyb$} \label{ss:TopologyHybridFamily} 
The precise definition of the topology on $\rsf^\hyb$ is quite technical, and the details are presented in Appendix~\ref{sec:AppendixTopology}. In the following, we give an informal outline and describe the convergent sequences.

Consider the fiber $\rsf_\thy^\hyb$ over a base point $\thy =(t,x) \in B^\hyb$. The topology on $\rsf^\hyb$ should ensure the continuity of the projection map $\pr \colon {\rsf}^{\hyb} \to B^{\hyb}$. Hence it suffices to clarify the topological relationship between $\rsf_\thy^\hyb$ and nearby fibers $\rsf^\hyb_\s$, that is, fibers whose base points $\s = (s, y) \in B^\hyb$ are close to $\thy$ in $B^\hyb$. By the properties of $B^\hyb$, we have the inclusion $E_s \subseteq E_t$ for the edge sets of the hybrid curves, that is, $\rsf_\thy^\hyb$ has more intervals than $\rsf_\s^\hyb$. Clearly, the Riemann surfaces parts and the  intervals of common edges $e \in E_s$ in $\rsf^\hyb_\s$ and $\rsf^\hyb_\thy$ can be related naturally. In order to treat the additional intervals in $\rsf_\thy^\hyb$, recall that $\rsf^\hyb_\s$ contains an analytic cylinder of the form 
\[W_{e,s} = \left\{(z^e_u, z^e_v) \in \C^2 \,\st \, |z^e_u|, | z^e_v| < 1, \,  z^e_u z^e_v = z_e(t)\right\}\] for every $e \in E_t \setminus E_s$. In the family $\rsf \to B$, the cylinder $W_{e,s}$ degenerates to two half-discs glued at the node $p^e(t)$, as $s \to t$. On the other hand, changing to \emph{logarithmic coordinates}
\[
z_u^e \, \,  \longleftrightarrow \, \,\Big (\arg(z_u^e), |\log(|z_u^e|)| \Big),
\]
$W_{e,s}$ becomes a real cylinder of length $|\log(|z_e(s)|)|$. Forgetting the angle and rescaling, we can relate it to the respective interval $\Ical_e$ in $\rsf^\hyb_\thy$.

These ideas are formalized by the following \emph{$\Log$-maps}. Let $p \in p^e(D_e)$ be a point on the section $p^e$, $e \in E$. Fix a standard coordinate neighborhood $(U,z)$ for $p$ in $\rsf$ (see \eqref{eq:standard_coordinates}). Then, we introduce the $\Log$-map
\begin{equation} \label{eq:DefLogMap}
\begin{array}{cccc}
\Log_e \colon & {U}\setminus \pr^{-1}(D_e) &\longrightarrow & [0,1] \\
& q &\mapsto & \frac{ \log|z_u^e(q)|}{\log|z_e(q)|}, 
\end{array}
\end{equation}
where  $z_e(q) = \pr(q)_e = z_u^e(q) z_v^e(q)$ and $\pr \colon \rsf \to B$ is the projection map. Using the Log-maps, we introduce the topology on $\rsf^\hyb$ in Appendix~\ref{sec:AppendixTopology}.

The following proposition characterizes the convergence of points on smooth Riemann surface fibers to points on hybrid curve fibers. In a similar way, one can characterize also the convergence of general sequences. 

\begin{prop} \label{prop:HybridFamilyConvergence}
Let $\p$ be a point in $\rsf^\hyb$ and $\thy =(t,x) \coloneqq  \pr(\p)$ its hybrid base point. Suppose further that $(p_n)_n$ is a sequence in $\rsf^\ast \subseteq \rsf^\hyb$ with base points $t_n \coloneqq  \pr(p_n)$ in $B^\ast \subseteq B^\hyb$. If $\lim_{n \to \infty} p_n = \p$ in $\rsf^\hyb$, then necessarily
\begin{equation} \label{eq:baseconv}
	\lim_{n \to \infty} t_n = \thy \qquad \text{in } B^\hyb.
\end{equation}
Moreover, the convergence can be characterized as follows.
\begin{itemize}[leftmargin = 2em]
\item [(i)] Suppose that $\p$ belongs to the Riemann surface part of its fiber $\rsf_\thy^\hyb$ (and is not an attachment point). That is, $\p \in \Pty_\pi^0$ and formally $\p = (p, x)$ for some non-node point $p$ in the fiber $\rsf_t$ (see ~\eqref{eq:P1}).

Then, $\lim_{n \to \infty} p_n = \p$ in $\rsf^\hyb$ exactly when~\eqref{eq:baseconv} holds and
\[
	\lim_{n \to \infty} p_n = p \qquad \text{in } \rsf.
\]
\item [(ii)] Suppose that $\p$ belongs to an interval of its fiber $\rsf_\thy^\hyb$.
That is, $\p = \lambda^e_u$ is a point $0 \le \lambda_u^e \le x_e$ in the closed interval $\Ical_e = [0,x_e] \subseteq \rsf_\thy^\hyb$ representing an edge $e=uv$ in the dual graph of $\rsf_\thy^\hyb$. The attachment points $\p_u^e(\thy)$ and $\p_v^e(\thy)$ correspond to $\lambda_u^e = 0$ and $\lambda_u^e = x_e$, respectively.

Then, $\lim_{n \to \infty} p_n = \p$ in $\rsf^\hyb$ exactly when~\eqref{eq:baseconv} holds, $\lim_{n \to \infty} p_n = p^e(t)$ in the original family $\rsf$ and
\[
	\lim_{n \to \infty} \Log_e(p_n) =  \lim_{n \to \infty} \frac{ \log|z_u^e(p_n)|}{\log|z_e(p_n)|} = \frac{\lambda_u^e}{x_e} 
\]
for the $\Log$-map on some (equivalently, on every) standard coordinate neighborhood $(U,z)$ of $p^e(t)$ in $\rsf$.

\end{itemize}
\end{prop}
\begin{proof}
The claims are easily derived from the definition of the topology on $\rsf^\hyb$, see Appendix~\ref{sec:AppendixTopology}.
\end{proof}

\begin{remark}
\begin{itemize}[leftmargin = 2em] The topological space $\rsf^{\hyb}$ has the following properties.
\item [(i)]  $\rsf^{\hyb}$ is a second countable, locally compact Hausdorff space and in particular, $\rsf^{\hyb}$ is metrizable.
\item[(ii)] The family $\pr \colon \rsf^\hyb \to B^\hyb$ of hybrid curves is continuous over $B^{\hyb}$.
\item[(iii)] On each subfamily $\rsf_{\pi}^{\hyb} \subseteq \rsf^\hyb$ (see \eqref{eq:S_hyb_decomposition}), the induced topology from $\rsf^\hyb$ coincides with the one introduced on $\rsf_{\pi}^{\hyb}$ in the previous section.
\end{itemize}
\end{remark}


\subsection{Hybrid moduli space} \label{sec:hybrid_moduli_space} Let $g$ and $n$ be two non-negative integers verifying  $3g-3+n\geq 0$. If $g=1$, we impose moreover that  $n\geq 1$. 
 
\subsubsection{Stratification of $\mgbar_{g,n}$ by stable graphs with markings} A nice reference for the materials in this short subsection is~\cite{ACP}, which already makes connections to the moduli space of tropical curves. See also~\cite{CCUW20} for a more refined treatment of the stack structure. Once again, we emphasize that our notion of tropical curve differs from the one considered in~\cite{ACP, CCUW20}.

 Let $\mg_{g,n}$ be the (coarse) moduli space of stable curves of genus $g$ with $n$ markings.  Denote by $\mggbar{g,n}$ the Deligne--Mumford compactification of $\mg_{g,n}$, obtained by adding stable curves of genus $g$ with $n$ marked points to $\mg_{g,n}$.

There is a stratification of $\mggbar{g,n}$ by combinatorial types of stable dual graphs $G$ (with $n$ marked points) of genus $g$. Denoting by $\mg!_G$ the associated stratum, we have
\[\mggbar{g,n} =\bigsqcup_{\substack{G:\, \textrm{stable marked graph of genus $g$}}} \mg!_G, \]
where the disjoint union is taken over isomorphism classes of stabled marked graphs of genus $g$. Here, an isomorphism between stabled marked graphs of genus $g$ is an isomorphism between the underlying graphs that additionally respects the genus functions and the marking functions from $[n]$ to the sets of vertices.

The stratum $\mg!_G$ can be described as follows. Consider a stable graph with $n$ marked points $G = (V, E, \genusfunction, \marking)$.  Let $\countmarking \colon V \to \N \cup \{0\}$ be the counting function that maps $v \in V$ to the number of labels placed at $v$.  For each vertex $v \in V$, consider the moduli space $\mg_{\genusfunction(v), \, \deg(v)+ \countmarking(v)}$, and the product
\begin{equation} \label{eq:ModuliStratumProduct} \ssub{\widetilde{\mg}}!_{G} \coloneqq  \prod_{v\in V} \mg_{\genusfunction(v), \, \deg(v)+ \countmarking(v)}
\end{equation}
Denote by $\aut(G)$ the automorphism group of the marked stable graph $G$, that is, the set of automorphisms of the underlying graph that respect both the genus  and the marking functions. The group $\aut(G)$ naturally acts on $\ssub{\widetilde{\mg}}!_{G}$ through the decomposition as a product of factors, by permuting the factors and respecting the marked points. We have $\mg!_G \simeq \ssub{\widetilde \mg}!_G /\aut(G)$.

We denote by $\unicurve_{g,n}$ the universal family of curves with $n$ markings over $\mg_{g,n}$. By an abuse of the notation, the universal family of stable curves with $n$ markings over $\mgbar_{g,n}$ is also denoted by $\unicurve_{g,n}$. We understand this in the sense of the theory of Deligne-Mumford stacks: in our setting, on étale charts given by the versal deformation spaces, we get a well-defined family of stable curves. When passing to the coarse moduli space, the fibers of the family are quotients of stable Riemann surfaces by the action of their automorphism groups. We refer to~\cite[Chapter XI]{ACGH} and~\cite{HK14} for more details.  

\subsubsection{Moduli space of tropical curves of given combinatorial type}\label{sec:moduli_space_tropical_curves}

Let $G = (V, E, \genusfunction, \marking)$ be a stable graph of genus $g$ with $n$ marked points. 
We define the \emph{moduli space of tropical curves of combinatorial type $G$} as follows. First, we introduce 
\[\ssub{\widetilde \mg}!_G^{\,\,\,\trop} \coloneqq  \bigsqcup_{\pi  \in \Pi(E)} \inn \sigma_{\pi},\]
where, we recall, for $\pi = (\pi_1, \dots, \pi_r)$ an ordered partition of $E$, we set $\inn \sigma_{\pi}  \coloneqq  \inn \sigma_{\pi_1} \times \inn \sigma_{\pi_2} \dots  \times \inn \sigma_{\pi_r}$. 

In Section~\ref{sec:tropical}, we endow $\ssub{\widetilde \mg}!_G^{\,\,\,\trop}$ with a natural topology, compatible with the hybrid topology on $B^\hyb$ from Section~\ref{ss:hybtop}. Here the vertices and the edges of the graph $G$ are all labeled.
The group $\aut(G)$ naturally acts on $\ssub{\widetilde \mg}!_G^{\,\,\,\trop}$ by permutation of the edges, and the moduli space of tropical curves of combinatorial type $G$ is defined by taking the quotient 
\begin{equation} \label{eq:ModuliSpaceFixedGraph}
\mg!_G^\trop \coloneqq  \ssub{\widetilde \mg}!_G^{\,\,\,\trop} /\aut(G).
\end{equation}

\subsubsection{Hybrid moduli space}

Let $G = (V, E, \genusfunction, \marking)$ be a stable graph of genus $g$ with $n$ marked points. The group $\aut(G)$ acts naturally on $\ssub{\widetilde \mg}!_G$ and on $\ssub{\widetilde \mg}!_G^{\,\,\,\trop}$.  From these two actions, we get a diagonal action of $\aut(G)$ on the product $\ssub{\widetilde \mg}!_G\times \ssub{\widetilde\mg}!_G^{\,\,\,\trop}$. We define the \emph{hybrid stratum} $\mg!_G^\hyb$ associated to $G$ as 
\[\mg!_G^\hyb \coloneqq  \contract{(\ssub{\widetilde \mg}!_G \times \ssub{\widetilde\mg}!_G^{\,\,\,\trop})}{\aut(G)}.\]

 By construction, the elements of the hybrid stratum $\mg!_G^\hyb$ correspond to isomorphism classes of hybrid curves with combinatorial type $G$ (here, we consider isomorphisms of the underlying marked stable Riemann surfaces respecting the ordered partition and edge lengths).

We now give an alternative characterization of these hybrid strata and obtain a refined decomposition into hybrid orbifolds.

Let $G = (V, E, \genusfunction, \marking)$ be a stable dual graph with $n$ markings. Consider an ordered partition $\pi = (\pi_1, \dots, \pi_r)$ of $E$ and denote by $(G, \pi)$ the corresponding layered stable graph with markings. Denote by $\aut(G, \pi)$ the subgroup of $\aut(G)$ consisting of those elements which preserve the ordered partition $\pi$, i.e., which respect the corresponding filtration $\filter^\pi$.  The hybrid orbifold associated with $\pi$ is the quotient
\[\mg!_{G, \pi}^\hyb \coloneqq \contract{(\ssub{\widetilde \mg}!_G  \times \inn \sigma_{\pi}) }{\aut(G, \pi)}. \]
Notice that the hybrid orbifold $\mg!_{G, \pi}^\hyb$ only depends on the orbit $[\pi]$ of the ordered partition $\pi \in \Pi(E)$ under the action of $\aut(G)$, that is, $\mg!_{G, \pi}^\hyb$ and $\mg!_{G, \pi'}^\hyb$ are isomorphic whenever $[\pi] =  [\pi']$. The hybrid stratum $\mg!_G^\hyb$ can then be identified as
\[\mg!_{G}^\hyb = \bigsqcup_{[\pi]\in \Pi(E) /\aut(G)} \mg!_{G, \pi}^\hyb \]
where the union is over all orbits of ordered partitions under the action of $\aut(G)$ (here, for each orbit in $\Pi(E) /\aut(G)$ we have fixed an arbitrary representative $\pi$).

\begin{defi}[Hybrid moduli spaces as a set] The  {\em moduli space $\mg_{g,n}^\hyb$ of hybrid curves of genus $g$ with $n$ markings}  is defined by
\[\mg_{g,n}^\hyb \coloneqq   \bigsqcup_{G} \bigsqcup_{[\pi] \in \Pi(E(G))/\aut(G)} \mg!_{G, \pi}^\hyb \]
where the first union is over all stable graphs $G$ of genus $g$ with $n$ markings  (again, $\pi$ is an arbitrary ordered partition representing the corresponding orbit).
\end{defi}
By construction, the points $\thy$ of $\mg_{g,n}^\hyb$ are in bijection with the isomorphism classes of marked hybrid curves. We denote by $\hcurvef_{g,n}$ the {\em universal hybrid curve} over $\mg_{g,n}^\hyb$. More precisely, as in the case of the Deligne-Mumford stack $\mgbar_{g}$, passing to the hybrid replacement $B^\hyb$ of étale charts $B$ given by the versal deformation spaces, we get a well-defined family of hybrid curves $\hcurvef_{g,n}$. Over the space $\mg_{g,n}^\hyb$, that we regard formally as the coarse moduli space of hybrid curves, the (set-theoretic) fiber $\hcurvef_\thy$ over some base point $\thy \in \mg_{g,n}^\hyb$ is the quotient of the hybrid curve represented by $\thy$ under the action of its automorphism group. Here, an automorphism of a hybrid curve is an automorphism of the underlying stable Riemann surface such that the induced automorphism on the dual graph respects the ordered partition and edge lengths.

\subsection{Hybrid topology} In this section, we describe the hybrid topology on the hybrid moduli spaces.

Consider a point $s_0$ in  $\mgbar_{g,n}$, and denote by $B$ the analytified versal deformation space of the corresponding stable curve (with $n$ marked points) $S_0 \coloneqq  \unicurve_{s_0}$. This way, we get an \'etale open chart (for the fine moduli space) $B \to \mgbar_{g,n}$ around $s_0$.  Let $B^\hyb$ be the hybrid deformation space defined in Section~\ref{sec:hybrid_base}. Denote by $G_{}$ the dual stable graph of $S_0$ with $n$ markings.

By construction, the automorphism group $\aut(S_0)$ of the marked stable Riemann surface $S_0$ acts on $B$. Since every automorphism of $S_0$ defines an automorphism of the dual stable graph $G$, this extends to a natural action of $\aut(S_0)$ on $B^\hyb$. We equip the quotient $B^\hyb/\aut(S_0)$ with the quotient topology induced from the hybrid topology of $B^\hyb$. Combining the above, we thus get a map from $B^{\hyb} / \aut(S_0)$ to $\mg_{g,n}^\hyb$.
\begin{prop} Notations as above, the collection of sets consisting of images of $B^{\hyb} / \aut(S_0)$ in $\mg_{g,n}^\hyb$ form a covering of $\mg_{g,n}^\hyb$. 
\end{prop}
\begin{proof} By definition of a covering, we must show that the union of these images is equal to $\mg_{g,n}^\hyb$. The points of $\mg_{g,n}^\hyb$ bijectively correspond to isomorphism classes of hybrid curves. Hence the claim is obvious, since (the isomorphism class of) each hybrid curve trivially belongs to the image of $B^{\hyb} / \aut(S_0)$ for its stable Riemann surface $S_0$. 
\end{proof}

\begin{defi}[Hybrid topology] We say that the images of sets of the form $B^{\hyb}/\aut(G)$ are all open in $\mg_{g,n}^\hyb$, and define the hybrid topology on $\mg_{g,n}^\hyb$ as the one generated by this open covering.
\end{defi}

\begin{thm} The space $\mg_{g,n}^\hyb$ is compact.
\end{thm}
\begin{proof} This follows from the compactness of $\mggbar{g,n}$ and local compactness of the spaces $B^\hyb$, stated in Proposition~\ref{prop:TopologicalPropertiesHybridSpaces}. 
\end{proof}

\subsection{Universal hybrid curve} \label{sec:universal_hybrid_curve}  It follows that locally, for the local charts $B^\hyb/\aut(S_0)$ considered in the previous section, the universal hybrid curve $\hcurvef_{g,n}$ restricts to the versal hybrid curve $\rsf^\hyb\to B^\hyb$ (see Section~\ref{sec:versal_hybrid_curve}) descended over $B^\hyb/\aut(S_0)$. The latter means that we first extend the action of $\aut(S_0)$ from the base $B^\hyb$ to the versal hybrid curve $\rsf^\hyb$ and then take the quotient. We endow the universal hybrid curve $\hcurvef_{g,n}$ with the topology induced by these local charts.

\subsection{Towers of hybrid moduli spaces} \label{sec:towermoduli} 
Using the constructions outlined in Section~\ref{sec:relations}, we can actually go further and define a tower of hybrid spaces 
\begin{align}\label{eq:tower-moduli}
\mggbar{g,n} \longleftarrow \mggbar{g,n}^{\hybr{1}} \longleftarrow \mggbar{g,n}^{\hybr{2}} \longleftarrow \dots \longleftarrow \mggbar{g,n}^{\hybr{\ndim-1}} \longleftarrow \mggbar{g,n}^{\hybr{\ndim}} = \mg_{g,n}^\hyb
\end{align}
\emph{interpolating} between $\mgbar_{g,n}$ and $\mg_{g,n}^{\hyb}$. Here $\ndim=3g-3+n$ is the dimension of $\mg_{g,n}$. We briefly discuss this here and refer to~\cite{AN-hybrid-green} for more details. 

The elements of $\mggbar{g,n}^{\hybr{k}}$ correspond to (equivalences classes of) marked hybrid curves of rank at most $k$, however the definition is slightly relaxed in order to ensure compactness. Namely, we allow the presence of some zero edge lengths in the last layer $\pi_k$ for hybrid curves of rank equal to $k$ (while at the same time keeping the normalization condition $\sum_{e \in \pi_k} \ell_e = 1$). 
The $k$-th hybrid space $\mggbar{g,n}^{\hybr{k}}$ compactifies $\mg_g$ by adding these (generalized) hybrid curves as its boundary part. The moduli space $\mg_{g,n}^{\hybr{k}}$ of (marked) hybrid curves of rank at most $k$ is an open subset of $\mggbar{g,n}^{\hybr{k}}$ and corresponds to those elements of $\mggbar{g,n}^{\hybr{k}}$ with all edge lengths non-zero. Note that for $k=\ndim$, we have the equality $\mggbar{g,n}^{\hybr{\ndim}} = \mg_{g,n}^\hyb$. The compactness of the spaces $\mggbar{g,n}^{\hybr{k}}$ can be deduced from the compactness of $\mggbar{g,n}$ and the local compactness of the spaces $B^{\hybr{r}}$ from Section~\ref{sec:relations}.

Moreover, over each member of the tower, we get the corresponding universal hybrid curve $\hcurve^{^{\Hyb(r)}}_{g,n}$, forming a tower 
\begin{align}\label{eq:tower-universalcurve}
\hcurve_{g,n} \longleftarrow \hcurve_{g,n}^{\hybr{1}} \longleftarrow  \hcurve_{g,n}^{\hybr{2}} \longleftarrow \dots \longleftarrow \hcurve_{g,n}^{\hybr{N-1}} \longleftarrow \hcurve_{g,n}^\hyb
\end{align}
compatible with~\eqref{eq:tower-moduli}, in the sense that via the maps from the universal hybrid curves to the moduli spaces, the  diagrams are commutative.


\section{Canonical measures}\label{sec:measures}
In this section, we introduce \emph{canonical measures} on different geometric objects: metric graphs, Riemann surfaces, tropical curves, and hybrid curves.

\subsection{Canonical measure on metric graphs} \label{ss:MeasuresGraphs}

Let $G=(V, E)$ be a finite graph of genus $h$ and $\ell  \colon E \to \R_+$ an edge length function. The corresponding metric graph is denoted by $\mgr$ (see Section \ref{ss:MetricGraphs}).  There are several equivalent ways to introduce the canonical measure on $\mgr$. The original definition of Zhang is based on the use of the metric graph Laplacian and its Green functions~\cite{Zhang}. In the following, we give three different expressions which are convenient for our purposes.

First of all, recall the following definition in terms of \emph{spanning trees}.
For a spanning tree $T = (V(T), E(T))$ of $G$, define its weight $\omega(T)$ as
\begin{equation} \label{eq:weight_spanning_tree}
 	\omega(T) = \prod_{e \in E \setminus E(T)} \ell_e.
 \end{equation}
The set of spanning trees of $G$ is denoted by $\cT(G)$. The \emph{Foster coefficients} $\mu(e)$, $e \in E$, of the metric graph $\mgr$ are defined as
 \begin{equation} \label{eq:graphmesspanning}
\mu(e) =  \frac{1}{\sum_{T \in \cT(G)} \omega(T)} \sum_{\substack{T \in \cT(G) \colon e \notin E(T) }} \omega(T),
\end{equation}
see, e.g., \cite[equation (TC1)]{BF11}. That is, $\mu(e)$ is the probability that a random spanning tree, sampled proportional to the weights~\eqref{eq:weight_spanning_tree}, does not contain the edge $e$.

The canonical measure on $\mgr$ is the measure
\begin{equation} \label{eq:CMFoster}
\mu_\Zh = \sum_{e \in E}\frac{\mu(e)}{\ell_e} \, d \theta_e  ,
\end{equation}
where $d\theta_e$  is the uniform Lebesgue measure on the edge $e \cong [0, \ell_e]$, so that $\frac 1{\ell_e}d \theta_e$ has total mass one on $e$. Since the complement of a spanning tree has $g$ edges, it follows that $\mu_\Zh$ has total mass $g$.

\smallskip

Alternatively, the canonical measure can be expressed using the \emph{cycle space} of $G$. Let $H_1(G, \Z)$ be the first homology of $G$. We fix an orientation for the edges of the graph. By definition, we have an exact sequence 
 \[ 0\to H_1(G, \Z) \longrightarrow \Z^E \longrightarrow \Z^{V} \to 0.
 \]
 For each edge $e$ of $G$, denote by $\langle\cdot\,,\cdot\rangle_e$ the bilinear form on $\R^E$ defined by
 \[\langle x\,,y\rangle_e \coloneqq  x_e y_e\]
 for any pair of elements $x = (x_f)_{f\in E}, y=(y_f)_{f\in E} \in \R^E$. We denote by $q_e$ the corresponding quadratic form.
The edge length function $\ell\colon E \to \R_+$ defines an inner product on $\R^E$ by 
 \[
 	\langle x,y \rangle_\ell \coloneqq  \sum_{e \in E} \ell_e \, \langle x\,,y\rangle_e = \sum_{e \in E} \ell_e\,  x_e y_e, \qquad x,y \in \R^E.
 \]
The associated quadratic form is denoted by $q_\ell$. Let $\pi \colon \R^E \to H_1(G, \R)$ be the orthogonal projection w.r.t. $\langle \cdot \,,\cdot\rangle_\ell$ of $\R^E$ onto its subspace $H_1(G, \R)$. The canonical measure $\mu_\Zh$ can be written as (e.g., \cite[Proposition 5.2]{SW19})
\begin{equation} \label{eq:graphmes}
\mu_\Zh = \sum_{e \in E} \frac{1}{\ell_e^2} \, q_\ell \Bigl(\pi(e)\Bigr) \,  d\theta_e.
\end{equation}

Finally, we derive a Hodge-theoretic expression for $\mu_\Zh$,  which resembles the canonical measure on Riemann surfaces (see \eqref{eq:cannonical_measure_RS} below). After fixing a basis $\gamma_1, \dots, \gamma_h$ of $H_1(G, \Z)$, we can identify the quadratic form $q_e$ restricted to $H_1(G, \R)$ with an $h \times h$ symmetric matrix $M_e$ so that, thinking of elements of $H_1(G, \R)$ as column vectors, we have
\begin{equation}\label{eq:1}
 q_e(x) = {x}^\transpose M_e x.
\end{equation}

 For an edge length function $\ell\colon E \to \R_+$, the restriction of $q_\ell$ to $H_1(G,\R)$ can be identified with the $h\times h$ symmetric matrix
 \begin{equation} \label{eq:M_l}
 	M_\ell = \sum_{e\in E}\ell_eM_e.
 \end{equation}
 We denote its inverse matrix by $M^{-1}_\ell$. Recall from Section~\ref{sec:basic-notations} that the $(i,j)$-coordinate of a matrix $A$ is denoted by $A(i,j)$. For an element $\gamma \in H_1(G, \R)$ and $e\in E$, denote by $\gamma(e)$ the $e$-coordinate of $\gamma$.

 \begin{thm}[Canonical measure via the graph period matrix] \label{thm:cmgraphscycles}The canonical measure on $\mgr$ is the measure
 \begin{equation} \label{eq:CMFoster2}
\mu_\Zh = \sum_{e \in E}\frac{\mu(e)}{\ell_e} \, d \theta_e,
\end{equation}
where $d\theta_e$ is the uniform Lebesgue measure on the edge $e \cong [0, \ell_e]$ and the Foster coefficients $\mu(e)$, $e \in E$, are given by 
 \[\mu(e) = \ell_e \sum_{i,j=1}^h M^{-1}_\ell(i,j) \, \gamma_i(e) \gamma_j(e), \]
with $\gamma_1,\dots,\gamma_h$ a basis of $H_1(G, \R)$.
 \end{thm}

 \begin{proof} Taking into account \eqref{eq:graphmes}, it suffices to show that for all edges $e \in E$,
 \[q_\ell \Bigl(\pi (e)\Bigr) = \ell_e^2 \sum_{i,j=1}^h \gamma_i(e) \gamma_j(e) M_\ell^{-1}(i,j). \]
Writing $\pi(e)$ in the basis $\gamma_1, \dots, \gamma_h$ gives 
\[\pi(e) = \sum_{i=1}^h a_{e,i}\gamma_i, \]
where the coefficents $a_{e,i}$ are determined by the system of linear equations
\[\sum_{j}a_{e, j} \langle\gamma_j, \gamma_i \rangle_\ell  = \langle\pi(e), \gamma_i \rangle_\ell = \langle e, \gamma_i \rangle_\ell = \ell_e \gamma_i(e)  \]
for $e \in E$ and $i\in[h]$. This implies that
\[\Bigl(a_{e,j}\Bigr)_{j=1}^h  = M_\ell^{-1} \Bigl( \ell_e \gamma_j(e)\Bigr)_{j=1}^h.\] 
We conclude by combining this expression with $q_\ell \Bigl(\pi(e)\Bigr) = \langle\pi(e), \pi (e) \rangle_\ell$.
 \end{proof}

\subsection{Canonical measure on a Riemann surface}
Let $S$ be a compact Riemann surface of positive genus $g$. Denote by $\Omega^1(S)$ the vector space of holomorphic one-forms $\omega$ on $S$ endowed with
the Hermitian inner product 
\[ \langle \omega_1, \omega_2\rangle \coloneqq \frac \imi{2} \int_S \omega_1 \wedge \bar{\omega}_2\]
for $\omega_1, \omega_2 \in \Omega^1(S)$. 

Choosing an orthonormal basis  $\eta_1, \dots,\eta_g$ of  $\Omega^1(S)$, the canonical  measure of $S$ is defined by
\[\mu_{\Ar} = \frac{\imi}{2}\sum_{j=1}^g \eta_j \wedge \bar{\eta}_j.\]

An alternate description of $\mu_\Ar$ can be given as follows. Choose a symplectic basis $a_1, \dots, a_g$, $b_1, \dots, b_g$  of $H_1(S, \Z)$, meaning that for the intersection pairing $\langle\cdot\,, \cdot \rangle$ between 1-cycles in $S$, we have for all $i, j = 1, \dots, g$,
\[\langle a_i, a_j \rangle = \langle b_i, b_j\rangle =0 \qquad \textrm{and} \qquad \langle a_i, b_j \rangle =\delta_{i,j}.\]
Let $\omega_1, \dots, \omega_g$ be an adapted basis of $\Omega^1(S)$ in the sense that 
\[\int_{a_j} \omega_i =\delta_{i,j}\]
and define the period matrix of $S$ by
\[\Omega = \Bigl(\Omega(i,j)\Bigr)_{i,j=1}^g = \Bigl(\int_{b_j}\omega_i\Bigr)_{i,j=1}^g.\]
The imaginary part $\Im(\Omega)$ of $\Omega$ is a symmetric positive definite real matrix.  Denote by $\Im(\Omega)^{-1}$ its inverse. Then,
\begin{equation} \label{eq:cannonical_measure_RS}
\mu_{\Ar} = \frac \imi{2} \sum_{i,j=1}^g \Im(\Omega)^{-1}(i,j) \, \omega_i \wedge \bar{\omega}_j
\end{equation}
where $\Im(\Omega)^{-1}(i,j)$ is the $(i,j)$-coordinate of the matrix $\Im(\Omega)^{-1}$.

\subsection{Tropical curves} \label{sec:cmtropicalcurve} Let $\hcurve^\trop = (\mgr,\pi)$ be a tropical curve with underlying finite graph model $G = (V, E)$, ordered partition $\pi = (\pi_1, \dots, \pi_r)$ on $E$ and edge length function $\ell \colon  E\rightarrow \R_+$ satisfying the layerwise normalization property. Let $\mgr$ be the associated augmented metric graph and $\genusfunction \colon V \to \N\cup\{0\}$ its genus function.

Consider the metric graph $\grm^j_\pi(\mgr)$ associated to the graded minor $\grm_\pi^j(G)$ and the edge length function $\ell^j \coloneqq  \ell\rest{\pi_j} \colon \pi_j \to \R_+$.  Denote by $\mu_\Zh^j$ the canonical measure on $\grm^j_\pi(\mgr)$.  By a slight abuse of the notation, we also denote by $\mu_\Zh^j$ the corresponding measure on $\mgr$ with support in the intervals $\Ical_e$ for $e\in \pi_j$. We can write 
\[\mu_\Zh^j = \sum_{e\in \pi_j} \frac{\mu^j(e)}{\ell_e} d\theta_e\]
for the uniform Lebesgue measure $d\theta_e$ on the interval $\Ical_e$ and the Foster coefficients $\mu^j(e)$, $e \in \pi_j$, of the metric graph $\grm^j_\pi(\mgr)$.

The canonical measure $\mu^{\can}$ of  $\hcurve^\trop$ is the measure on $\mgr$ defined by
\begin{equation} \label{eq:TropicalCanonicalMeasure}
\mu^\can \coloneqq  \sum_{v\in V} \genusfunction(v) \delta_v + \sum_{j=1}^r\sum_{e\in \pi_j} \frac{\mu^j(e)}{\ell_e}  d\theta_e.
\end{equation}
It follows from the considerations in Section~\ref{ss:LayeredSpanningTrees} and Section~\ref{ss:MeasuresGraphs} that $\mu^j(e)$ is the probability that a random layered spanning tree of the tropical curve does not contain the edge $e$ (when sampled proportional to the product of the lengths of edges not present in the tree).

\begin{example} Consider first the tropical curve of rank two with underlying layered graph from Example~\ref{fig:layered_spanning_trees}, and with edge lengths all equal to  $1/2$. The corresponding probability distribution is uniform on four layered spanning trees. The resulting canonical measure is the Lebesgue measure. The total mass is two.

Consider now the same graph with ordered partition consisting of a single set, namely the set of all edges $\pi_1 = E$, and with edge lengths all equal to $1/2$. The probability distribution on five spanning trees is uniform. In a random spanning tree, an edge is absent with probability either $2/5$ or $3/5$. In the associated canonical measure, two edges have a mass of $2/5$ and the other two a mass of $3/5$. The total mass is again two.
\end{example}
\subsection{Hybrid curves} \label{ss:cmhybridcurves}
Consider a hybrid curve $\hcurve^\hyb$ with underlying metrized complex  $\mc$, graph $G=(V,E)$, ordered partition $\pi$ on $E$, and edge length function $\ell \colon E \to\R_+$ satisfying the normalization property. Denote by $\hcurve^\trop$ the associated tropical curve. For each vertex $v \in V$, let $C_v$ be the corresponding smooth compact Riemann surface of genus $\genusfunction(v)$.

We define the canonical measure $\mu^\can$ on $\hcurve^\hyb$ as the measure which restricts to the canonical measure of $\hcurve^\trop$ on the intervals $\Ical_e$, for each edge $e\in E$, and which coincides with the canonical measure on each Riemann surface $C_v$, $v \in V$ (on components $C_v$ of genus $\genusfunction(v) = 0$, we set $\mu^\can \equiv 0$).

The following is straightforward.

\begin{prop} \label{prop:PushoutMeasures} Let $\forget \colon \hcurve^\hyb \to \hcurve^\trop$ be the natural projection map which contracts each Riemann surface $C_v$ to the vertex $v$ of $\hcurve^\trop$. Then the pushout of the canonical measure on $\hcurve^\hyb$ by $\forget$ is the canonical measure of $\hcurve^\trop$.
\end{prop}

Note that canonical measures on hybrid curves provide a common generalization of all previously introduced notions of canonical measures.

\begin{figure}[!t]
\centering
    \scalebox{.3}{\input{canonical_measures.pspdftex}}
\caption{Example of two hybrid curves with the same underlying stable Riemann surface. The canonical measures have the same Archimedean parts. The non-Archimedean parts are however different.  Note that changing the order of $\pi_1, \pi_2$ in the first example would drastically change the non-Archimedean part.}
\label{fig:canonical_measures}
\end{figure}

\section{Continuity of the universal canonically measured tropical curve} \label{sec:tropical}

The aim of this section is to study the continuity of canonical measures on the universal curve over $\mg!_G^\trop$, the moduli space of tropical curves with combinatorial type $G$. We start by recalling the definition of $\mg!_G^\trop$, making precise its topology, and defining the universal tropical curve $\hcurve!_G^\trop$.  We then obtain the following result that will be needed in the proof of Theorem~\ref{thm:mainglobal-intro}.

\begin{thm}[Continuity: tropical curves] \label{thm:mainmetricgraphs}
The universal family of canonically measured tropical curves $\hcurve!_G^\trop$ of given combinatorial type $G=(V, E, \genusfunction)$ over $\mg!_G^\trop$ is a continuous family. That is, for every continuous function $f \colon \hcurve!_G^\trop \to \R$, the function $F \colon \mg!_G^\trop \to \R$ defined by integration
\begin{align*}
F(x) \coloneqq  \int_{\hcurve^\trop_x}  f_{|_{\hcurve^\trop_x}} \, d\mu^\can_x, \qquad  x \in \mg!_G^\trop,
\end{align*}
is continuous on $\mg!_G^\trop$.
\end{thm}

In the above statement, $\hcurve^\trop_x$ denotes the fiber over a point $x \in \mg!_G^\trop$ in the universal family $\hcurve!_G^\trop$. It coincides with the quotient of the tropical curve $\widetilde {\mathscr C}_x^{\,  \trop}$ associated to $x$ by the action of its automorphism group. Each fiber $\hcurve^\trop_x$ is equipped with its canonical measure $\mu^\can_x$, the pushout of the canonical measure on the tropical curve $\widetilde {\mathscr C}_x^{\,  \trop}$ through the quotient map.

\subsection{The moduli space $\mg!_G^\trop$} Fix an augmented graph $G=(V, E, \genusfunction)$ and let $\mg!_G^\trop$ be the moduli space of tropical curves of combinatorial type $G$ (see Section~\ref{sec:moduli_space_tropical_curves}). Recall that, as a set, $\mg!_G^\trop$ is defined as the quotient
\begin{equation} \label{eq:mgtrop}
\mg!_G^\trop = \ssub{\widetilde \mg}!_G^{\,\,\,\trop}/\aut(G)
\end{equation}
of the set
\begin{equation} \label{eq:SpaceTropicalCurvesFixedGraph}
\ssub{\widetilde \mg}!_G^{\,\,\,\trop} = \bigsqcup_{\pi \in \Pi(E)} \inn \sigma_{\pi}.
\end{equation}
In the following, we endow the spaces $\mg!_G^\trop$ and $\ssub{\widetilde \mg}!_G^{\,\,\,\trop}$ with suitable topologies. We define a topology on $\ssub{\widetilde \mg}!_G^{\,\,\,\trop}$ by specifying for each $x \in \ssub{\widetilde \mg}!_G^{\,\,\,\trop}$ its system of neighborhoods $\mathcal{U}(x)$, and then put the quotient topology on $\mg!_G^\trop$.

It suffices to define a neighborhood base $\bigl(\,U^\varepsilon(x)\,\bigr)_{\varepsilon >0}$ around each point $x \in \ssub{\widetilde \mg}!_G^{\,\,\,\trop}$. Hence, assume that $x \in \ssub{\widetilde \mg}!_G^{\,\,\,\trop}$ and $\varepsilon > 0$ are fixed. Suppose that  $x\in \inn \sigma_\pi$ for the ordered partition $\pi$ of $E$.
 We will introduce $U^\varepsilon(x)$ as a union
 \[
	U^\varepsilon (x) \coloneqq  \bigcup_{\pi' \preceq \pi} U_{\pi'}^\varepsilon,
\]
where the subset $U_{\pi'}^\varepsilon \subseteq \inn \sigma_{\pi'}$ is defined as follows: \\
Assume that $\pi' = (\pi_i')_{i=1}^r$ is an ordered partition of $E$ with $\pi' \preceq \pi$. Then, $\pi$ has the following form 
\begin{equation} \label{eq:preceqdecomp_tropical}
\pi = \Big( \varrho^1, \varrho^2, \dots, \varrho^r \Big ) = \Big( (\varrho^i)_{i=1}^r \Big ),
\end{equation}
where $\varrho^i =(\varrho_k^i)_{k=1}^{s_i}$ is an ordered partition of $\pi_i'$ for $i=1,\dots,r$. \\
A point $y \in \inn\sigma_{\pi'}$ belongs to $U_{\pi'}^\varepsilon$ if the following condition is satisfied. For each $i \in \{1,\dots,r\}$, consider the ordered partition $\varrho^i =(\varrho_k^i)_{k=1}^{s_i}$ of $\pi_i'$. Then,
\begin{equation} \label{eq:condpart1_tropical}
	\max_{k =1, \dots, s_i-1} \left\{ \frac{\sum_{e \in \varrho_{k+1}^i} y_e}{ \sum_{e \in \varrho_{k}^i } y_e} \right \}  < \varepsilon \qquad  \text{ and }  \qquad \max_{k =1, \dots, s_i} \max_{e \in \varrho_k^i} \Big|\frac{y_e}{\sum_{\hat e \in \varrho_k^i} y_{\hat e}} -x_e \Big | < \varepsilon.
\end{equation}
As stated in the following theorem, this defines a topology on $\ssub{\widetilde \mg}!_G^{\,\,\,\trop}$.
\begin{thm} 
There is a unique topology on $\ssub{\widetilde \mg}!_G^{\,\,\,\trop}$ such that $\mathcal{U}(x)$ coincides with the  system of neighborhoods for any  $x \in \ssub{\widetilde \mg}!_G^{\,\,\,\trop}$. 
\end{thm}

This topology is compatible with the hybrid topology on $\mg_{g,n}^\hyb$. Making this precise requires the introduction of logarithm maps, which is out of the scope of this paper. We refer instead to our paper~\cite[Section 13]{AN-hybrid-green}.

It is also clear that $\ssub{\widetilde \mg}!_G^{\,\,\,\trop}$ is a second countable and compact Hausdorff space. Similar to Proposition~\ref{prop:convseq}, we can describe the convergent sequences.

\begin{prop} \label{prop:convseq_tropical}
Let $x \in \ssub{\widetilde \mg}!_G^{\,\,\,\trop}$ and assume $x \in \inn \sigma_\pi$ for an ordered partition $\pi \in \Pi(E)$. Suppose that $(x_n)_n$ is a sequence in $\ssub{\widetilde \mg}!_G^{\,\,\,\trop}$. Then, the following statements hold.
\begin{itemize}[leftmargin = 2em]
\item [(a)] If $x_n$ converges to $x$ in $\ssub{\widetilde \mg}!_G^{\,\,\,\trop}$, then, almost all $x_n$ belong to strata $\inn \sigma_{\pi'}$ of ordered partitions $\pi'$ of $E$ with $\pi' \preceq \pi$.

\item [(b)] Assume that $(x_n)_n \subseteq \inn \sigma_{\pi'}$ for some fixed ordered partition $\pi' \preceq \pi$. In particular, $\pi$ is of the form \eqref{eq:preceqdecomp_tropical}. Then, $x_n$ converges to $x$ in $\ssub{\widetilde \mg}!_G^{\,\,\,\trop}$ if and only if the following conditions hold:
\end{itemize}
\begin{itemize}
\item  Let $\varrho^i =(\varrho^i_k)_{k=1}^{s_i}$ be one of the ordered partitions in \eqref{eq:preceqdecomp_tropical}. Then, for all $k$ and $e \in \varrho_k^i$,
\begin{equation} \label{eq:conv3_tropical}
\lim_{n\to\infty} \frac{x_{n,e}}{\sum_{\hat e \in \varrho_k^i} x_{n, \hat e}} =x_e.
\end{equation}
Moreover, if $e \in \varrho_k^i$ and $e' \in \varrho_{k'}^i$ with $k < k'$, then
\begin{equation} \label{eq:conv4_tropical}
\lim_{n\to\infty}\frac{x_{n,e'}}{x_{n,e}}=0
\end{equation}

\end{itemize}
 \end{prop}
\begin{proof}
The claims follow from the definition of the topology on $\ssub{\widetilde \mg}!_G^{\,\,\,\trop}$.
\end{proof}

\subsection{Universal tropical curve of given combinatorial type} \label{ss:thm-tropicalcurve}
Each point $x \in \mg!_G^\trop$ or $x \in \ssub{\widetilde \mg}!_G^{\,\,\,\trop}$ clearly represents a tropical curve $\widetilde {\mathscr C}_x^{\, \trop}$ of combinatorial type $G=(V, E, \genusfunction)$. Namely, $\widetilde {\mathscr C}_x^{\,  \trop}$, $x \in  \ssub{\widetilde \mg}!_G^{\,\,\,\trop}$, consists of the augmented metric graph $\mgr_x$ (i.e., the metric realization of the pair $(G, x)$) and the ordered partition $\pi$ of the edge set $E$ associated to $x$ (i.e., such that $\inn \sigma_\pi \ni x$), and similar for $x \in \mg!_G^\trop$.

The space $\ssub{\widetilde \mg}!_G^{\,\,\,\trop}$ supports a natural family of tropical curves $\ssub{\widetilde {\mathscr C}}!_G^{\,  \trop}$, whose fibers are the respective tropical curves (see below). The action of the automorphism group $\aut(G)$ on $\ssub{\widetilde \mg}!_G^{\,\,\,\trop}$ extends naturally to an action on $\ssub{\widetilde {\mathscr C}}!_G^{\,  \trop}$. We define the \emph{universal tropical curve} $\hcurve!_G^\trop$ over  $\mg!_G^\trop$ as the topological quotient
\[
\hcurve!_G^\trop \coloneqq  \ssub{\widetilde {\mathscr C}}!_G^{\,  \trop} / \aut(G).
\]
There is a natural projection map from $\hcurve!_G^\trop$ to $\mg!_G^\trop = \ssub{\widetilde \mg}!_G^{\,\,\,\trop}/\aut(G)$. The topological fiber $\hcurve^\trop_x$ over $x \in \mg!_G^\trop$ is the quotient of the tropical curve $\widetilde {\mathscr C}_x^{\,  \trop}$ associated to $x$ by the action of its automorphism group.

The construction of the family $\ssub{\widetilde {\mathscr C}}!_G^{\,  \trop}$ is done by standard techniques and included here only for completeness (see~\cite{ACP, CCUW20} for similar considerations). We first combine the interval fiber bundles $\Ical_\pi^e$ over $\inn \sigma_\pi$ (see \eqref{eq:Ical_pi^e}) to an interval fiber bundle $\Ical_e = \bigsqcup_{\pi \in \Pi(E)} \Ical_\pi^e$, $e \in E$, over $\ssub{\widetilde \mg}!_G^{\,\,\,\trop}$. The  family of tropical curves $\ssub{\widetilde {\mathscr C}}!_G^{\,  \trop}$ over $\ssub{\widetilde \mg}!_G^{\,\,\,\trop}
$ is the topological quotient
\[
\ssub{\widetilde {\mathscr C}}!_G^{\,  \trop} \coloneqq   \bigsqcup_{e \in E} \Ical_e / \sim,
\]
where, above any point in the base $\ssub{\widetilde \mg}!_G^{\,\,\,\trop}
$, we  glue the extremities of the intervals $\Ical_e$  using the incidence relations in the graph $G$.

\subsection{Proof of Theorem~\ref{thm:mainmetricgraphs}} The proof of Theorem~\ref{thm:mainmetricgraphs} is based on the following lemma on convergence of Foster coefficients.
\begin{lem}\label{lem:sptreesct}
Let $\ell$ be a point in $\ssub{\widetilde \mg}!_G^{\,\,\,\trop}
$ and $\pi = (\pi_j)_{j=1}^r$ its ordered partition, that is, $\ell \in \inn \sigma_\pi$. Suppose that $(\ell_n)_n \subseteq  \ssub{\widetilde \mg}!_G^{\,\,\,\trop}
$ is a sequence with $\lim_{n \to \infty} \ell_n=\ell$ in $\ssub{\widetilde \mg}!_G^{\,\,\,\trop}
$ and $\ell_n \in \inn \sigma_{(E)}$ for all $n$. Consider the Foster coefficients $\mu_n(e)$, $ e \in E$, of the metric graph $\mgr_n$ associated to $(V, E, \ell_n)$.

Then, for all layers $\pi_j$, $j=1,\dots,r$ and all edges $e \in \pi_j$,
\begin{equation} \label{eq:sptreesct}
	\lim_{n \to \infty} \mu_{n} (e)  = \mu_{\infty, j} (e),
	\end{equation}
where $\mu_{\infty, j}(e)$, $e \in \pi_j$, are the Foster coefficients of the metric graph $\grm^j_\pi(\mgr)$ associated to the graded minor $\grm_\pi^j(G)$ and the edge lengths $\ell_e$, $e \in \pi_j$.
\end{lem}

\begin{proof}
Recall the definition of the Foster coefficients using spanning trees \eqref{eq:graphmesspanning}. Fix a spanning tree $T \in \cT(G)$ and denote its weight in $\lg_n$ by
\[\omega_n(T) = \prod_{e \notin E(T)} \ell_{n,e} = \prod_{i=1}^r  \prod_{e \in \pi_i \setminus E(T)} \ell_{n,e}.
\]
We first give an asymptotic description of the weights $\omega_n(T)$ as $n \to \infty$.

Since $G_\pi^ i \cap T$ contains no cycles, Proposition~\ref{prop:genusformula} gives the inequality
\begin{equation} \label{eq:estgenusedges}
	 \bigl|\,E^i_\pi  \setminus E(T)\bigr|  = \sum_{j=i}^r\, \bigl|\,\pi_j \setminus E(T)\,\bigr|  \ge h(G^i_\pi) = \sum_{j=i}^r h^j_\pi, \qquad i = 1, \dots, r.
\end{equation}
If equality holds true in \eqref{eq:estgenusedges} for all $i$, then $(V, E^i_\pi \cap E(T))$ is a spanning tree of $G^i_\pi$ for all $i$. However, then $T_i \coloneqq  \pi_i \cap T$ is a spanning tree of $\grm^i_\pi$ for all $i$. Hence, $T = \bigcup_{i=1}^r T_i$ belongs to the set $\cT_\pi(G)$ of spanning trees of the layered graph $(G, \pi)$ (see Definition~\ref{def:LayeredSpanningTree}). By the equality $ | \pi_i \setminus E(T)| = h^i_\pi$  and Proposition~\ref{prop:convseq_tropical}, we infer that, as $n \to \infty$,
\begin{equation} \label{eq:LimitSpanningTrees1}
\omega_n(T) = \prod_{i=1}^r L_i(\lg_n)^{h_\pi^i} \, \Big (  \prod_{i=1}^r   \nu_i (T_i) + o(1) \Big )
\end{equation}
where $L_i(\lg_n) \coloneqq  \sum_{e \in \pi_i} \ell_{n,e}$ and $\nu_i (T_i) =\prod_{e\in \pi_i \setminus E(T)} \ell_e$ is the weight of the spanning tree $T_i$ in the metric graph $\grm^i_\pi(\mgr)$, which is obtained by equipping the $i$-th graded minor $\grm_\pi^i(G)$ of $(G, \pi)$ with the edge length function $\ell \rest{\pi_i}$. 

Conversely, assume that the inequality in \eqref{eq:estgenusedges} is strict for some $i_0 \in [r]$. As follows from \eqref{eq:estgenusedges}, the edge set $E(G) \setminus E(T)$ can be decomposed as
\[
	E \setminus E(T) = F_1 \sqcup F_2 \sqcup \dots \sqcup F_r,
\]
where each $F_i$ is a subset of $\bigcup_{j \ge i} \pi_j \setminus E(T)$ containing exactly $h^i_\pi$ edges. By our assumption on $T$, at least one edge set $F_i$ must contain an edge from some $\pi_j$ with $j > i$. Applying Proposition~\ref{prop:convseq_tropical}, we conclude that for $T \notin \cT_\pi(G)$,
\begin{equation} \label{eq:LimitSpanningTrees2}
\lim_{n \to \infty} \prod_{i=1}^r L_i(\lg_n)^{ - h_\pi^i}  \,  \omega_n(T) = 0.
\end{equation}
In order to prove \eqref{eq:sptreesct}, fix an edge $e$ in the $j$-th set  $\pi_j$, $j \in [r]$. Combining the representation \eqref{eq:graphmesspanning} with \eqref{eq:LimitSpanningTrees1} and \eqref{eq:LimitSpanningTrees2}, we finally obtain that 
\begin{align*}
	\lim_{n \to \infty}  \mu_n(e)   =  \sum_{\substack{T = (T_i)_i \in \cT_\pi(G) \colon \\ e \notin E(T)}} \,  \prod_{i=1}^r \nu_i (T_i) \,\bigl/{\sum_{T = (T_i)_i \in \cT_\pi(G)}}  \,  \prod_{i=1}^r \nu_i (T_i) = \mu_{\infty, j} (e),
\end{align*}
where the last equality follows by applying \eqref{eq:graphmesspanning} to the graded minor $\gr^j_\pi$.
\end{proof}

\begin{proof}[Proof of Theorem~\ref{thm:mainmetricgraphs}]
We will prove that the family $(\widetilde {\mathscr C}_x^{\,  \trop}, \mu^\can_x)$, $x \in \ssub{\widetilde \mg}!_G^{\,\,\,\trop}
$, where each tropical curve fiber $\widetilde {\mathscr C}_x^{\,  \trop}$ is equipped with its canonical measure $ \mu^\can_x$, is a continuous family of measured spaces over $\ssub{\widetilde \mg}!_G^{\,\,\,\trop}
$. Recalling that $\mg!_G^\trop$ and $\hcurve!_G^\trop$ are obtained as topological quotients, Theorem~\ref{thm:mainmetricgraphs} follows.

Since the canonical measure is defined by Foster coefficients (see \eqref{eq:TropicalCanonicalMeasure}), Lemma~\ref{lem:sptreesct} proves the continuity at each point $\ell \in \inn \sigma_\pi$ through the stratum $\ssub{\inn \sigma}!_{(E)}$. To prove the continuity through general strata, one can fix an ordered partition $\pi' \preceq \pi$ and apply Lemma~\ref{lem:sptreesct} to each graded minor $\gr^j_{\pi'}(G)$. 
\end{proof}

\section{Monodromy} \label{sec:monodromy}

This section provides basic results about degenerations of Riemann surfaces.  We recall an interpretation of the rank one symmetric matrices $M_e$ introduced in \eqref{eq:1} in terms of the monodromy of a degenerating family of Riemann surfaces. We then define admissible basis for layered graphs and hybrid curves. In the next section, we relate these results to the asymptotics of the period mapping.   We use the notations of Section~\ref{sec:deformations}.

Let  $S_0$ be a stable curve of genus $g$ with dual graph $G=(V, E, \genusfunction)$, and consider the analytic family $\rsf \to \base$ of stable Riemann surfaces over a polydisc $\base$ of dimension $\ndim=3g-3$ with fiber $\rsf_0 = S_0$ over $0\in B$. Let $B^* \coloneqq B \setminus \bigcup_{e \in E} D_e$ be the locus of points whose fibers in the family $\pr \colon \rsf \to B$ are smooth, and set $\rsf^* \coloneqq  \pr^{-1}(B^*)$.

Fix a base-point $\bp \in \base^*$. The fundamental group $\pi_1\bigl( \base^* , \bp\bigr)$ is isomorphic to $\Z^E$, with one generator $\lambda_e$ per edge $e\in E$ corresponding to a simple loop $\lambda_e\subset \base^*$.  The loop $\lambda_e$ is based at $\bp$ and turns once around the divisor $D_e$. Moreover, $\lambda_e$ is contractible in the space $\base \setminus \bigcup_{e'\neq e}D_{e'}$.

From the family $\pr \colon \rsf^* \to \base^*$, we get the local system $\mathcal H$ on $B^*$ with $\mathcal H_t\coloneqq  H_1(\rsf_{t},\Z)$, $t\in B^*$. This local system is classified by the \emph{monodromy action} of $\pi_1(B^*, \bp)$ on $H_1(\rsf_{\bp},\Z)$. The elements $\lambda_e$, $e\in E$, form a system of generators for the fundamental group, and 
{Picard-Lefschetz theory} describes  the action of $\lambda_e$ on $H_1(\rsf_{\bp},\Z)$ as follows (see~\cite[Section 6]{Lam81} and~\cite[Chapter 4]{CMP17}). There is a \emph{vanishing cycle $a_e$ in $H_1(\rsf_{\bp},\Z)$} associated to the singular point $p^e$ of $S_0$ such that the monodromy action of $\lambda_e$ is given by
\begin{align}
\lambda_e \colon H_1(\rsf_{\bp}, \Z) \longrightarrow H_1(\rsf_{\bp}, \Z)\\
\label{eq:PL} \beta \mapsto \beta-\langle \beta,a_e\rangle a_e.
\end{align}
Here $\langle\cdot\,, \cdot \rangle$ denotes the intersection pairing between one-cycles in $H_1(\rsf_{\bp})$.

\subsection{Specialization} After possibly shrinking the polydisc $B$, the inclusion $S_0 \hookrightarrow \rsf$ admits a deformation retraction $\rsf \rightarrow S_0$. It follows that $S_0$ and $\rsf$ have the same homology groups. From the composition of the inclusion $\rsf_{\bp} \hookrightarrow \rsf \simeq_{\textrm{homotopy}} S_0$ with the homotopy equivalence, we get the \emph{specialization map}
\begin{equation*}
\mathrm{sp} \colon H_1(\rsf_{\bp},\Z) \longrightarrow H_1(\rsf,\Z) \simeq H_1(S_0,\Z).
\end{equation*}

We have the following proposition, cf.~\cite{ABBF}.
\begin{prop}The specialization map $\mathrm{sp}$  is surjective. 
\end{prop}
\begin{proof} For the sake of completeness, we provide a proof. We take a one-dimensional disk $\Delta$ which goes through the point $\bp$ and $0$ in $B$ such that the restriction of the family over $\Delta$ is algebraic. We get a variation of Hodge structures $H^1(\rsf_t)$, $t\in \Delta^*$, that we can assume to have unipotent monodromy by going to a finite cover of $\Delta^*$. Denote by $N$ the logarithm of the monodromy map on the cohomology associated to the family over $\Delta^*$. Clemens-Schmid exact sequence in asymptotic Hodge theory~\cite{PS08} then gives the following exact sequence 
\begin{equation*}
 0 \to H^1(S_0) \xrightarrow{\mathrm{sp}^*} H^1(\rsf_{\bp}) \xrightarrow{N} H^1(\rsf_{\bp}).
\end{equation*} 
Passing to the dual, gives the surjectivity stated in the proposition. 
\end{proof}

Let $A\subset H_1(\rsf_{\bp},\Z)$ denote the kernel of the specialization map.  It corresponds precisely to the subspace spanned by the
vanishing cycles $a_e$, $e \in E$. We have an exact sequence
\begin{equation*}
0\to A \to H_1(\rsf_{\bp}, \Z) \xrightarrow{\mathrm{sp}} H_1(S_0, \Z) \to 0.
\end{equation*}

From the inclusion of subspaces $C_v \hookrightarrow S_0$, we get a short exact sequence 
\[0 \to \bigoplus_{v\in V} H_1(C_v, \Z) \hookrightarrow H_1(S_0, \Z) \to H_1(G, \Z) \to 0.\]

Define 
\[A'\coloneqq A + \mathrm{sp}^{-1}\Bigl(\bigoplus_{v \in V} H_1(C_v, \Z)\Bigr)\subseteq H_1(\rsf_{\bp}, \Z).\]

It follows that
\begin{equation}
 \label{eq:aux}
H_1(\rsf_{\bp}, \Z)/A' \simeq H_1(S_0, \Z) / \bigoplus_{v \in V} H_1(C_v, \Z)\simeq H_1(G, \Z) .
\end{equation}

Recall that a subspace $H\subseteq H_1(S_\bp,\Z)$ is called \emph{isotropic} if for any pair of elements $a, b \in H$, the intersection pairing $\langle a, b\rangle$ vanishes.

\begin{prop}\label{prop:isotropic} The subspace $A\subset
  H_1(\rsf_{\bp},\Z)$ of vanishing cycles has rank equal to $\graphgenus$ and is isotropic. Moreover, we have $\langle A,A'\rangle =0$.
  \end{prop}
\begin{proof}
By surjectivity of the specialization map, we have 
  \begin{align*}
    \rk A&= \rk H_1(\rsf_{\bp}, \Z) - \rk H_1(S_0, \Z) = 2g - \graphgenus - 2\sum_{v\in V}\genusfunction(v)= \graphgenus.
  \end{align*}
The vanishing cycles $a_e$ are disjoint for $\bp$ close to $0\in B$, and their self-intersections are zero, which shows the second claim.  The space $A'$ contains $A$ and all cycles which can be deformed to a cycle on a component of $S_0$. The elements of $A$ vanish in $S_0$, the pairing between $A$ and cycles in any $C_v$ is thus trivial, and the last claim follows. 
\end{proof}

It follows that the intersection pairing involving vanishing cycles reduces to a pairing
$A \times \Bigl(H_1(\rsf_{\bp}, \Z)/A'\Bigr) \rightarrow \Z$.

\subsection{Description of the monodromy action} For every edge $e \in E$, define 
\[N_e \coloneqq  \lambda_e - \id.
\] 
By \eqref{eq:PL} and the previous proposition, we get that $N_e\circ N_e =0$, that is, $N_e= \log(\lambda_e)$. Note that \eqref{eq:PL} implies that actually $N_{e}$ vanishes on $A'$.

It follows that $N_e\colon H_1(\rsf_{\bp}, \Z) \to A$ passes to the quotient by $A'$ and induces 
\begin{equation} \label{eq:nm}
H_1(
G, \Z) \simeq H_1(\rsf_{\bp}, \Z)/A' \xrightarrow{N_e} A \simeq
(H_1(\rsf_{\bp}, \Z)/A')^\vee \simeq H_1(G, \Z)^\vee.
\end{equation}

We recall the following result from~\cite[Proposition~3.4]{ABBF}.
\begin{prop} \label{prop:BilinearForm} The bilinear form on $H_1(G, \Z)$ given by the composition of maps in \eqref{eq:nm} coincides with the bilinear form $\langle\cdot\,,\cdot\rangle_e$ restricted to $H_1(G, \Z)$.
\end{prop}

\subsection{Admissible symplectic basis} \label{sec:admissible_basis_standard}
Enumerate the vertices of the graph $G$ as $v_1, \dots, v_n$. Define $g_0, \dots, g_n$ by 
\[g_i = \graphgenus+ \sum_{j=1}^{i} \genusfunction(v_i).\]
Note that $g_0 = \graphgenus$ and $g_n =g$.

An \emph{admissible symplectic basis} $a_1,\dotsc,a_{g},b_1,\dotsc,b_{g}$ for
 $H_1(\rsf_{\bp}, \Z)$ is a symplectic basis such that
 
 \begin{enumerate}
 \item[(1)] $a_1, \dots, a_\graphgenus$ form a basis of the space of vanishing cycles
 $A$, and

 \item[(2)] for any $r=1, \dots, n$, the collection of elements $a_{j}$ and $b_j$ for $j=g_{r-1}+1, \dots, g_r$ gives a basis of $H_1(C_{v_r}, \Z)$ in 
 \[H_1(C_{v_r}, \Z) \hookrightarrow H_1(S_0, \Z) \simeq H_1(\rsf_{\bp}, \Z) /A.\]
 \end{enumerate} 
 
 Note that in (2) the elements $a_{j}$ and $b_j$  seen in $H_1(C_{v_r}, \Z)$ necessarily  form a symplectic basis  of $H_1(C_{v_r}, \Z)$. Indeed, the pairing between $A$ and $A'$ is trivial in $H_1(\rsf_{\bp},\Z)$, and so  the intersection pairing in $H_1(C_{v_r}, \Z)$ coincides with the one induced from $H_1(\rsf_{\bp},\Z)$ via the specialization map.

 For any edge $e\in E$, the vanishing cycle $a_e$ can be decomposed as a linear combination of the basis elements $a_i$ of $A$: 
\begin{equation}\label{eq:vanishing_cycle}
a_e = \sum_{i=1}^\graphgenus c_{e,i}a_i.
\end{equation}
Let $\Bg \subset H_1(\rsf_{\bp}, \Z)$ be the subspace generated by the basis elements $b_1,
\dots, b_{\graphgenus}$. As we saw above, projection onto $H_1(\rsf_{\bp}, \Z) / A'$ provides an isomorphism $\Bg \simeq 
H_1(G, \Z)$. The pairing $\langle \cdot\,, \cdot \rangle$ gives an isomorphism $\Bg \simeq A^{\vee}$, and the monodromy operators $\lambda_e$ give maps
\begin{displaymath}
  N_{e}\colon \Bg \longrightarrow A.
\end{displaymath}
The following proposition is straightforward to prove.

\begin{prop}\label{prop:matricial} In terms of the basis $b_1,\dots, b_\graphgenus$ for $\Bg\simeq H_1(G, \Z)$, we can write
\[M_e=\Bigl(c_{e,i} \, c_{e,j}\Bigr)_{1\leq i,j\leq \graphgenus}
\]
for $M_e$ defined in Equation~\eqref{eq:1}.
\end{prop}

In the following, we will write $\basisa$ (respectively $\basisb$) for the subspace of
$H_1(\rsf_{\bp}, \Z)$ generated by $a_{1},\dots ,a_{g}$ (respectively
$b_{1},\dots,b_{g}$). Then $\basisa$ is a maximal isotropic subspace
with $A\subset \basisa\subset A'$.

We will adapt the definition of an admissible symplectic basis to hybrid curves in Section~\ref{sec:admissible_basis_hybrid_curves} (by adding a condition regarding the layering).

\subsection{Admissible basis for layered graphs and hybrid curves} \label{sec:admissible_layered}

In this section, we refine the notion of an admissible symplectic basis in context with hybrid curves. Later on, this will allow to control the asymptotics of period matrices under degeneration to a hybrid curve.

Let $G=(V, E, \pi)$ be a layered graph of genus $\graphgenus$ with vertex set $V$, edge set $E$, and ordered partition $\pi = (\pi_1, \dots, \pi_r)$ on $E$. Denote by $\dfilter_\pi^\bullet$ the decreasing filtration of $E$ introduced in Section~\ref{ss:GradedMinors},
\[\dfilter_\pi^\bullet\colon \quad E^1_\pi =E \supsetneq E^2_\pi \supsetneq \dots \supsetneq E^r_\pi \supsetneq E^{r+1}_\pi =\varnothing.\]
Consider the associated decreasing sequence of spanning subgraphs of $G$,
\[G \eqqcolon G^1_\pi \supset G^2_\pi \supset  \dots \supset G_\pi^{r} \supset G_\pi^{r+1}= (V, \varnothing)\]
with $G^j_\pi = (V, E^j_\pi)$, and let $\grm^j_\pi(G) = \contract{G^j_\pi}{E^{j+1}_\pi}$ be the $j$-th graded minor.

Recall that $\cont{j}\colon G^{j}_\pi \to \grm^j_\pi(G)$ is the contraction map. By an abuse of the notation, we denote by $\proj_j \colon H_1\bigl(G^{j}_\pi, \Z\bigr) \to H_1\bigl(\grm^j_\pi(G), \Z\bigr)$ the corresponding map on the level of homology groups. For future use, we record the following basic proposition.
\begin{prop}\label{prop:surjectivity} The map $\proj_j \colon H_1\bigl(G^{j}_{\pi}, \Z\bigr) \to H_1\bigl(\grm^j_\pi(G), \Z\bigr)$ is surjective.
\end{prop}

\subsubsection{Admissible basis for layered graphs} \label{sec:admissible_basis_layered_graphs} Notations as above, let $G=(V,E)$ be a layered graph of genus $\graphgenus$ with ordered partition $\pi=(\pi_1, \dots, \pi_r)$ on $E$.

By the genus formula in Proposition~\ref{prop:genusformula}, we have $\graphgenus = h^1_\pi + \dots + h^r_\pi$.  Consider the partition 
\begin{equation} \label{eq:ordered_partition_edges}
[\graphgenus] = J^1_\pi \,\sqcup\, J^2_\pi \sqcup \dots \sqcup J^r_\pi
\end{equation}
into intervals $J^j_\pi$ of size $h^j_\pi$ given by 
\[J^j_\pi \coloneqq  \Bigl\{1+\sum_{i=1}^{j-1}h^i_\pi,\,\, 2+\sum_{i=1}^{j-1}h^i_\pi,\,\, \dots,\,\,\sum_{i=1}^{j}h^i_\pi\Bigr\}.\]

An \emph{admissible basis} for the layered graph $(G, \pi)$ is a basis $\gamma_1, \dots, \gamma_h$ of $H_1(G, \Z)$ that satisfies the following conditions for all layers $\pi_j$, $j = 1, \dots, r$:
\begin{itemize}[leftmargin = 2em]
\item[$(i)$] For each $k \in J^j_\pi$, the cycle $\gamma_k$ lies in the spanning subgraph $G^{j}_\pi$ of $G$. In other words, all edges of $\gamma_k$ belong to $E^{j}_\pi \subseteq E$.
\item[$(ii)$] The cycles $\proj_j(\gamma_k)$ for $k\in J^j_\pi$ form a basis of $H_1(\grm^j_\pi(G), \Z)$.
\end{itemize}

\begin{prop}\label{prop:ExistenceAdmissibleBasis} Any layered graph $G=(V, E, \pi)$ admits an admissible basis. 

\end{prop}

\begin{proof} This is a consequence of the genus formula, cf. Proposition~\ref{prop:genusformula}, and Proposition~\ref{prop:surjectivity}.  We can indeed construct an admissible basis by first taking, for each $j=1, \dots, r$, a basis $\gamma^j_k$, $k\in J^j_\pi$, of $\grm^j_\pi(G)$, and then lifting these cycles to cycles $\gamma_k$ in $G$ via the (surjective) projection maps $\proj_j \colon G^{j}_\pi \to \grm^j_\pi(G)$. It is easy to see that the cycles $\gamma_k$ form a basis of $H_1(G, \Z)$, and they obviously verify the above properties $(i)$ and $(ii)$. 
\end{proof}

\subsubsection{Admissible basis for hybrid curves} \label{sec:admissible_basis_hybrid_curves} Let $\hcurve$ be a hybrid curve with underlying stable curve $S_0$ and layered dual graph $G = (V, E, \pi)$. Denote as in Section~\ref{sec:deformations} by $\genusfunction\colon V \to \N\cup\{0\}$ the genus function on vertices. 

Enumerate the vertices of the graph as $v_1, \dots, v_n$. Define $g_0, \dots, g_n$ by 
\begin{equation} \label{eq:DefineGi}
g_i = \graphgenus+ \sum_{j=1}^{i} \genusfunction(v_j).
\end{equation}
Note that $g_0 = \graphgenus$ and $g_n =g$.

Let ${\bp}$ be a fixed base point in $B^*$. An \emph{admissible symplectic basis} for $\hcurve$ is a symplectic basis $a_1, \dots,a_g, b_1, \dots, b_g$ of $H_1(\rsf_{\bp}, \Z)$ that verifies the following properties.

 \begin{enumerate}
 \item The collection of elements $a_1, \dots, a_{\graphgenus}$ forms a basis of the space of vanishing cycles $A \subseteq H_1(\rsf_{\bp}, \Z)$;

 \item[(2)] for any $k=1, \dots, n$, the collection of elements $a_{j}$ and $b_j$ for $j=g_{k-1}+1, \dots, g_k$ gives a basis of $H_1(C_{v_k}, \Z)$ in 
 \[H_1(C_{v_k}, \Z) \hookrightarrow H_1(S_0, \Z) \simeq H_1(\rsf_{\bp}, \Z) /A,\quad \text{  and }\]

 \item[(3)] the elements $b_1, \dots, b_{\graphgenus} \in \Bg \simeq H_1(G, \Z)$ form an admissible basis for the layered graph $G = (V,E, \pi)$, in the sense of the previous section. 
 \end{enumerate} 
In particular, Proposition~\ref{prop:ExistenceAdmissibleBasis} implies that an admissible symplectic basis for $\hcurve$ always exists.


\section{Period map} \label{sec:period}

The aim of this section is to describe the period map for the variation of Hodge structures $H_1(\rsf_t)$ for $t \in B^*$, for the family $\pr \colon \rsf^\ast\to B^\ast$ described in the previous section. We refer to~\cite{ABBF, Hof84} for more details. 

Recall first that $B=\Delta^{3g-3}$ and that $B^* = B \setminus \cup_{e\in E} D_e$. Shrinking the polydisk if necessary and making a choice of local parameters around $0$ for $D_e$, we can write  $B^* \simeq (\Delta^*)^E \times \Delta^{3g-3 -|E|}$. Let $\widetilde {B^*}$ be the universal cover of $B^*$. We get an isomorphism $\widetilde{B^*} \simeq \H^E \times \Delta^{3g-3 -|E|}$ with $\H$ the Poincar\'e half-plane. The map 
\[ \widetilde{B^*} \to B^*\] is given by sending $\zeta_e\in \H$ to $\exp(2\pi \imi \zeta_e) \in \Delta^*$.

Denote by $\widetilde{\rsf}^*$ the family of Riemann surfaces over $\widetilde {B^*}$ obtained by pullback of the family $\rsf^*\to B^*$. In the following, we will use the notations of the preceding section. We fix a point $\bp \in 
B^\ast$ together with a lift $\tilde \bp \in \widetilde{B^*}$. We also choose an admissible symplectic basis
\begin{equation*}
a_1,\dotsc,a_{g},b_1,\dotsc,b_{g} \in H_1(\rsf_{\bp},\Z) =
\basisa \oplus \basisb.
\end{equation*} 
Recall from Section~\ref{sec:admissible_basis_standard} that the space of vanishing cycles $A$ is generated by
$a_1,\dots, a_\graphgenus \in A$, and $b_1, \dots, b_\graphgenus$ generate $H_1(\rsf_{\bp},
\Z)/A' \simeq H_1(G, \Z)$ as in \eqref{eq:aux}, 
and the corresponding subsequent interval blocks of size $\genusfunction(v_k)$ give a symplectic basis of $H_1(C_{v_k}, \Z)$ for $k=1, \dots, n$.

Let $\mathcal H$ be the local system on $B^*$ given by the first integral homology of the fibers of $\rsf^* \to B^*$. Let $\widetilde{\mathcal H}$ be the local system on $\widetilde {B^*}$ obtained by the pullback of $\mathcal H$. We have $\widetilde{\mathcal H}_{\tilde t} = H_1(\widetilde{\rsf}_{\tilde t}, \Z) = H_{1}(\rsf_{t},\Z) $ for any $t\in B^*$ and $\tilde t\in \widetilde {B^*}$ over it. The local system $\widetilde{\mathcal H}$ being trivial, we can spread out the above symplectic basis at $\bp$ to a symplectic basis 
\begin{equation*}
a_{1,\tilde t}, \dots, a_{g,\tilde t}, b_{1,\tilde t}, \dots,
b_{g,\tilde t} 
\end{equation*}
of $H_{1}(\widetilde{\rsf}_{\tilde t},\Z)$ for $\tilde t\in
\widetilde {B^*}$ living above any point $t \in B^*$.

Since $\basisa$ is isotropic (see Proposition~\ref{prop:isotropic}) and contains the subspace of vanishing
cycles, the Picard-Lefschetz formula \eqref{eq:PL} implies that
the elements $a_{i,\tilde t}$ only depend on $t$ and not on $\tilde
t$. Thus, we will  denote them merely by $a_{i,t}$. The same remark actually applies for the basis elements $b_{j,\tilde t}$ with $j>\graphgenus$ (see again Proposition~\ref{prop:isotropic}). So the only non-trivial part of the symplectic basis regarding the monodromy action is $b_{1, \tilde t}, \dots, b_{\graphgenus, \tilde t}$. 

If there is no risk of confusion, we drop $\tilde t$, and simply use
$a_i$ and $b_i$ for the elements, having in mind the dependency on $t$ or $\tilde t$, with the above remarks.

The local system $\mathcal H$ on $B$ leads to an admissible variation of polarized Hodge structures. Denote by $\mathcal H_\C$ the vector bundle on $B^*$ associated to the locally free sheaf $\mathcal H \otimes_{\Z}\Ocal_{B^*}$ of $\Ocal_{B^*}$-modules.  The Hodge filtration on $\mathcal H_\C$ is given by $\filt^0 = \mathcal H_\C \supset \filt^1 = \pr_*\omega_{\rsf^*/B^*}$, the pushout of the relative dualizing sheaf $\omega_{\rsf^*/B^*}$. By the work of Deligne~\cite{Del70}, the holomorphic vector bundle
$\mathcal H_\C$ admits a canonical extension to a holomorphic vector bundle over $B$ that lives in $\jmath_*(\mathcal H_\C)$ where $\jmath\colon B^*\hookrightarrow B$.  We denote this holomorphic bundle by $\mathcal V$, and note that we have $\mathcal V \rest{B^*} = \mathcal H_\C$. Sections of $\mathcal V$ are locally of the form $\exp(-\sum_{e \in E} \zeta_e N_e)$ applied to a multi-valued section of $\mathcal H_\C$ over $U \cap B^*$, for an open subset $U$ in $B$~\cite[Section II.5]{Del70}. By nilpotent orbit theorem, the Hodge filtration $\filt^\bullet$
 can be extended holomorphically to a filtration $\filt^0=\mathcal V \supset \filt^1$. Since $\pr_*\omega_{\rsf^*/B^*}$ extends canonically to the pushout $\pr_*\omega_{\rsf/B}$ of the relative dualizing sheaf to $B$ (which is locally free of rank $g$), the sheaf of $\mathcal O_B$-modules associated to $\filt^1$ is simply  $\pr_*\omega_{\rsf/B}$, which is a coherent $\mathcal O_B$-module by functoriality of the dualizing sheaf on families. 
 
  We define a collection of holomorphic one-forms
$(\omega_{i})_{i=1, \ldots, g}$ on $\rsf^*$ by requiring that, for each $t\in B^*$, the forms $\omega_{i, t}\coloneqq
\omega_{i}\rest{\rsf_{t}}$, $i=1, \ldots, g$, verify the equations
\begin{equation} \label{eq:7}
 \int_{a_{i, t}}
\omega_{j,t}=\delta_{i, j}.
\end{equation}
These equations uniquely determine $\omega_{1, t},\dots, \omega_{g,t}$, $t\in B^*$. Moreover, by coherence property, $\omega_{j,t}$ extends holomorphically to $\rsf$. Although we do not need this here, we note that for each $t\in D \subset B$, the collection $(\omega_{j,t})_{j=1}^g$ verifies the same equations given above where for those cycles $a_i$ which become vanishing cycles for $t$, the corresponding equations in \eqref{eq:7} fix the residues of the restriction of $\omega_{j,t}$ on each component of the stable curve $\rsf_t$ at the nodes appearing on that component.  

 It follows that the period matrix for each fiber $\widetilde{\rsf}^*_{\tilde t}$ in the 
family  $\widetilde{\rsf}^* \to \widetilde {B^*}$ is given by
\[\Omega_{\tilde t} = \left(\,\int_{b_{i, \tilde t}}\omega_{j,t}\,\right)_{i,j=1}^g.\]

Let $\H_g$ be the Siegel domain defined as    
\begin{equation*}
\H_g \coloneqq  \bigl\{\text{ $g\times g$ complex symmetric matrix } \Omega\,\,\st\
\rm{Im}(\Omega)>0\,\bigr\}.
\end{equation*}

\begin{defi}\label{def:period_map} The \textit{period map} of the variation of  Hodge
  structures $H_1(\rsf_t, \Z)$ over $B^*$ is given by
\begin{align}\label{eq:periodmap}
{\widetilde{\Phi}} \colon \widetilde{B^*} &\longrightarrow \mathbb
H_{g} \\
\tilde{t} &\longmapsto \Bigl( \,\int_{b_{i, \tilde t}}\omega_{j,t}\,\Bigr)_{i,j=1}^g. \qedhere
\end{align}
\end{defi}

The group $\mathrm{Sp}_{2g}(\R)$ acts on $\H_g$ by 
$$\begin{pmatrix}A & B \\ C & D\end{pmatrix}\Omega = (A\Omega+B)(C\Omega+D)^{-1}. 
$$
The Siegel moduli space is  a
smooth Deligne--Mumford stack defined as the quotient  $\Acal_g= \mathrm{Sp}_{2g}(\Z)\,\bigl\backslash\, \H_g$ and 
 parametrizes principally polarized abelian
varieties of dimension $g$. As we explain below, the image of the monodromy map for the local system $\mathcal H$ lives in $\mathrm{Sp}_{2g}(\Z)$.

Since $\omega_{j,t}$ for $t\in D \subset B$ might have poles at nodes appearing in the fiber $\rsf_t$, the (multivalued) coordinates $\int_{b_{i,s}}\omega_{j,s}$ can tend to infinity when $s$ tends to $t$. This may happen when the specialization of $b_{i,s}$ to $\rsf_t$ goes through nodes in $\rsf_t$.  In order to describe the precise asymptotic of these coordinates, we need  to explicitly describe the action of the logarithm of monodromy maps $N_e$, $e\in E$, on the entries in \eqref{eq:periodmap}. 

In what follows, we view the entries of the period matrix as functions on $\widetilde{B^*}$, and note that $N_e = \lambda_e -\id$ acts on the space of functions $\varphi \colon \widetilde B^* \to \C$ by $N_e(\varphi)(\tilde t) = \varphi(\lambda_e(\tilde t)) - \varphi(\tilde t)$.

 As in \eqref{eq:vanishing_cycle}, for each edge $e$, we can write 
\begin{equation} \label{eq:2}
  a_{e}=\sum_{i=1}^{g}c_{e,i} a_{i}.
\end{equation}
Moreover, since $a_e \in A$, the coefficients  $c_{e,i}$ are zero for $i>\graphgenus$.

By the Picard-Lefschetz formula \eqref{eq:PL}, we deduce that 
\begin{align}
  N_{e}(b_{i})
  &=\langle a_e , b_i \rangle a_{e} = c_{e,i}a_{e} \label{eq:10}.
  \end{align}
The forms $\omega _{j}$ are defined
globally, so they are invariant under monodromy. The integral of these
forms with respect to the vanishing cycles is computed by
\begin{equation}\label{eq:mono-forms}
  \int_{a_{e}}\omega _{j}= \sum_{i=1}^\graphgenus c_{e,i} \int_{a_i}\omega_j = c_{e,j}.
\end{equation}

For each pair of integers $i, j \in [g]$, we view the integration of $\omega_j$ along $b_i$ as a function 
\begin{align*}
\int_{b_{i,\bullet}}\omega_j \colon \widetilde{B^*} \to \C, \qquad  \tilde t \mapsto \int_{b_{i,\tilde t}}\omega_j.
\end{align*}
 Applying \eqref{eq:10} and \eqref{eq:mono-forms}, we deduce that
\begin{align*}
  N_e\Bigl(\int_{b_{i}}\omega_{j}\Bigr) =  \int_{\lambda_e(b_{i})}\omega_{j} - \int_{b_{i}}\omega_{j} = \int_{N_e(b_{i})}\omega_{j} = \langle a_e ,b_i \rangle\int_{a_e}\omega_{j} =
    c_{e,i}\,c_{e,j}.
  \end{align*}
We introduce the $g\times g$ matrices $\widetilde M_{e} = \Bigl(\widetilde M_{e}(i,j)\Bigr)_{i,j=1}^g$  by
\begin{gather}\label{eq:matrixmtilde}
  {\widetilde{M}}_{e} (i,j)\coloneqq  \begin{cases} c_{e,i} \, c_{e,j} & \qquad \textrm{if $i,j \leq \graphgenus$}\\
  0 & \qquad \textrm{otherwise}. 
  \end{cases}
\end{gather}
By Proposition~\ref{prop:matricial}, the matrix ${\widetilde{M}}_{e}$ is the
$\graphgenus \times \graphgenus$ matrix $M_{e}$ from \eqref{eq:1}, which we have filled up with zeros to a $g \times g$
matrix. In particular, the entries of this matrix are all integers.

Altogether, we have proved that the logarithm $N_e$ of (the monodromy map) $\lambda_e$ is given by the following element of the Lie
algebra $\mathrm{sp}_{2g}(\Z)$
\begin{equation}\label{eq:logmono}
N_{e}=
  \begin{pmatrix}  0 & {\widetilde{M}}_{e}\\
     0 & 0 
     \end{pmatrix}. 
\end{equation}

It follows that the image of the monodromy lives in $\mathrm{Sp}_{2g}(\Z)$ and the
map $\widetilde{\Phi}$ descends to $B^*$, resulting in the following commutative diagram.
\begin{equation}
\label{diag-phi}
\begin{CD} 
{\widetilde{B^*}} @>\widetilde\Phi>> 
\H_{g} \\
@VVV @VVV \\
B^* @>\Phi>> 
\mathrm{Sp}_{2g}(\Z)\bigl\backslash \H_{g}.
\end{CD}
\end{equation}
The map $\Phi$ sends a point $t \in B^*$ to the point of the Siegel moduli space $\Acal_g$ corresponding to the Jacobian $\mathrm{Jac}(\rsf_t)$.

\subsection{Asymptotics of the period map}
\label{sec:asympt-peri-map}
The asymptotics of the period map can be described by the nilpotent orbit
theorem for variations of polarized Hodge structures.

 We denote the coordinate variables corresponding to the edges as $
 z_E$, and write the coordinates of any point $t$ of $B^*$ as $
z_E \times z_{E^c}$. The coordinates in the
 universal cover $\widetilde{B^*}$ 
 will be denoted by $\zeta_{e}$, $e \in  [N]$. In these coordinates, the projection $\widetilde{B^*}\to B^*$ is
 given by
 \begin{equation}
   \label{eq:16}
   z_{e}=
   \begin{cases}
     \exp(2\pi \imi \zeta_{e}),&\text{ for }e\in E,\\
     \zeta_{e},&\text{ for }e\not \in E.
   \end{cases}
 \end{equation}
In these coordinates, the action of $\lambda_e$ on $\widetilde{B^*}$, for $e\in E$, corresponds to $\zeta_e \mapsto \zeta_e+1$.

The \emph{twisted period map} $\widetilde\Psi$ on $\widetilde {B^*}$ given by
\begin{equation}
  \label{eq:18}
\widetilde\Psi(\tilde t)=\exp(-\sum_{e \in E} \zeta_e N_e)\widetilde\Phi(\tilde t)
\end{equation}
takes values in the compact dual $\Check{D}$ of the period domain $D = \H_g$, which is 
(essentially) a flag
variety parametrizing filtrations $\filt^\bullet\mathbb C^{2g}$ satisfying the conditions of being a Hodge filtration of weight $1$, with Hodge numbers $h^{1,0}=h^{0,1}=g$~\cite[Section 3]{Sch73}. The space $\check{D}$ contains $\H_g$ as
an open subset and the action of the symplectic group $\mathrm{Sp}_{2g}(\Z)$ can be extended to $\check D$. 

The map $\widetilde \Psi$ is invariant under the transformation $\zeta_e \mapsto \zeta_e+1$ and so descends to a map $\Psi \colon B^* \to \check{D}$.   We have the following result \cite{Sch73, CKS86}.

\begin{thm}[Nilpotent orbit theorem] \label{thm:nilpotentorbit} 
After shrinking the radius of $\Delta $ if necessary, the map $\Psi$ extends
to a holomorphic map
\begin{displaymath}
  \Psi \colon B \longrightarrow \check{D}. 
\end{displaymath}
Moreover, there exists a constant $T_0 >0$ such that
\begin{displaymath}
  \exp(\sum_{e\in E} \zeta_eN_e)\Psi(t)\in \H_g
\end{displaymath}
for all  $t\in B$, provided that $\Im(\zeta_{e})\ge T_{0}$ for any $e\in E$.

Furthermore, there are constants $C, \beta >0$ so that the following estimate
\[\mathrm{dist}\Bigl(\widetilde \Phi(\tilde s) ,  \exp(\sum_{e\in E} \zeta_e(\tilde s)N_e)\Psi(t)\Bigr) \leq C \sum_{e\in E} \Im(\zeta_e(\tilde s))^\beta \exp(-2\pi \Im(\zeta_e(\tilde s)))\]
holds for any $s\in B^*$, for the point $t$ with $z_{E}(t)=0$ and $z_{E^c}(t) =z_{E^c}(s)$, and any $\tilde s  \in \widetilde {B^*}$ above $s$ satisfying $\Im(\zeta_{e}(\tilde s ))\ge T_{0}$ for all $e \in E$.
\end{thm}

In what follows we denote the period map $\widetilde\Phi$ from Definition~\ref{def:period_map} as
\begin{align*}
{\widetilde \Phi} \colon {\widetilde{B^*}} &\longrightarrow \mathbb
H_{g} \\
\tilde{t} &\longmapsto \Omega_{\tilde t} = \Bigl( \,\int_{b_{i, \tilde t}}\omega_{j,t}\,\Bigr)_{i,j=1}^g.
\end{align*}

\begin{prop} Notations as in the previous section, we have that
\[\widetilde \Psi(\tilde t) = \Omega_{\tilde t} - \sum_{e\in E} \zeta_e(\tilde t) \widetilde M_e.\] 
\end{prop}

It follows from the nilpotent orbit theorem that, shrinking the polydisk $B$ around the origin, if necessary, we can write $\widetilde \Psi(\tilde t) = \Lambda_t \in \H_g$ for any $t\in B^*$ and $\tilde t \in \widetilde{B^*}$ above $t$. Moreover, the family of matrices $(\Lambda_t)_{t \in B^\ast}$ can be extended to a family over $B$ by setting $\Lambda_t = \Psi(t)$ for the extension $\Psi \colon B \to \check D$.

We will be interested in the imaginary part of $\Omega_{\tilde t}$. By the description of the monodromy, $\Im(\Omega_{\tilde t})$ is invariant under monodromy and  descends to a map from $B^*$ to the space of positive definite symmetric $g\times g$ matrices. For a point $t\in B$, we denote this matrix by $\Im(\Omega_t)$.

The above results imply the following estimate for the imaginary parts. Given $s \in B^\ast$, we introduce the log parameters $\ell_e(s) \coloneqq  -\log|z_e(s)|$, $e\in E$.
 
\begin{thm} \label{thm:nilpotentorbit2} For $s \in B^\ast$ with $\ell_e(s)$, $e \in E$, large enough, the uniform estimate
\[\mathrm{dist} \Bigl(\Im(\Omega_s), \Im(\Lambda_t) + \frac{1}{2\pi} \sum_{e\in E} \ell_e(s) \widetilde M_e\Bigr) \leq C \sum_{e\in E} \ell_e(s)^\beta \exp(- \ell_e(s)) \]
holds, where $t \in B$ is the point with $z_E(t)=0$ and $z_{E^c}(t)=z_{E^c}(s)$.
\end{thm}

\begin{defi} We denote by $\Lambda_0$ the value of $\Psi$ at $0\in B$ and call it the (\emph{regularized}) \emph{limit period matrix}, or simply, the \emph{limit period matrix} at the origin in $B$. We have 
\[\Lambda_0 = \lim_{\substack{t\in B^*\\ t \to 0}} \bigl(\Omega_{\tilde t} - \sum_{e\in E} \zeta_e(\tilde t) \widetilde M_e\bigr)\]
where $\tilde t$ is any point of $\widetilde{B^*}$ living above $t$.
\end{defi}

\begin{remark} The limit period matrix $\Lambda_0$ at the origin $0 \in B$ depends on the choice of parameters, and is defined only up to a sum of the form $\sum_{e \in E} \lambda_e\widetilde M_e$.
\end{remark}

\subsection{Description of the limit period matrix $\Lambda_0$ at origin} As before, let $\pr \colon \rsf \to B$ be the family of Riemann surfaces over the polydisk $B$. The pullbacks $\pr^*(D_e)$ of divisors $D_e$ for $e\in E$ form a simple normal crossing divisor in $\rsf$. The dualizing sheaf $\omega_{\rsf/B}$ coincides with the sheaf $\Omega^1_{\rsf/B}\bigl(\log(\pr^*D)\bigr)$ of holomorphic 1-forms with logarithmic singularities along the divisor $\pr^*D$. Our forms $\omega_1, \dots, \omega_g$ form a basis of  $H^0\Bigl(B, \pr_*\omega_{\rsf/B}\Bigr)$. 

Let $a_1, \dots, a_g, b_1, \dots, b_g$ be an admissible symplectic basis of $H_1(\rsf_{\bp}, \Z)$. For $k\in [n]$, and for each pair $g_{k-1}+1\leq j \leq g_k$ (see \eqref{eq:DefineGi} for the definition of $g_k$), and each $i\in [g]$, we have 
\[\int_{a_i} \omega_j =\delta_{i,j}.\]
We infer that the restrictions
$\omega!_{g_{k-1}+1, 0}, \dots,\omega!_{g_{k},0}$ to $C_{v_k}$ form a basis of the space of holomorphic differentials on $C_{v_k}$. Moreover, the function $\int_{b_i} \omega_{j}$ is holomorphic on $B$ for any $j\in[g]$.

For each $k=0, 1, \dots, n$, let $I_k$ be the index set $I_k = \{g_{k-1}+1, \dots, g_k\}$ (we set $g_{-1}=0$). For a $g\times g$ matrix $P$, denote by $P[I_k, I_{k'}]$ the matrix with rows in $I_k$ and columns in $I_{k'}$.  We deduce the following theorem.
\begin{thm} The limit period matrix $\Lambda_0$ at the origin $t = 0$ has the following shape.
\begin{itemize}[leftmargin = 2em]
\item Let $1\le k \le n$. Then the matrix $\Lambda_0[I_k, I_k] \in \C^{\genusfunction(v_k)\times \genusfunction(v_k)}$ is the period matrix $\Omega_{v_k}$ of the component $C_{v_k}$ for the symplectic basis $a_j$, $b_j$, for $j \in I_k$, and the holomorphic forms $\omega_{j, 0}$, $j \in I_k$, restricted to $C_{v_k}$.
\item All matrices $\Lambda_0[I_k, I_{k'}] \in \C^{\genusfunction(v_k)\times \genusfunction(v_{k'})}$ for indices $1\le k, k' \le n$ with $k \neq k'$ are vanishing.
\end{itemize}
 \end{thm}
 In other words, the limit period matrix $\Lambda_0$ has the following form
\begin{equation} \label{eq:Lambda_0}
 \Lambda_0 = \left(
 \begin{matrix}
\Omega_G & * & * & * & \cdots & *\\
* & \Omega_{v_1} & 0 & 0 & \cdots & 0\\
* & 0 & \Omega_{v_2} & 0 & \ddots & 0\\ 
\vdots & \vdots & \vdots & \ddots & \ddots &0 \\
* & 0 & \cdots & \cdots & \Omega_{v_{n-1}} & 0 \\
* & 0 & \cdots & \cdots & 0 & \Omega_{v_n}
\end{matrix} \right)
\end{equation}
for some $\graphgenus\times\graphgenus$ matrix $\Omega_G$.

\begin{proof} By definition, $\Lambda_0$ is the limit of $\Omega_{\tilde t} - \sum_{e\in E} \zeta_e(\tilde t) \widetilde M_e$ as $t\in B^*$ tends to 0 and $\tilde t$ is a point of $\widetilde{B^*}$ above $t \in B^\ast$. By definition of the matrix $\widetilde M_e$ \eqref{eq:matrixmtilde}, the entries of $\Lambda_0[I_k, I_{k'}]$ for $k,k'\in [n]$ correspond to the limit of the corresponding entries in $\Omega_{\tilde t}$. By what preceded, for $k=k'$, these are precisely the entries of the period matrix $\Omega_{v_k}$ of the component $C_{v_k}$. This proves the first statement, Consider now the case $k\neq k'$. For any $t\in B^*$, we have $\int_{a_{i,t}}\omega_{j,t} =0$ for any $i\in I_k$ and $j\in I_{k'}$. Passing to the limit, we get $\int_{a_{i,0}} \omega_{j,0}=0$, that is the restriction of $\omega_{j,0}$ to the component $C_{v_{k}}$ has vanishing integrals along all the elements $a_{i,0}$, $i\in I_{k}$. Since $a_{i,0}, b_{i,0}$, $i\in I_{k}$, form a symplectic basis of $C_{v_{k}}$, and the restriction of $\omega_{j,0}$ to $C_{v_{k}}$ is holomorphic, we infer that $\omega_{j,0}$ vanishes on $C_{v_{k}}$, and so all the entries in $\Lambda_{0}[I_{k}, I_{k'}]$ are zero.  
\end{proof}


\section{Generic continuity}\label{sec:generic}

Let $S_0$ be a stable Riemann surface with dual graph $G = (V, E, \genusfunction)$, and denote by $\rsf\to B$ the analytic versal deformation space and the versal family of Riemann surfaces over it. We fix a base point $\bp \in B^*$.

Consider the family of canonically measured hybrid curves $(\rsf_\thy^\hyb,\mu_\thy)_{\thy \in B^\hyb}$ over $B^\hyb$. That is, each hybrid curve $\rsf^\hyb_\thy$, $\thy \in B^\hyb$, is endowed with its canonical measure $\mu_\thy = \mu_\thy^\can$. 

The aim of this section is to prove the continuity of canonical measures when approaching hybrid points \emph{through the open subset $B^* \subset B^\hyb$}. More precisely, we show the following result.
\begin{thm} \label{thm:mainlocal_open} The family of canonically measured hybrid curves $(\rsf_\thy^\hyb,\mu_\thy)_{\thy \in B^\hyb}$ over the hybrid space $B^\hyb$ is continuous through the open subset $B^* \subset B^\hyb$. That is, for every continuous function $f \colon \rsf^\hyb \to \R$, the function $F \colon B^\hyb \to \R$ defined by integration along fibers
\begin{align*}
&F(\thy) \coloneqq  \int_{\rsf^\hyb_{\thy}} f_{|_{\rsf^\hyb_{\thy}}} \, d\mu_\thy, \qquad  \thy \in B^\hyb,
\end{align*}
satisfies the continuity condition
\[
\lim_{\substack{t \to \thy \\ t \in B^\ast}} F(t) = F(\thy)
\]
for all points $\thy \in B^\hyb$.
\end{thm}
\begin{remark}\label{rem:first-reduction} Without loss of generality,  the proof of the continuity can be reduced to the case where $\thy=(t,x)$ with $t=0 \in B$. This allows to simplify the presentation in the following.
\end{remark}

In the following, we fix an ordered partition $\pi=(\pi_1, \dots, \pi_r)$ on $E$, and a point of the form $\thy=(0,x)$ in the hybrid stratum $D_\pi^\hyb =  \inn D_\pi \times \inn \sigma_\pi$. We also fix an admissible basis for $H_1(S_{\bp}, \Z)$ with respect to the layering $\pi$.

Let $(\omega_i)_{i=1}^g$ be the corresponding family of holomorphic one-forms on $\rsf$ (see Section~\ref{sec:period}). For each point $t \in B^\ast$, we write the canonical measure on the smooth Riemann surface $\rsf_t$ as
\begin{equation} \label{def:decomposition}
\mu_t = \sum_{i,j=1}^g \mu_{i,j,t},
\end{equation}
where, for $i,j\in[g]$, the complex-valued measure $\mu_{i,j,t}$ on $\rsf_t$ is given by 
\[\mu_{i,j,t} \coloneqq  \frac {\imi}{2}\Im(\Omega_t)^{-1}(i,j) \, \omega_i \wedge \bar \omega_j.\]

We prove Theorem~\ref{thm:mainlocal_open} by describing the limit of measures $\mu_{i,j,t}$, as $t \in B^\ast$ approaches the point $\thy=(0,x)$. There will be three regimes, each treated in a separate section below. 

For notational convenience, we introduce the following index sets  
\begin{equation}\label{eq:ik}
I_k \coloneqq  \{g_{k-1}+1, \dots, g_k\}, \qquad k=1, \dots, n,
\end{equation}
where $g_k = \graphgenus+ \sum_{j=1}^{k} \genusfunction(v_j)$ is given as in \eqref{eq:DefineGi}, with $\genusfunction(v_j)$ the genus of the Riemann surface component $C_{v_j}$ of $\rsf^\hyb_\thy$, for $j \in [n]$. 
We will establish the following limiting behavior (see below for more precise statement):
\begin{enumerate}[leftmargin = 2em]
\item If $\graphgenus< i,j $, then, $\mu_{i,j,t}$ converges
\begin{itemize}[leftmargin = 1em]
\item[-] either, to a piece of the canonical measure on the Riemann surface component $C_{v_k}$ of $\rsf^\hyb_\thy$, according to the analogue decomposition to~\eqref{def:decomposition}; this happens in the case $i,j$ belong to the same block $I_k$ for $k \in [n]$,
\item[-] or, to zero; when $i,j$ belong to two different blocks among $I_k$s, $k \in [n]$,
\end{itemize}

\item If $1 \le i,j \le \graphgenus$, then $\mu_{i,j,t}$ converges to a piece of the canonical measure on the underlying tropical curve of the hybrid curve $\rsf^\hyb_\thy$, according to the decomposition deduced from Theorem~\ref{thm:cmgraphscycles}.

\item If either, $i\leq \graphgenus <j$, or, $j\leq \graphgenus <i$, then $\mu_{i,j,t}$ converges to zero.
\end{enumerate}

\subsection{Inverse lemma} The proof of Theorem~\ref{thm:mainlocal_open} requires an understanding of the asymptotics of the inverse of period matrices. Our key tool is the \emph{inverse lemma} stated below. Its proof is given in Section \ref{sec:inverse_lemma_proof}

\begin{lem}[Inverse lemma] \label{lem:inverse_lemma} Let $X$ be a topological space and fix a point $ \thy \in X$. Let $y_1, \dots, y_r$ be  a collection of non-vanishing complex valued functions on $X \setminus \{\thy\}$ such that
\begin{equation} \label{eq:inverse_lemma_hyp}
	\lim_{t \to \thy} \frac{y_{k+1}(t)}{y_k(t)} = 0, \qquad k = 1, \dots, r-1.
\end{equation}
Let $M \colon  X \setminus \{\thy\} \to \C^{n \times n}$ be a matrix-valued function. Assume that $M(t)$ has an $(r,r)$ block decomposition of the form
\[
	M(t) = \Big(\mA_{kl}(t) \Big )_{1 \le k,l \le r},
\]
where, as $t$ converges to $\thy$ in $X$, the blocks $\mA_{kl} \colon X \setminus \{\thy\} \to \C^{n_k \times n_l}$ are asymptotically given by
\[
	\mA_{kl}(t) = \ssub{y}!_{\max\{k,l\}}(t) \big( \widehat{\mA}_{kl} + o(1) \big )
\]
for matrices $\widehat{\mA}_{kl} \in \C^{n_k \times n_l}$, and all the diagonal matrices $\widehat{\mA}_{kk}$, $k=1,\dots,r$, are invertible.

Then, the matrix $M(t)$ is invertible for $t$ in a punctured neighborhood of $\thy$. Moreover, as $t$ converges to $\thy$, the inverse $M(t)^{-1}$ has an $(r,r)$ block decomposition with asymptotic behavior 
\[
	M(t)^{-1} = \Big( \ssub{y}!_{\min\{k,l\}}(t)^{-1} \big (\mB_{kl} + o(1) \big) \Big)_{1 \le k,l \le r}
\]
for some matrices $\mB_{kl} \in \C^{n_k \times n_l}$. The matrices in the asymptotics of on-diagonal blocks are given by
\[
	\mB_{kk} = \widehat{\mA}_{kk}^{\, -1}, \qquad k=1, \dots r.
\]
\end{lem}

\subsection{The inverse of the period matrix}
In this section we apply the inverse lemma~\ref{lem:inverse_lemma} to $\Im(\Omega_t)$, the imaginary part of the period matrix for the Riemann surface $\rsf_t$, $t \in B^\ast$, and describe the asymptotic behavior of $\Im(\Omega_t)^{-1}$ as $t \in B^\ast$ converges to $\thy$.

{
In order to apply Lemma \ref{lem:inverse_lemma}, we need a detailed description of the asymptotics of $\Im(\Omega_t)$. Recall from Theorem \ref{thm:nilpotentorbit2} that, as $t\in B^\ast$ converges to the origin $0 \in B$,}
\begin{equation}  \label{eq:first_asymptotics}  \mathrm{dist} \Bigl(\Im(\Omega_t), \Im(\Lambda_0) +  \frac{1}{2\pi} \sum_{e\in E} \ell_e(t) \widetilde M_e\Bigr) \to 0.
\end{equation}
 Thus, choosing an admissible basis for $H_1(\rsf_{\bp}, \Z)$ as in Section \ref{sec:admissible_basis_hybrid_curves}, $\Im(\Omega_t)$ has the following form 
 \begin{equation*} \Im(\Omega_t)  
= \left(
 \begin{matrix}
\Im(\Omega_G) + M_{\ell(t)} & *  & * & \cdots & *\\
* & \Im(\Omega_{v_1})  & 0 & \cdots & 0\\
\vdots & \vdots  & \ddots & \ddots &0 \\
* & 0 & \cdots  & \Im(\Omega_{v_{n-1}}) & 0 \\
* & 0 &  \cdots & 0 & \Im(\Omega_{v_n})
\end{matrix} \right)+ o(1),  
\end{equation*}
as $t \in B^\ast$ tends to the origin $0 \in B$ in the standard topology on $B$. Here, $M_{\ell(t)} = \sum_{e\in E} \ell_e(t) M_e$ is the matrix defined in~\eqref{eq:M_l} for the edge lengths $\ell_e(t) =-\log|z_e(t)|$, $e \in E$.

Next, we rewrite $\Im(\Omega_t)$ as a square block matrix with $r +1$ blocks in each row/column. Here, $r$ is the rank of the ordered partition $\pi = (\pi_i)_{i=1}^r$ with $\thy = (0,x) \in D_\pi^\hyb$. Using the index decomposition $[h] = J_\pi^1 \sqcup \dots \sqcup J^{r}_\pi$ from \eqref{eq:ordered_partition_edges}, we set
\[ 
	J^{r+1}_\pi \coloneqq  \{ h+1, \dots, g \}
\] such that altogether, we get a partition
\begin{equation} \label{eq:decomposition_indices}
[g] = J^1_\pi \sqcup \dots \sqcup J^{r+1}_\pi.
\end{equation}
For a matrix $\mA \in \C^{g \times g}$, we write $\mA_{kl}\coloneqq  \mA\rest{J^k_\pi \times J^l_\pi}$, $1 \le k, l \le r+1$, for its $(k,l)$-th block with respect to the decomposition \eqref{eq:decomposition_indices}. We stress that the notation $\mA(i,j)$ is used for the $(i,j)$-th matrix entry.

Let $\mA(t) \coloneqq  \Im(\Omega_t)$. We obtain the following $(r+1, r+1)$ block decomposition of $\Im(\Omega_t)$,
\[
\Im(\Omega_t)  = \mA(t) = \Big ( \mA_{kl}(t) \Big )_{1 \le k,l \le r+1} \quad \text{where} \quad \mA_{kl}(t) \coloneqq   \Im(\Omega_t)_{kl}  = \Im(\Omega_t)\rest{J^k_\pi \times J^l_\pi}. 
\]

To describe the asymptotics of the blocks $\mA_{kl}(t)$, we need the following functions $y_1, \dots y_{r+1}$ on $B^\ast$. For $k = 1, \dots, r$, define $y_k\colon B^\ast \to \R_+$ by
\begin{equation}\label{eq:yk}
	y_k(t) \coloneqq  - \frac{1}{2\pi} \sum_{e \in \pi_k} \log|z_e(t)| = \frac{1}{2\pi} \sum_{e \in \pi_k} \ell_e(t) , \qquad t \in B^\ast,
\end{equation}
and set $y_{r+1}(t) \equiv 1$ on $B^\ast$. By the topological properties of $B^{\hyb}$, we have
\[
	\lim_{\substack{t \to \thy \\ t \in B^\ast }} \frac{y_{k+1}(t)}{y_k(t)} = 0, \qquad k = 1, \dots, r.
\]
We obtain the following asymptotic behavior when $t \in B^\ast$ tends to $\thy$ in $B^{\hyb}$:

\begin{itemize}[leftmargin = 2em]
\item Assume first that $1\le k \le l \le r$. Then,
\begin{align*}
	\mA_{kl}(t) = \mA_{lk}(t)^\transpose &= \Im(\Lambda_0)_{kl} + \frac{1}{2\pi} \sum_{e \in E} \ell_e(t) (M_e)_{kl} + o(1) \\
	&= \frac{1}{2\pi} \sum_{e \in \pi_l \sqcup \dots \sqcup \pi_r} \ell_e(t) (M_e)_{kl} + O(1) = y_l(t) \Big( \sum_{e \in \pi_l} \frac{\ell_e(t)}{2 \pi \, y_l(t)} (M_e)_{kl} + o(1) \Big) \\
	& = y_l(t) \Big( \sum_{e \in \pi_l} x_e (M_e)_{kl} + o(1) \Big), 
\end{align*}
where we have used that $(M_e)_{kl}  = 0$ for all $e \in \pi_1 \sqcup \dots \pi_{l-1}$. The latter holds true since our fixed basis of $H_1(\rsf_{\bp}, \Z)$ is admissible. Note also that, as $t$ tends to $\thy =(0,x)$, we have $-\log|z_e(t)| = 2\pi y_l(t) (x_e +o(1))$ for $e \in \pi_l$ by the definition of convergence in $B^\hyb$.
\item On the other hand, if $k= r+1$ or $l=r+1$, then by \eqref{eq:first_asymptotics},
\begin{align*}
	\mA_{kl}(t) &= \Im(\Lambda_0)_{kl} + o(1). 
\end{align*}
\end{itemize}

The inverse lemma~\ref{lem:inverse_lemma} now allows to describe $\mA(t)^{-1}=\Im(\Omega_t)^{-1}$. Denote by
\begin{equation} \label{eq:graded_matrix}
M^k_{\pi,x} \coloneqq \sum_{e\in \pi_j}  x_e M_e \, \in \R^{h^{k}_\pi \times h^{k}_\pi}
\end{equation}
the matrix \eqref{eq:M_l} on the $k$-th graded minor $\grm_\pi^k(G)$ of $G$, equipped with the edge lengths in the simplicial part $x$ of the point $\thy=(0,x)$, and with respect to the basis of $H_1(\grm^k_\pi(G), \Z)$ given by the contracted cycles $\proj_k(b_i)$, $i \in J^k_\pi$ (see Section \ref{sec:admissible_basis_hybrid_curves} for details).

\begin{thm} \label{thm:inverse_period_asymptotics}
Let $\Im(\Omega_t) \in \R^{g \times g}$, $t \in B^\ast$, be the imaginary part of the period matrix of $\rsf_t$. Consider the block matrix decomposition of its inverse $\Im(\Omega_t)^{-1}$. If $t \in B^\ast$ converges to $\thy = (0,x)$ in $B^{\hyb}$, then
\begin{equation}
	\Im(\Omega_t)^{-1} = \Big ( \ssub{y}!_{\min\{k,l\}}(t)^{-1} (\mB_{kl} + o(1)) \Big )_{1 \le k, l \le r+1}
\end{equation}
 for matrices $\mB_{kl} \in \R^{h^{k}_\pi \times h^{l}_\pi}$ (where $h^{r+1}_\pi \coloneqq  g -h$). Moreover,
 \[
	\mB_{k k} =  ( M_{\pi, x}^k )^{-1}, \qquad k =1, \dots, r,
\]
and 
\[
	\mB_{r+1, r+1} =  \left ( \begin{matrix}
  \Im(\Omega_{v_1})^{-1} & 0 & 0 & \cdots & 0\\ 
  0 & \Im(\Omega_{v_2})^{-1} & 0 & \ddots & 0\\ 
  \vdots & \vdots & \ddots & \ddots &0 \\
  0 & \cdots & \cdots & \Im(\Omega_{v_{n-1}})^{-1} & 0 \\
 0 & \cdots & \cdots & 0 & \Im(\Omega_{v_n})^{-1}
\end{matrix} \right).
\]
\end{thm}

\begin{proof}
The claim is an immediate consequence of the preceding discussion and Lemma \ref{lem:inverse_lemma}. Indeed, the structure of $\mB_{r+1, r+1}$ follows from \eqref{eq:Lambda_0} and it only remains to notice that  $\sum_{e \in \pi_k} x_e (M_e)_{kk} = M_{\pi, x}^k$ for each $k=1,\dots, r$. 
\end{proof}

\subsection{Continuity I} \label{ss:ctI}
For each $k=1, \dots, n$, consider the index set $I_k = \{g_{k-1}+1, \dots, g_k\}$ defined in \eqref{eq:ik}. We prove the following result.

\begin{thm}\label{thm:cont1} Assume that $t \in B^\ast$ converges to $\thy=(0,x)$. Then the following holds true for each pair of indices $(i,j)$ with  $\graphgenus+1\leq i, j \leq g$.

\begin{itemize}[leftmargin = 2em]
\item If $i, j$ belong to the same set $I_k$ for some $1\leq k\leq n$, then $\mu_{i,j,t}$ converges to the measure 
\[
\mu_{i,j,\thy} \coloneqq  \frac {\imi} {2} \Im(\Omega_{v_k})^{-1}  (i,j) \, \omega_i \wedge \bar \omega_j
\]
on {$\rsf^\hyb_{\thy}$}, supported on $C_{v_k} \subset {\rsf^\hyb_{\thy}}$. Here $\Omega_{v_k}$ denotes the period matrix of the Riemann surface component $C_{v_k}$ for the symplectic basis $\{a_i,b_j\}_{i,j\in I_k}$ (see Section~\ref{sec:admissible_basis_hybrid_curves}).

\item If $i, j$ belong to distinct sets $I_k \neq I_l$, then $\mu_{i,j,t}$ converges to the zero measure on {$\rsf_{\thy}^\hyb$}.
\end{itemize}
\end{thm}
\begin{proof}
The relative holomorphic forms $\omega_i, \omega_j$ on $\rsf^* \to B^*$ extend holomorphically to $\rsf\to B$. In particular, the complex valued-measure $\omega_i\wedge \bar \omega_j$ extends continuously over all points of $B$. It follows that when $t \in B^\ast$ approaches the hybrid point $\thy=(0,x)$, the family of measured spaces $(\rsf_t, \omega_{i, t}\wedge \bar \omega_{j,t})$ converges to $(\rsf^\hyb_{\thy}, \omega_{i}\wedge \bar \omega_{j})$.  This measure is supported in $C_{v_k}$, for $i\in I_k$, and is zero if $j\not\in I_k$.

It will be thus enough to prove that
\begin{equation}\label{eq:limit1}
\lim_{\substack{t \to \thy \\ t \in B^\ast}} \Im(\Omega_t)^{-1} (i,j) = \Im(\Omega_{v_k})^{-1} (i,j)
\end{equation}
if $i$ and $j$ belong to the same interval $I_k$. This follows immediately from Theorem \ref{thm:inverse_period_asymptotics}.
\end{proof}

\subsection{Continuity II} \label{ss:ctII}
Consider a pair of indices $1\leq i,j \leq \graphgenus$. Let $M^k_{\pi,x} \in \R^{h_\pi^k \times h_\pi^k}$, $k=1, \dots r$,  be the matrices from~\eqref{eq:graded_matrix} and $[h] = J_\pi^1 \sqcup \dots \sqcup J^{r}_\pi$ the index decomposition from~\eqref{eq:ordered_partition_edges}.

\begin{thm}\label{thm:cont2} Assume that $t \in B^\ast$ converges to $\thy=(0,x)$. Then, the following holds true for each pair of indices $(i,j)$ with $1 \le i,j \le \graphgenus$.

\begin{itemize}[leftmargin = 2em]
\item If $i,j$ belong to the same set $J_\pi^k \subset [h]$, then $\mu_{i,j,t}$ converges to a measure $\mu_{i,j, \thy}$ on {$\rsf^\hyb_\thy$} supported on the intervals of {$\rsf^\hyb_\thy$}. More precisely, the limit measure is given by 
\[\mu_{i,j,\thy} = \sum_{e\in \pi_k}( M^k_{\pi,x})^{-1} (i,j) \, \gamma_i(e) \gamma_j(e) \, d\theta_e\]
where $\gamma_i = \proj_k(b_i)$ and $\gamma_j = \proj_k(b_j)$ are the cycles in $H_1(\grm^k_\pi(G), \Z)$ corresponding to the elements $b_i, b_j \in \Bg\simeq H_1(G, \Z)$ of the fixed admissible basis, $\proj_k$ is the projection map $G^{k}_\pi\to \grm_\pi^k(G)$, and $d\theta_e$ is the uniform Lebesgue measure on the interval $\Ical_e \subseteq{\rsf^\hyb_\thy}$ representing the edge $e$. 

\item If $i,j$ belong to distinct sets $J_\pi^k \neq J_\pi^l$, then $\mu_{i,j,t}$ converges to the zero measure on {$\rsf_{\thy}^\hyb$}.
\end{itemize}
\end{thm}

The rest of this section is devoted to the proof of this theorem.

\subsubsection{The behavior near singular points} \label{sec:behavior_singular_points}
Fix a small neighborhood $U_0$ of the origin $0 \in B$. Let $e =uv$ be an edge of the graph $G$, and consider the singular point $p^e(0)$ of the fiber $\rsf_0=S_0$. We find a small neighborhood $U_{e}$ of $p^e(0)$ in $\rsf$ lying above $U_0$ and put coordinates $\underline{z} = \Bigl((z_i)_{i \neq e}, z^e_u, z^e_v\Bigr)$ on $U_{e}$ with the equation $z^e_u z^e_v =z_e$, using a standard coordinate neighborhood $(U_e, z)$ (see \eqref{eq:standard_coordinates} for details).

Since $\omega_i$ and $\omega_j$ are global sections over $\rsf$ of the dualizing sheaf $\omega_{\rsf/B}$, locally in a small neighborhood $U_{e}$ as above, we can write
\begin{equation} \label{eq:1forms}
\omega_i  = \frac 1{2\pi \imi}f_i(\underline z) \frac {d z_u^e}{z_u^e}, \qquad \omega_j = \frac 1{2\pi \imi}f_j(\underline z) \frac {d z_u^e}{z_u^e}
\end{equation}
for holomorphic functions $f_i$ and $f_j$ on $U_{e}$.

For the vanishing cycle $a_e$, we have $\int_{a_e} \omega_i =\gamma_i(e)$ and $\int_{a_e} \omega_j =\gamma_j(e)$, see~\eqref{eq:mono-forms}. It follows by the residue formula that $f_i(\underline z) = \gamma_i(e)$ if $z^e_u = z^ e_v = 0$. The same holds true for $f_j$.

Now, we write $z_u^e = \exp(2\pi \imi\zeta^e_u)$, and pass to the polar coordinates $\ell^e_u \coloneqq - \log|z^e_u|$ and $\tau^e_u \coloneqq\operatorname{Re}(\zeta^e_u)$ on $U_{e} \setminus \, \pr^{-1}(D_e)$. Note that in particular, $\ell^e_u+ \ell^e_v = \ell_e(t) = -\log|z_e(t)|$.

 In these coordinates, we have
\[ \frac {d z_u^e}{z_u^e} = - d\ell^e_u + (2\pi \imi)d\tau^e_u.\]

It follows that on $U_{e} \setminus \, \pr^{-1}(D_e)$,
\[\omega_i \wedge \bar \omega_j =  \frac{1}{\pi \textrm{i}} \, f_i(\underline z) \bar f_j(\underline z) d\tau^e_u \wedge d\ell^e_u.\]

Let now $t \in B^\ast$ be a point close to $\thy$ in $B^{\hyb}$. Restricting the measure $\mu_{i,j,t}$ to $U_{e} \cap \rsf_t$, we get the expression 
\[\ \mu_{i,j,t} = \frac{1}{2 \pi} \Im(\Omega_t)^{-1} (i,j) f_i(\underline z) \bar f_j(\underline z) d\tau^e_u \wedge d\ell^e_u \qquad \text{on $U_{e} \cap \rsf_t$}.\]

Assume now that $i$ and $j$ belong to the index sets $J_\pi^k$ and $J_\pi^l$ with $1 \le k \le l \le r$. Then, Theorem \ref{thm:inverse_period_asymptotics} implies that
\[
 \mu_{i,j,t} =  \frac{ f_i(\underline z) \bar f_j(\underline z) }{2 \pi \, y_k(t)} \Big( \mB_{kl}(i,j) + o(1) \Big)  d\tau^e_u \wedge d\ell^e_u \qquad \text{on } U_{e} \cap \rsf_t,
\]
and the $o(1)$-term tends to zero uniformly on $U_e \cap S_t$ as $t \in B^{\ast}$ tends to $\thy$.

Finally, suppose that the edge $e$ belongs to the $m$-th set $\pi_m$ of the ordered partition $\pi$. Normalizing the coordinates by the respective lengths, we get
\[ \mu_{i,j,t}  = \frac{y_m(t)}{y_k(t)} f_i(\underline z) \bar f_j(\underline z)  \Big( \mB_{kl}(i,j) + o(1) \Big)   d\tau^e_u \wedge d\Theta_e\]
in $U_{e} \cap \rsf_t$, where
\[\Theta_e  \coloneqq  \frac{1}{2 \pi} \frac{\ell_u^e}{ y_m(t)}  = \frac{\ell_u^e}{\sum_{\hat e \in \pi_m} \ell_{\hat e}(t)} .
\]
Moreover, $f_i$ is holomorphic on $U_e$ with $f_i(\underline z) = \gamma_i(e)$ for $z^e_u= z^e_v =0$, and hence
\[
	f_i(\underline z) = \gamma_i(e) + O(|z_u^e|) + O(|z_v^e|) \qquad \text{in } U_e.
\]
The same holds true for $f_j$. Altogether, we have shown that
\begin{equation} \label{eq:measure_expression_1}
	\mu_{i,j,t} = c_t( p ) \, d\tau^e_u \wedge d\Theta_e \qquad \text{in } \rsf_t \cap U_e.
\end{equation}\label{eq:measure_expression_2}
where the function $c_t \colon \rsf_t \cap U_e \to \C$ has the form
\begin{equation} \label{eq:measure_expression_2}
	c_t(p) = \frac{y_m(t)}{y_k(t)} \Big( \gamma_i(e) \gamma_j(e) + O(|z_u^e|) + O(|z_v^e|) \Big )   \Big( \mB_{kl}(i,j) + o(1) \Big),
\end{equation}
as $t \in B^\ast$ tends to $\thy$ with uniform estimates for the error terms. More precisely,  the $o(1)$ term tends to zero uniformly on $U_e \cap \rsf_t$ for $t \to \thy$ and the constant upper bounds in the $O(|z_u^e|)$ and $O(|z_u^e|)$ terms can be chosen independent of $t$.

\subsubsection{Proof of Theorem~\ref{thm:cont2}}
Let $f$ be a continuous function on $\rsf^{\hyb}$.  We are concerned with the limit behavior of
\begin{equation} \label{eq:limit_in_question}
	\int_{ \rsf_{t}} f_{|_{\rsf_{t}}} \, d\mu_{i,j,t} 
\end{equation}
as $t \in B^\ast$ converges to $\thy=(0,x)$ in $B^{\hyb}$. Using the same notation as in the preceding section, we consider small neighborhoods $U_e$ around the singular points of $\rsf_0 = S_0$ and choose coordinates as above. As before, we suppose that $i \in J^\pi_k$ and $j \in J^\pi_l$ with $1 \le k \le l \le r$. The case $1 \le l < k \le r$ can be treated by symmetry, using $\overline{\mu_{i,j,t}} = \mu_{j,i,t}$.

Outside the open neighborhoods $U_e$, the measures $\mu_{i,j,t}$ extend continuously by zero. Indeed, the measures $\omega_i \wedge \bar \omega_j$ extend continuously outside these open sets and $\Im(\Omega_t)^{-1}(i,j) \to 0$ as $t \to \thy$ by Theorem~\ref{thm:inverse_period_asymptotics}. More precisely, we have $\Im(\Omega_t)^{-1}(i,j) = O(1/y_k(t))$ by Theorem~\ref{thm:inverse_period_asymptotics} and $y_k(t) \to + \infty$ as $t \to \thy$. In particular, 
\begin{equation} \label{eq:lim_zero_holomorphic}
	\lim_{\substack{t \to \thy \\ t \in B^\ast}} \int_{ \rsf_t \setminus (\bigcup_e U_e \cap \rsf_t)} f_{|_{\rsf_{t}}} \, d\mu_{i,j,t}  = 0.
\end{equation}

Hence it remains to analyze the integrals of $f_{|_{\rsf_{t}}}$ over the sets $U_e \cap \rsf_t$, $e \in E$. Fix an edge $e \in E$ and suppose that $e$ belongs to $m$-th set $\pi_m$ of the ordered partition $\pi$.

If $m < k$, then $\gamma_i(e) = \gamma_j(e) = 0$ since our fixed basis of $H_1(\rsf_{\bp}, \Z)$ is admissible (see Section \ref{sec:admissible_layered}). In this case, $\omega_i$ and $\omega_j$ do not have logarithmic poles at the node $p^e(0)$ of $\rsf_0$. Using the same argument as in \eqref{eq:lim_zero_holomorphic}, we see that, for any $e \in \pi_m$ with $m<k$, we have 
\[
	\lim_{\substack{t \to \thy \\ t \in B^\ast}} \int_{ \rsf_t\cap  U_e } f_{|_{\rsf_{t}}} \, d\mu_{i,j,t}  = 0.
\]
Therefore, in the following, we may suppose that $m \ge k$.

Next, we compute the range of the coordinate $\Theta_e = \ell^e_u/(2 \pi \,   y_m(t))$ on $U_e \cap \rsf_t$. Suppose that the coordinates $z_u^e$ and $z_v^e$ on $U_e$ have range $0 \le |z_u^e|, |z_v^e| \le \varepsilon$, for $\varepsilon>0$. Taking into account that $z^e_u z^e_v = z_e(t)$ on $U_e \cap \rsf_t$, one readily computes that $\Theta_e$ has range the interval
\[
\Theta_e(U_e \cap \rsf_t) = [h_1(t), h_2(t)] \eqqcolon \Ical_{e,t}
\]
where for $t \in B^\ast$, we set
\[
h_1(t) \coloneqq \frac{- \log(\varepsilon) }{2 \pi \, y_m(t)}, \qquad h_2(t) \coloneqq \frac{ - \log( |z_e(t)| / \varepsilon) }{ 2 \pi \,  y_m(t)} = \frac{ \ell_e(t) +\log(\varepsilon)}{ 2 \pi \,  y_m(t)},
\]
with $\ell_e(t)=-\log|z_e(t)|$.

We then decompose $U_e \cap \rsf_t = U_{e,t}^1 \sqcup U_{e,t}^2$ as the disjoint union of the sets
\begin{align*}
	U_{e,t}^1 &\coloneqq \Big \{ p \in U_e \cap \rsf_t \, \st \, | z_u^e(p)| \ge \frac{1}{\ell_e(t)}  \text{ or }  | z_v^e(p)| \ge \frac{1}{\ell_e(t)} \Big \}, \textrm{ and }  \\
	U_{e,t}^2 &\coloneqq \Big \{ p \in U_e \cap \rsf_t \, \st \, | z_u^e(p)| < \frac{1}{\ell_e(t)}  \text{ and }  | z_v^e(p)| < \frac{1}{\ell_e(t)} \Big \}.
\end{align*}
 A direct verification shows that on $U^1_{e,t}$ and $U^2_{e,t}$, the coordinate $\Theta_e$ has range
\[
\Theta_e(U^1_{e,t})	= \bigl[ h_1(t),  \tilde h_1(t)  \bigr ] \cup \bigl[\tilde h_2(t) , h_2(t)  \bigr] \eqqcolon \Ical_{e,t}^1, \qquad \Theta_e(U^2_{e,t})	= \bigl (\tilde h_1(t),  \tilde h_2(t)  \bigr ) \eqqcolon \Ical_{e,t}^2
\]
where 
\begin{align*}
&\tilde h_1(t) \coloneqq   \frac{\log\big (\ell_e(t)\big ) }{2 \pi \,  y_m(t) }, &\tilde h_2(t) \coloneqq  \frac{ \ell_e(t)-\log\left(\ell_e(t)\right)}{ 2 \pi \,  y_m(t) }.
\end{align*}
By the definition of convergence in $B^\hyb$ and~\eqref{eq:yk}, as $t \in B^\ast$ tends to $\thy=(0,x)$ in $B^\hyb$, we get
\begin{equation} \label{eq:LimitCoordinateRange}
\lim_{\substack{t \to \thy \\ t \in B^\ast}} \tilde h_1(t) = \lim_{\substack{t \to \thy \\ t \in B^\ast}} h_1(t) = 0, \qquad \lim_{\substack{t \to \thy \\ t \in B^\ast}} \tilde h_2(t) = \lim_{\substack{t \to \thy \\ t \in B^\ast}} h_2(t) = x_e.
\end{equation}
Thus, in the limit $t\to \thy$, the sets $\Ical^1_{e,t}$ shrink to $\{0\} \cup \{x_e\}$, and the intervals $\Ical^2_{e,t}$ expand to the interval $(0, x_e)$. Using the boundedness of $f$ and the coefficient $c_t$, see \eqref{eq:measure_expression_2}, and using that $m \ge k$, this implies that
\[
	\lim_{\substack{t \to \thy \\ t \in B^\ast}} \int_{ U_{e,t}^1} f_{|_{\rsf_{t}}} \, d\mu_{i,j,t} = 0.
\]

It remains to understand the limit of the integrals over the sets $U_{e,t}^2$. 

We first observe that 
\begin{equation} \label{eq:limit_coordinates_integral}
	\lim_{\substack{t \to \thy \\ t \in B^\ast}}  \int_{U_{e,t}^2} f_{|_{\rsf_{t}}} \, d\tau_u^e \wedge d\Theta_e   = \int_{\mathcal{I}_e} f(\lambda) \, d \theta_e (\lambda)
\end{equation}
where $\Ical_e$ is the interval of length $x_e$ representing the edge $e \in E$ in the hybrid curve ${\rsf^\hyb_\thy}$ and $\theta_e$ denotes the uniform Lebesgue measure on $\Ical_e$. Indeed, we can explicitly write
\[
\int_{U_{e,t}^2} f_{|_{\rsf_{t}}} \, d\tau^e_u \wedge d\Theta_e = \int_{\tilde h_1(t)}^{\tilde h_2(t)} \int_0^1 f \left( (z_i)_{i \neq e}, \tau_u^e, \Theta_e \right) \, \,  d \tau_u^e \, d \Theta_e\,.
\]
Then, \eqref{eq:limit_coordinates_integral} is a direct consequence of the definition of the topology on $\rsf^\hyb$ (namely, Proposition~\ref{prop:HybridFamilyConvergence}) and~\eqref{eq:LimitCoordinateRange}.

In view of \eqref{eq:measure_expression_1}, it finally remains to describe the behavior of the function $c_t(p)$. If $k=l=m$, then it follows from \eqref{eq:measure_expression_2} that 
\[ 
	\lim_{\substack{t \to \thy \\ t \in B^\ast}}  \sup_{p \in U_{e,t}^2} \| c_t(p) - \gamma_i(e) \gamma_j(e)  (M^k_{\pi,x})^{-1} (i,j) \| = 0.
\]
Here, we use that $|z_u^e|, |z_v^ e| \le 1 / |\log(|z_e(t)|)|$ on $U_{e,t}^2$ and hence $z^e_u, z^e_v \to 0$ uniformly on $U^2_{e,t}$ for $t \to \thy$. (Remark that this argument fails for $U_{e,t}^1$ instead of $U_{e,t}^2$ and hence the above decomposition was necessary.)
The above implies in turn that
\begin{equation*}
	\lim_{\substack{t \to \thy \\ t \in B^\ast}}  \int_{U_{e,t}^2} f_{|_{\rsf_{t}}} \,  d\mu_{i,j,t}   = \gamma_i(e) \gamma_j(e)  (M^k_{\pi,x})^{-1} (i,j)  \int_{\mathcal{I}_e} f(\lambda)\, d \theta_e (\lambda).
\end{equation*}
On the other hand, if $m>k$, then $y_m(t) / y_k(t)$ goes to zero as $t \in B^\ast$ converges to $\thy$ in $B^{\hyb}$. 
This means that
\[
\lim_{\substack{t \to \thy \\ t \in B^\ast}}\, \sup_{p \in U_{e,t}^2} \| c_t (p) \| = 0
\]
and in particular
\[
\lim_{\substack{t \to \thy \\ t \in B^\ast}}  \int_{U_{e,t}^2} f_{|_{\rsf_{t}}} \, d\mu_{i,j,t}   = 0.
\]
Since  $k \le l \le m$, we have $m>k$ provided that $l > k$. We infer that non-zero contributions of the integral over $U_e\cap \rsf_t$ to the limit~\eqref{eq:limit_in_question} only occur if $k=l=m$.

Finally, combining all the above considerations, we can compute the limit in \eqref{eq:limit_in_question} and arrive at Theorem~\ref{thm:cont2}. \qed

\subsection{Continuity III} \label{ss:ctIII}
Consider a pair of indices $(i,j)$ with $1\leq i\leq \graphgenus$ and $j > \graphgenus$.
\begin{thm}\label{thm:cont3} If $t \in B^\ast$ converges to $\thy=(0,x)$, then the measure $\mu_{i,j,t}$ converges to the zero measure on $\rsf^\hyb_\thy$. \end{thm}
The proof proceeds in the exact same way as in the previous section.

Taking the open sets $U_{e}$ around singular points $p^e$ of the fiber $\rsf_0= S_0$, one sees that outside the union of the $U_{e}$'s the measure $\omega_i \wedge \bar \omega_j$ extends continuously while $ \Im(\Omega_t)^{-1} (i,j)$ converges to zero (see Theorem \ref{thm:inverse_period_asymptotics}). Reasoning as in the preceding section, and supposing that $i \in J_\pi^k$, this remains true for the restriction of measures to the subsets $U_e$ for  $e\in \pi_m$ with $m < k$. 

 On the other hand, using the same arguments and notations as in the previous section, on the open set $U_e$ of an edge $e \in  \pi_m$, we have
\[ 
\mu_{i,j, t} = \frac{y_m(t)}{y_k(t)} \bar f_j(\underline z) \varphi(\underline z) \, d\tau^e_u \wedge d\Theta_e \qquad \text{on } U_e \cap \rsf_t,
\]
where $\varphi\colon U_e \to \C$ is a bounded function such that $\varphi(\underline z) = 0$ if $z^e_u=z^e_v = 0$. The latter property follows from the fact that $\omega_j$ is holomorphic at the singularities $p^e$, $e\in E$, of $S_0$. Proceeding as in the proof of Theorem~\ref{thm:cont2}, we arrive at Theorem~\ref{thm:cont3}.  \qed

\subsection{Proof of Theorem~\ref{thm:mainlocal_open}} For points $t\in B^*$, we decomposed the canonical measure as $\mu_t = \sum_{i,j=1}^g \mu_{i,j,t} $ in~\eqref{def:decomposition}. Using Remark~\ref{rem:first-reduction}, we then reduced the proof of Theorem~\ref{thm:mainlocal_open} to understanding the limits of the measures $\mu_{i,j,t}$ when $t \in B^\ast$ approaches a limit point of the form $\thy=(0,x)$ in $B^{\hyb}$. This was achieved in the three preceding subsections.  The link between the measures $\mu_{i,j,\thy}$ in Theorem~\ref{thm:cont2} and the canonical measure on $\rsf^\hyb_\thy$ is given by Theorem~\ref{thm:cmgraphscycles}. \qed

\subsection{Proof of the inverse lemma} \label{sec:inverse_lemma_proof}
We prove Lemma \ref{lem:inverse_lemma} proceeding by induction on $r$, the number of blocks in each row/column. 

If $r=1$, then $M(t) = \mA_{11}(t)$ and the claim is trivial. So suppose that $M(t)$ is an $(r,r)$ block matrix and we have already proven the claim for $(r- 1, r-1)$ block matrices. Then we can write
\[
M(t) = \left (  \begin{matrix} \Pi_{11}(t) &  \Pi_{12}(t) \\
\Pi_{21}(t) &  \Pi_{22}(t)
\end{matrix}  \right ),
\]
where $\Pi_{11} (t) = \mA_{11}(t) \in \C^{n_1 \times n_1}$, $\Pi_{22}(t)$ is the $(r-1, r-1)$ block matrix
\[
\Pi_{22}(t) = (\mA_{kl}(t))_{2 \le k, l \le r} \, \in \C^{n' \times n'}
\]
with $n' = n-n_1=\sum_{j=2}^r n_k$ and $\Pi_{21}(t) \in \C^{n' \times n_1}$, $\Pi_{12}(t) \in \C^{n_1 \times n'}$ are given by
\begin{align*}
\Pi_{21} &=  \left(  \begin{matrix} \mA_{21} \\  \vdots \\  \mA_{r1} \end{matrix}  \right ) =  \left  (  \begin{matrix} y_2 ({\widehat{ \mA}}_{21} +o(1))  \\  \vdots \\  y_r ({\widehat{\mA}}_{r1} +o(1)) \end{matrix}   \right ), \qquad \Pi_{12} =   \left (  \begin{matrix} \mA_{12} \\ \dots \\  \mA_{1r} \end{matrix}   \right  )^\transpose =   \left   (  \begin{matrix} y_2 ({\widehat{\mA}}_{12} +o(1)) \\ \dots \\   y_r ({\widehat{\mA}}_{1r} +o(1)) \end{matrix}   \right )^\transpose.
\end{align*}
Consider the Schur complements $S(t) \in \C^{n_1 \times n_1} $ and $\tilde S (t) \in \C^{n' \times n'} $ given by
\begin{align*}
	S(t) &\coloneqq  \Pi_{11}(t) -  \Pi_{12}(t) \Pi_{22}(t)^{-1} \Pi_{21}(t), \\
	\tilde{S} (t) &\coloneqq  \Pi_{22}(t) - \Pi_{21}(t) \Pi_{11}(t)^{-1} \Pi_{12}(t).
\end{align*}
By the Schur complement formula, the inverse of $M(t)$ can be written as
\begin{align}  \label{eq:SchurComplementFormula}
M(t)^{-1} &= \left (  \begin{matrix} \Psi_{11}(t) &  \Psi_{12}(t) \\
\Psi_{21}(t) &  \Psi_{22}(t) \end{matrix}  \right ) = \left( 
\begin{matrix}
S(t)^{-1} &   -S(t)^{-1} \Pi_{12}(t)\Pi_{22}(t)^{-1}\\
-\Pi_{22}(t)^{-1}\Pi_{21}(t)S(t)^{-1} & \tilde{S}(t)^{-1}
\end{matrix} \right),
\end{align}
provided that $\Pi_{11}(t)$, $\Pi_{22}(t)$, $S(t)$, and $\tilde{S}(t)$ are invertible. In the following, we prove that the invertibility condition holds and then use \eqref{eq:SchurComplementFormula} to compute the asymptotics of $M(t)^{-1}$.

The matrix $\Pi_{11} (t) = \mA_{11}(t)$ is invertible, since $\widehat{\mA}_{11}$ is invertible by assumption. 

The $(r-1, r-1)$ block matrix $\Pi_{22}(t)$ satisfies the assumptions of the induction hypothesis. In particular, $\Pi_{22}(t)$ is invertible and
\[
	\Pi_{22}^{-1} (t)=  \left ( \frac{1}{\ssub{y}!_{\min\{k,l\}} (t)} (\mathcal N_{kl} + o(1) ) \right )_{2 \le k,l \le r}
\]
for some matrices $\mathcal N_{kl} \in \C^{n_k \times n_l}$, $2 \le k,l \le r$. 

Turning to $S(t)$, note that, since $\lim_{t \to \thy} \frac{y_k (t)}{\ssub{y}!_{\min\{k,l\} }(t) }$ exists, the limit $\mathcal N \coloneqq  \lim_{t \to \thy} \Pi_{22}^{-1} \Pi_{21} \in \C^{n' \times n_1}$ exists as well. Indeed, we have $\mathcal N= (\mathcal N_i)_{i=2}^r$ where $\mathcal N_i \in \C^{n_i \times n_1} $, $i =2, \dots, r$, is given by
\begin{align*}
	 \mathcal N_i	= \lim_{t \to \thy} \sum_{j = 2}^r \frac{y_j(t) }{\ssub{y}!_{\min\{i,j\}} (t)} (\mathcal N_{ij} + o(1) ) ({\widehat{\mA}}_{j1} + o(1)) = \sum_{j = 2}^i \mathcal N_{ij} {\widehat{\mA}}_{j1}.  
\end{align*}
Combining this with the fact that $\Pi_{12} (t)= y_2(t) O(1)$, we obtain
\begin{align} \label{eq:SchurComplementAymptotics}
S(t) = \Pi_{11}(t) +  y_2(t) O(1) = y_1(t) \Bigl( {\widehat{\mA}}_{11} +  \frac{y_2(t)}{y_1(t)} O(1) \Bigr) = y_1(t) \Bigl( {\widehat{\mA}}_{11} + o(1) \Bigr).
\end{align}
In particular, $S(t)$ is invertible for $t$ close to $\thy$ by invertibility of ${\widehat{\mA}}_{11}$. 

It remains to treat $\tilde S (t) $. In order to do so, we consider $\tilde S (t)$ as a $(r-1, r-1)$ block matrix. Note that
\[
	\Pi_{21} \Pi_{11}^{-1} \Pi_{12} = \Big ( \mA_{k1} \mA_{11}^{-1} \mA_{1l} \Big )_{2 \le k,l \le r}
\]
and it follows from \eqref{eq:inverse_lemma_hyp} that for each index pair $(k,l)$ with $2 \le k,l \le r$, 
\begin{align*}
	\lim_{t \to \infty} \frac{1}{\ssub{y}!_{\max\{k,l\}} (t) } &\mA_{k1} (t) \mA_{11} (t) ^{-1} \mA_{1l} (t) = \\ &\lim_{t \to \infty} \frac{y_k (t) y_l (t)}{y_1 (t) \ssub{y}!_{\max\{k,l\}} (t) }( {\widehat{\mA}}_{k1}  {\widehat{\mA}}_{11}^{-1}  {\widehat{\mA}}_{1l} + o(1)) = 0.
\end{align*}
This means that the $(k,l)$-th block of $\tilde{S} (t) $ behaves like
\[
 \tilde{S}_{kl} (t)  = \ssub{y}!_{\max\{k,l\}}(t) ( {\widehat{\mA}}_{kl} + o(1)).
\]
We can thus apply the induction hypothesis to $\tilde{S}(t)$, and deduce in particular that it is invertible.

It remains to determine the asymptotic behavior of $M(t)^{-1}$ using~\eqref{eq:SchurComplementFormula}. From \eqref{eq:SchurComplementAymptotics}, we obtain that $\Psi_{11}(t) = y_1(t)^{-1}({{\widehat{\mA}}_{11}}^{-1} + o(1))$, as claimed. 
In addition, it is clear that
\[
\Psi_{21}(t) = - \frac{1}{y_1(t)} (N + o(1)) ({\widehat{\mA}}_{11}^{\, -1} + o(1)) = - \frac{1}{y_1(t)} (N {\widehat{\mA}}_{11}^{\,  -1} + o(1)).
\]
This implies the claimed asymptotics for $\Psi_{21}(t)$. The second off-diagonal matrix $\Psi_{12}(t)$ can be treated similarly. Finally, the claimed asymptotics for $\Psi_{22}(t)$ follow by again applying the induction hypothesis to the $(r-1, r-1)$ block matrix $\tilde{S}(t)$. The proof is complete. \qed


\section{Proof of the main theorem}\label{sec:proof}

In this section, we present the proof of our main theorem. 
Let $S_0$ be a stable curve of genus $g$ with dual graph $G = (V, E, \genusfunction)$. As in the previous sections, let $\rsf \to B$ be the versal analytic family of stable curves over a polydisc $B$ of dimension $\ndim=3g-3$ with $\rsf_0 = S_0$. Consider the hybrid space $B^\hyb$ and the family of hybrid curves $\rsf^\hyb \to B^{\hyb}$. We equip each fiber $\rsf^\hyb_\thy$, $\thy \in B^\hyb$, with its canonical measure $\mu_\thy$.%

\begin{thm}[Continuity: local case] \label{thm:mainlocal_intro}  Notations as above, the family of canonically measured hybrid curves $(\rsf^\hyb_\thy,\mu_\thy)_{\thy\in B^\hyb}$ is continuous.
\end{thm}

By Remark~\ref{rem:first-reduction}, it suffices to prove the continuity at hybrid points $\thy_0 \in B^\hyb$ of the form $\thy_0 = (0,x)$, that is, their underlying complex point equals $t = 0 \in B$. Fix such a hybrid point $\thy_0$ and let $\pi = (\pi_1, \dots, \pi_r)$ be the ordered partition of the edge set $E$ with  $\thy_0 \in D_\pi^\hyb$. By Theorem~\ref{thm:mainlocal_open}, the canonical measures are continuous at $\thy_0$ through the open part $B^*$ of $B^\hyb$. In order to establish Theorem~\ref{thm:mainlocal_intro}, we need to prove the continuity at $\thy_0$ through all hybrid strata $D_{\pi'}^\hyb$ of ordered partitions $\pi'$ with $\pi' \preceq \pi$.

More precisely, fix a subset of edges $F \subseteq E$ and an ordered partition $\pi'$ of $F$ with $\pi' \preceq \pi$. We need to show that the family of canonically measured hybrid curves $\rsf^\hyb \rest{D_{\pi'}^\hyb \cup \{\thy_0\}}$ is continuous at $\thy_0$ for the induced topology on $D_{\pi'}^\hyb \cup \{\thy_0\} \subseteq B^\hyb$.  We will reduce this claim to Theorem~\ref{thm:mainlocal_open} and Theorem~\ref{thm:mainmetricgraphs}.

 Let $G' = (V', F)$ be the stable dual graph of the stable curve $\rsf_t$ for any $t$ in the stratum $D_{\pi'}^\hyb$. The graph $G'$ is obtained from $G$ by contracting all edges in $E \setminus F$. Since $\pi' \preceq \pi$, there exists $r_0\leq r$ such that the initial part $\pi_I \coloneqq  (\pi_1, \dots, \pi_{r_0})$ of $\pi$ forms an ordered partition of $F$.

 We have a continuous forgetful map $\forget \colon D_{\pi'}^\hyb \sqcup D_\pi^\hyb \to \ssub{\widetilde \mg}!_{G'}^{\,\,\,\trop}$, which sends any hybrid point $\thy = (t,x)$ to the point $\forget(\thy) \coloneqq  x\rest{F}$ either in $\inn \sigma_{\pi'}$ or in $\inn \sigma_{\pi_I}$, depending on whether $\thy \in D_{\pi'}^\hyb$ or $\thy \in D_{\pi}^\hyb$. By Proposition~\ref{prop:PushoutMeasures}, for any point $\thy \in D_{\pi'}^\hyb \sqcup D_\pi^\hyb$ and any edge $e\in F$, the restriction of the canonical measure $\mu_\thy$ to the interval $\Ical_{e}\rest{\thy}$ of the hybrid curve $\rsf^\hyb_\thy$ coincides with the restriction of the canonical measure $\mu^\can_{\forget(\thy)}$ on the tropical curve 
   $\widetilde {\mathscr C}_{\forget(\thy)}^{\, \trop}$ to the interval $\Ical_{e}\rest{\forget(\thy)}$ of $\widetilde {\mathscr C}_{\forget(\thy)}^{\, \trop}$. 
By Theorem~\ref{thm:mainmetricgraphs}, the universal family of canonically measured tropical curves $\ssub{\widetilde {\mathscr C}}!_{G'}^{\, \trop}$ is continuous over $\ssub{\widetilde \mg}!_{G'}^{\,\,\,\trop}$. It follows that for any edge $e\in F$, the family of measured intervals $(\Ical_{\thy,e},\mu_{\thy} \rest{{\Ical_{\thy,e}}})$ is continuous over $D_{\pi'}^\hyb \sqcup D_\pi^\hyb$. This proves half of the continuity claim, namely, the continuity of canonical measures on the common intervals of the respective hybrid curves.

To prove the continuity on the Riemann surface parts, let $\proj_{\pi\succ \pi'}\colon G   \to G' $ be the contraction map, which contracts all edges in  $E \setminus F$. Denote by $V_{v}$, $v \in V'$, the set of all vertices of $G$ that are mapped to $v$ by $\proj_{\pi\succ \pi'}$. For a vertex $v\in V'$, consider the subcurve $S_{v} \subseteq S_0 = \rsf_{0} $ consisting of the union of all irreducible components $C_{0,u} = C_u$ of $S_0$,  for $u \in V_{v}$.  The subcurve $S_{v}$ comes with a marking given by those nodes of $S_0$ that lie on $S_v$ and correspond to edges in $F$. The marked curve $S_{v}$ is stable and we denote by $\prescript{}{v}B$ its analytic versal deformation space. Let $\prescript{}{v}\rsf^\hyb \to \prescript{}{v}B^\hyb$ be the family of hybrid curves over the associated hybrid space $\prescript{}{v}B^\hyb$, endowed fiberwise with the canonical measure.

We have a natural projection map from $D_F = \bigsqcup_{F \subseteq \hat F \subseteq E} \inn D_{\hat F}$ to $B_v$, which gives rise to a hybrid projection map $\prescript{}{v}\pr \colon B^\hyb_{[\pi', \pi]} \to \prescript{}{v}B^\hyb$. Recall that the subspace $B^\hyb_{[\pi', \pi]} \subseteq B^\hyb$ was introduced in \eqref{eq:BInterval} and corresponds to ordered partitions $\varrho$ with $\pi' \preceq \varrho \preceq \pi$. Under this last map, the inclusion
\[ D_{\pi'}^\hyb \subseteq  \prescript{}{v}\pr^{-1}\big (\prescript{}{v}B^*\big ) \]
holds true. In addition, the pullback $\prescript{}{v}\pr^*\Bigl(\prescript{}{v}\rsf^\hyb\Bigr)$ can be identified naturally with a measured subspace of the hybrid family $\rsf^\hyb$.  By Theorem~\ref{thm:mainlocal_open} (that extends verbatim to the marked setting), the family of canonically measured hybrid curves $\prescript{}{v}\rsf^\hyb_\thy$, $\thy \in \prescript{}{v}B^\hyb$, is continuous at the point $\prescript{}{v}\pr(\thy_0)$ through the open part $\prescript{}{v}B^*$.
Hence the family of measured spaces $\prescript{}{v}\pr^*\Bigl(\prescript{}{v}\rsf^\hyb\Bigr)$ is continuous over the subspace $D_{\pi'}^\hyb \cup\{\thy_0\} \subset B^\hyb$. Since this holds for all vertices $v \in V' $, we get the second half of the continuity, and Theorem~\ref{thm:mainlocal_intro} follows.
\qed

\appendix

\section{The topology on the versal hybrid family}\label{sec:AppendixTopology}
Let $S_0$ be a stable (marked) Riemann surface and $\rsf \to B$ its versal deformation family. Consider the family of hybrid curves $\pr \colon \rsf^\hyb \to B^\hyb$ over the hybrid deformation space $B^\hyb$. In what follows, we make the topology on $\rsf^\hyb$ precise. We use the terminology and notations from Section~\ref{sec:versal_hybrid_curve}.

The following construction is needed later on. For $e \in E$, we obtain a normalized family
\begin{equation} \label{eq:normalized_e}
 \breve{\rsf}_e \to D_e,
\end{equation}
by resolving the singular points $p^e(t)$, $t \in B$, in the versal family $\rsf\rest{D_e} \to D_e$. That is, in each fiber $\rsf_t$, $t \in D_e$, we replace the node $p^e(t)$ of the edge $e = uv$ in the dual graph $G_t$ by two different points $p_u^e(t)$ and  $p_v^e(t)$ lying on the components $C_{u,t}$ and $C_{v,t}$, thereby making them disjoint. We stress that $\breve{\rsf}_e$ is defined over $D_e$.

We will define a topology on $\rsf^\hyb$ by specifying for each point $\p \in \rsf^{\hyb}$ a neighborhood subbase $\mathcal{V}_0(\p)$. The total neighborhood system $\mathcal{V}(\p)$ at $\p$ is then obtained by taking first all intersections of finitely many sets in $\mathcal{V}_0(\p)$ and then all supersets of such sets.

Fix a point $\p \in S^{\hyb}$. Assume that $\p \in S_\pi^{\hyb}$ for the ordered partition $\pi$ of the subset $F \subseteq E$ (see~\eqref{eq:S_hyb_decomposition}; here, $S_{\pi_\varnothing}^{\hyb} \coloneqq  \rsf^\ast $ for the empty partition $\pi_\varnothing$ of $F = \varnothing$). Let $\thy = (t,x) \coloneqq  \pr(\p) \in B^\hyb$ be the hybrid base point of the fiber $\rsf^\hyb_\thy$ containing $\p$. The neighborhood subbase $\mathcal{V}_0(\p)$ will contain two different types of sets.  We first add the preimage $\pr^{-1}(U)$ of every neighborhood $U$ of the base point $\thy$ in $B^\hyb$ to $\mathcal{V}_0(\p)$. This ensures the continuity of the hybrid projection map $\pr \colon \rsf^\hyb \to D^\hyb$ at $\p$.

The second type of neighborhoods in $\mathcal{V}_0(\p)$ clarifies the topological relation between the subfamilies in \eqref{eq:S_hyb_decomposition}. These sets have the structure
\begin{equation} \label{eq:V}
	W = \bigsqcup_{\pi' \preceq \pi} W_{\pi'}
\end{equation}
where each $W_{\pi'}$ is a subset of $\rsf_{\pi'}^\hyb$. Their precise form depends on the type of the point $\p$ according to the decomposition \eqref{eq:point_types}. We will proceed by case distinction. Assume that...

$(i)$ $\p$ is a smooth point on some Riemann surface component of $\rsf^\hyb_\thy$, that is, $\p$ belongs to $\Pty_\pi^0$. In this case, the set $W_{\pi'}$ will simply reflect the topology on the original family of Riemann surfaces $\rsf$. By \eqref{eq:P1}, we can formally write $\p = (p,x)$ for some point $p \in P^\circ$ (see \eqref{eq:P^0}).

Suppose $W_p$ is a neighborhood of $p$ in $\rsf$ containing no nodes, that is, $W_p \subseteq P^\circ$. Then for any ordered partition $\pi'$ with $\pi' \preceq \pi$, we set
\begin{equation} \label{eq:clear_definition}
	W_{\pi'} \coloneqq   ( \bigcup_{s \in \inn D_{F'}} W_p \cap \rsf_s  ) \times \inn \sigma_{\pi'} =   (W_p \cap \rsf_{F'}) \times \inn \sigma_{\pi'} \subseteq \rsf_{\pi'}^\hyb
\end{equation}
where $F' = E_{\pi'}$. By \eqref{eq:smooth_points_normalized}, we can view $W_{\pi'}$ as a subset of $\rsf_{\pi'}^\hyb$. Using the sets $W_{\pi'}$, we define a neighborhood $W$ of $\p$ in $\rsf^\hyb$ by \eqref{eq:V}. It contains only smooth points of the corresponding hybrid curves.

$(ii)$ $\p$ is an interior point of an interval $\mathcal{I}_e$ in $\rsf^\hyb_\thy$, that is, $\p$ belongs to $\Pty_\pi^2$. Formally, we can write the point $\p$ as $\p = (t, x, \lambda_u^e, \lambda_v^e)$ for some edge $e \in F$ and $0 < \lambda_u^e, \lambda_v^e < x_e$ with $\lambda_u^e + \lambda_v^e = x_e$ (see \eqref{eq:P3}).   Fix an ordered partition $\pi'$ of some subset $F'$ with $\pi' \preceq \pi$. Suppose that...

\begin{itemize}[leftmargin = 2em]
\item $\pi'$ contains the edge $e$ (that is, $e \in F'$). Then, we can naturally relate the respective intervals in the hybrid curve fibers. Namely, we can view
\[
	W_{\pi'} \coloneqq  \left \{ \bigl(s, y, \tilde{\lambda}_u^e, \tilde{\lambda}_v^e\bigr) \in \inn D_{F'} \times \Ical_{\pi'}^e \,\, \st \,\,\, |\tilde{\lambda}_u^e / y_e -  \lambda_u^e / x_e | < \varepsilon \right \}
\]
as a subset of $\rsf^\hyb_{\pi'}$ by \eqref{eq:HybridNormalizedFamily}. It is contained in the intervals representing the edge $e$ in the respective hybrid curve fibers, that is, $W_{\pi'} \subseteq \Pty_{\pi'}^2$. 

\item $\pi'$ does not contain the edge $e$ (that is, $e \notin F'$).  In this case, we use the $\Log$-maps introduced in \eqref{eq:DefLogMap}. Fix a standard coordinate neighborhood $(U,z)$ of the point $p^e(t)$ in the original family $\rsf \to B$. Moreover, take a small neighborhood $W_p$ of $p^e(t)$ in $\rsf$ and some $\varepsilon > 0$. Using the $\Log$-map on $U$, we define
\[
	W_{\pi'} \coloneqq  \left \{ q \in  \rsf_{F'} \cap U\,\,\st \,\, |\Log_e(q) - \lambda_u^e / x_e| < \varepsilon \text{ and } q \text{ belongs to } W_p \right  \} \times \inn \sigma_{\pi'}.
\]
That is, we take those points in fibers over $D_{\pi'}^\hyb$ that lie in the Riemann surfaces parts, are close to the node $p^e(t)$ in the original family $\rsf$, and whose logarithmic images are close to the point in question on the interval. By definition, we have $W_{\pi'}  \subseteq \Pty_{\pi'}^0$ (see \eqref{eq:P1}).

\end{itemize}
Altogether, we obtain the neighborhood $W$ for $\p$ from the sets $W_{\pi'}$ by \eqref{eq:V}.

$(iii)$ $\p$ is the attachment point of an interval $\mathcal{I}_e$ in $\rsf^\hyb_\thy$, that is, $\p$ belongs to $\Pty_\pi^1$ (see \eqref{eq:P2}). Parametrizing the interval $\mathcal{I}_e$ starting from $\p$, we can assume that $\p = \p_u^e(\thy)$ is the left endpoint $0$ of the interval $\Ical_e = [0, x_e]$ representing the edge $e = uv \in F$ in the hybrid curve $\rsf^\hyb_\thy$. Fix an ordered partition $\pi'$ of some subset $F'$ with $\pi' \preceq \pi$. Suppose that... 

\begin{itemize}[leftmargin = 2em]
\item $\pi'$ does not contain the edge $e$. Then, we use the $\Log$-maps in the same way as in case (ii). Namely, we fix a standard coordinate neighborhood $(U,z)$ of the point $p^e(t)$ in the original family $\rsf$, a small neighborhood $W_p$ of $p^e(t)$ in $\rsf$ and some $\varepsilon > 0$. Using the $\Log$-map on $U$, we then set
\[
	W_{\pi'} \coloneqq  \left \{ q \in \rsf_{F'} \cap U \, \st \,\Log_e(q) < \varepsilon \text{ and } q \text{ belongs to } W_p \right\} \times \inn \sigma_{\pi'}.
\]
That is, we take those points in fibers over $D_{\pi'}^\hyb$ that lie in the Riemann surfaces parts, are close to the node $p^e(t)$ in the original family $\rsf$, and have logarithmic image close to the left endpoint of the interval.  Again, $W_{\pi'}$ is viewed as a subset of $\rsf_{\pi'}^\hyb$ by \eqref{eq:P1}. It contains only smooth points of the corresponding fibers, that is, $W_{\pi'} \subseteq \Pty_{\pi'}^0$.

\item $\pi'$ contains the edge $e$ (that is, $e \in F'$).  Then, the singularity along the section $p^e$ was resolved also in the family $\rsf_{\pi'}^\hyb$. Thus, locally, $\rsf_{\pi'}^\hyb$ looks like the normalized family $\breve{\rsf}_e$ from \eqref{eq:normalized_e} with some intervals attached. In this case, we want the set $W_{\pi'}$ to contain both parts of stable Riemann surfaces and intervals. 
To define the first part, note that we can view
\[
\breve{\rsf}_e \rest{\inn D_{F'}} \times \inn \sigma_{\pi'}
\]
 as a subset of $\rsf_{\pi'}^\hyb$ (similar to \eqref{eq:clear_definition}). In particular, for a neighborhood $\breve{W}$ of $p^e_u(t)$ in $\breve{\rsf}_e$, we can see
\[
	W_{\pi', 1} \coloneqq  \left(\bigcup_{s \in \inn D_{F'}} \breve{W} \cap (\breve{\rsf}_e\rest{s})\right)  \times \inn \sigma_{\pi'} =\left( \breve{W}  \cap (\breve{\rsf}_e\rest{\inn D_{F'}})\right)    \times \inn \sigma_{\pi'}
\]
as a subset of $\rsf_{\pi'}^\hyb$. By \eqref{eq:HybridNormalizedFamily}, the second part of $W_{\pi'}$ is the subset 
\[
	W_{\pi', 2} \coloneqq  \left \{ (s, y, \lambda_u^e, \lambda_v^e) \in \inn D_{F'} \times \Ical_{\pi'}^e \, \st \, \lambda_u^e / y_e < \varepsilon \right \} \subseteq \rsf_{\pi'}^\hyb,
\]
which corresponds to taking interval parts $[0, \varepsilon]$ in the edge $e$ of the hybrid curves. Altogether, the subset $W_{\pi'} \subseteq \rsf_{\pi'}^\hyb$ is defined as the union $W_{\pi'} = W_{\pi', 1} \cup W_{\pi', 2}$. It contains points of all three types in \eqref{eq:point_types}.

\end{itemize}
Again, the neighborhood $W$ of $\p$ is obtained from the sets $W_{\pi'}$ by \eqref{eq:V}.

\smallskip

We can now finish the definition of the neighborhood subbase $\mathcal{V}_0(\p)$ of the point $\p \in \rsf\hyb$. Distinguishing the three cases $(i)$-$(ii)$-$(iii)$, we add to $\mathcal{V}_0(\p)$ all sets $W$ that can be constructed in the way described above (ranging over all $\varepsilon > 0$, and neighborhoods $U$, $W_p$, $\breve{W}$, etc.).

Finally, the topology on $\rsf^{\hyb}$ is the unique one such that $\mathcal{V}_0(\p)$ is a neighborhood subbase for each $\p \in \rsf^\hyb$. A direct verification shows that the related neighborhood system $(\mathcal{V}(\p))_{\p}$, obtained by taking finite intersections and supersets of sets in $\mathcal{V}_0(\p)$, satisfies the axioms of a neighborhood system on $\rsf^\hyb$, concluding the definition of the hybrid topology on $\rsf^\hyb$. 

\bibliographystyle{alpha}
\bibliography{bibliography}

\appendix

\end{document}